\documentclass[11pt]{article}
\usepackage[a4paper,margin=1in]{geometry}
\usepackage{amsmath,amssymb,amsfonts,amsthm,mathtools}
\usepackage{graphicx}
\usepackage{booktabs,array}
\usepackage{caption,subcaption}
\usepackage{xcolor,needspace,etoolbox,float}
\usepackage[section]{placeins}

\pretocmd{\subsection}{\Needspace{5\baselineskip}}{}{}
\usepackage[backend=biber,style=numeric]{biblatex}
\usepackage[colorlinks=true,linkcolor=blue!50!black,citecolor=blue!50!black,urlcolor=blue!50!black]{hyperref}

\numberwithin{equation}{section}
\theoremstyle{definition}
\newtheorem{theorem}{Theorem}[section]
\newtheorem{proposition}{Proposition}[section]
\newtheorem{lemma}{Lemma}[section]

\newtheorem{remark}{Remark}[section]
\newtheorem{example}{Example}[section]
\newtheorem{assumption}{Assumption}[section]

\newcommand{\ZZ}{\mathbb{Z}}
\newcommand{\RR}{\mathbb{R}}

\newcommand{\TT}{\mathbb{T}}
\newcommand{\Aa}{\mathcal{A}}
\newcommand{\Bb}{\mathcal{B}}
\newcommand{\Cc}{\mathcal{C}}
\newcommand{\Ee}{\mathcal{E}}
\newcommand{\Ff}{\mathcal{F}}
\newcommand{\Gg}{\mathcal{G}}
\newcommand{\Ii}{\mathcal{I}}
\newcommand{\Ll}{\mathcal{L}}
\newcommand{\Ss}{\mathcal{S}}
\newcommand{\Tt}{\mathcal{T}}
\newcommand{\Oo}{\mathcal{O}}
\renewcommand{\a}{\alpha}
\newcommand{\s}{\sigma}
\renewcommand{\t}{\theta}
\newcommand{\x}{\chi}
\newcommand{\eps}{\epsilon}
\newcommand{\z}{\zeta}
\renewcommand{\d}{\mathrm{d}}

\DeclareMathOperator{\supp}{supp}

\DeclareMathOperator{\rank}{rank}
\newcommand{\A}{\tilde{A}}
\newcommand{\bb}{\tilde{b}}

\title{Data-driven discovery of sparse governing equations with parameterized operators}
\author{Mike Wendels$^*$, J. Nathan Kutz$^{*,\dag}$\\[.2in]
$^*$ Department of Applied Mathematics, University of Washington, Seattle, USA\\
$^\dag$ Autodesk Research, London, UK}
\date{\today}

\begin{document}
\maketitle

\begin{abstract}
We study parameterized SINDy (pSINDy), an extension of the Sparse Identification of Nonlinear Dynamics (SINDy) framework for the discovery of parsimonious governing equations involving operators with unknown parameters such as fractional-order Caputo derivatives, fractional Laplacians, and parameterized kernel operators. The framework jointly identifies the unknown operator parameters and sparse library coefficients through a sparsity-promoting optimization problem. To solve this problem, we develop and analyze the parameterized sequentially thresholded least-squares (pSTLSQ) algorithm. In each iteration, it removes terms with small coefficients, updates the operator parameters by a fixed number of optimization steps that improve the fit of the remaining library terms, and refits their coefficients by least squares. Under suitable assumptions, we prove that pSTLSQ converges locally to a local minimizer and give sufficient conditions for identifying the terms present in the governing equation. Numerical experiments demonstrate recovery of governing equations from noisy synthetic data for fractional ODEs, fractional PDEs, and nonlocal diffusion equations. We also examine the influence of the hyperparameters, discretization scheme, and noise, and compare pSTLSQ with two benchmark methods. All three methods identify the correct terms at low noise levels, but differ in parameter accuracy, cost, and sensitivity to initialization and the candidate library.
\end{abstract}

\section{Introduction}\label{sec:intro}

Inferring governing equations from observational data is a central problem in applied mathematics and data science, particularly when accurate models are difficult to derive from first principles. Since multiple candidate models may provide comparable descriptions of finite and noisy observations, structural assumptions are essential to avoid overfitting. A commonly adopted constraint is parsimony, as for many physical systems the governing equations involve a relatively small number of terms when expressed in appropriate coordinates \cite{kutz2022parsimony}.

The Sparse Identification of Nonlinear Dynamics (SINDy) framework introduced in \cite{brunton2016discovering} builds on this principle by identifying governing equations from data as a sparse linear combination of prescribed library terms with constant coefficients. For a scalar equation, this takes the form
\begin{align}\label{eq:nonfrac}
    \Bb(u) = \sum_{i = 1}^p \xi^i \Aa^i(u),
\end{align}
where $\Bb$ is a known operator, $\Aa^i$, $i = 1, \dots, p$, are prescribed candidate operators, and $\xi \in \RR^p$ is an unknown coefficient vector assumed to be sparse. For ordinary differential equations (ODEs), the library consists of candidate functions of the state, and for partial differential equations (PDEs) it additionally includes partial derivatives and potentially nonlinear terms \cite{rudy2017data}. For systems, this representation is applied componentwise, with library terms that may depend on all state components. Discretizing and evaluating the library and $\Bb$ at $u$ gives a matrix $A$ and a vector $b$, respectively, from which the library coefficients $\xi$ for each component equation are then approximated through the sparse regression problem
\begin{align}\label{eq:sindy0}
    \min_{\xi \in \RR^p}\|A\xi - b\|_2^2 + \lambda^2\|\xi\|_0,
\end{align}
where $\|\xi\|_0 := |\supp(\xi)|$ counts the nonzero entries in $\xi$, and $\lambda > 0$ balances the misfit to the observed data with sparsity. Problem \eqref{eq:sindy0} can be approached through the sequentially thresholded least-squares (STLSQ) algorithm, which alternates thresholding a candidate library coefficient vector with least-squares regression. If $A$ has full column rank, STLSQ finds a local minimizer for optimization problem \eqref{eq:sindy0} in finitely many steps \cite{zhang2019convergence}.

Several adaptations and extensions of the SINDy framework have been proposed. Sparse Relaxed Regularized Regression (SR3) \cite{zheng2018unified} provides an alternative formulation of sparse regression problem \eqref{eq:sindy0}. Weak formulations have been developed to avoid numerical differentiation of the observed data \cite{messenger2021weak}. The framework has also been extended to controlled systems \cite{kaiser2018sparse}, implicit equations \cite{kaheman2020sindy}, and stochastic differential equations \cite{boninsegna2018sparse}, as well as PDEs with spatially or temporally varying coefficients \cite{rudy2019data} and time-varying dynamical systems \cite{li2019discovering}. Furthermore, statistical criteria have been used to compare and select among candidate SINDy models \cite{mangan2017model}.

Governing equations may involve unknown parameters arising from scaling and nondimensionalization, such as the Reynolds number \cite{buckingham1914physically, segel1972simplification, bakarji2022dimensionally}. When these parameters enter as coefficients of the left-hand side operator or the library terms, they can be identified by SINDy. However, in the case where the unknown parameters cannot be absorbed into the coefficients, joint estimation of the parameters and sparse library coefficients is required. Several extensions of SINDy address this problem. For example, \cite{champion2020unified} combines sparse coefficient estimation with Newton updates of nonlinear library parameters, while ADAM-SINDy \cite{viknesh2026adam} uses the Adam optimization method \cite{kingma2015adam} to jointly estimate these parameters and the library coefficients, incorporating sparsity regularization and thresholding. In both these formulations, the unknown parameters enter the library, while the left-hand side operator $\Bb$ is known and independent of these parameters. In this work, we instead consider equations of the form
\begin{align}\label{eq:eq_to_be_discovered}
    \Bb_{\a}(u) = f(u),
\end{align}
where $\Bb_{\a}$ is a parameterized operator with unknown parameter $\a \in \RR^m$, and $f$ is an unknown function which we approximate by a sparse linear combination of candidate library terms. Equations of this form arise, for example, in fractional differential equations and nonlocal diffusion equations with parameterized kernels. We remark that the parameterized operator $\Bb_{\a}$ is prescribed as part of the equation, so its presence does not depend on the selection of library terms in contrast to the formulations of \cite{champion2020unified, viknesh2026adam}.

Several approaches have been proposed to discover equations of the form \eqref{eq:eq_to_be_discovered}. These differ in how they estimate the library coefficients and operator parameters, and how they select the library terms in a sparsity-promoting way. The method of \cite{zhang2024sparse} uses a constrained $\ell^0$ formulation and applies regression followed by thresholding small coefficients over a candidate parameter grid. It selects among the candidate models by balancing an error against the number of nonzero coefficients. The error obtained by substituting the data into the equation (model error) and the error between simulated trajectories and the data (trajectory error) are compared. The approaches in \cite{singh2022data,vats2024new} jointly estimate the coefficients and operator parameters by minimizing the model error with $\ell^1$ regularization of the coefficients and an additional penalty on the operator parameter. The method of \cite{fan2023data} alternates coefficient regression and thresholding with updates of the operator parameter, where the coefficients are fitted by minimizing the model error with the parameter fixed while the parameter is updated by minimizing the trajectory error with the thresholded coefficients fixed. The method of \cite{xiangnan2025data} first reconstructs the data with a neural network to mitigate the effects of noise, after which it alternates coefficient regression and thresholding with derivative-free parameter optimization using the model error obtained by substituting the reconstructed solution into the equation. The method of \cite{thanasutives2026robust} instead reduces sensitivity to measurement noise by fitting the equations in weak form and searches over candidate terms and parameter values, estimating the coefficients by regression for each candidate model. Sparsity is then imposed by selecting the number of terms through a criterion balancing the model error on validation data and model size. We note that the formulations in \cite{vats2024new,xiangnan2025data,thanasutives2026robust} also allow the library terms to depend on unknown parameters.

In this work, we follow an alternating approach similar to \cite{fan2023data,xiangnan2025data} to estimate the operator parameter $\a$ and sparse library coefficients.
In particular, writing $f(u) = \sum_{i = 1}^p\xi^i\Aa^i(u)$, we consider the following sparsity-promoting optimization problem
\begin{align}\label{eq:psindy0}
    \min_{\a \in \Omega_{\a}, \xi \in \RR^p}\|A\xi - b_{\a}\|_2^2 + \lambda^2\|\xi\|_0,
\end{align}
where $b_{\a}$ is the discretization of $\Bb_{\a}(u)$, and $\Omega_{\a} \subseteq \RR^m$ is the set to which $\a$ is constrained, e.g., to restrict fractional derivative orders to a prescribed range. We refer to this formulation as parameterized SINDy (pSINDy). To solve the pSINDy optimization problem, we introduce the parameterized STLSQ (pSTLSQ) algorithm, which alternates sparse regression as in STLSQ with optimization over the operator parameter $\a$. To this end, we assume that $b_\a$ is continuously differentiable in $\a$ and use its derivatives to compute the parameter updates.

We analyze how the parameter updates and coefficient thresholding in pSTLSQ affect convergence and recovery of the governing terms. As a starting point, we consider exact pSTLSQ, in which the parameter update is done by exactly solving the corresponding optimization problem. Under suitable assumptions, we show that every fixed point of exact pSTLSQ is a local minimizer and establish finite termination when the parameter optimization problems have unique solutions. Afterwards, we study inexact pSTLSQ, in which each parameter update consists of a fixed finite number of steps of an iterative optimization method. For this version of pSTLSQ, we establish local convergence to a local minimizer provided the coefficients are strictly separated from the threshold and the iterative method contracts around the local minimizer. Additionally, we discuss conditions for the threshold $\lambda$ and the initial guess of $\a$ in pSTLSQ that allow for the recovery of the true model, and we show that nondegenerate interior local minimizers of the pSINDy optimization problem \eqref{eq:psindy0} are stable under small perturbations of the data. Numerical experiments on fractional ODEs, fractional PDEs, and nonlocal diffusion equations with parameterized kernels examine the influence of noise, discretization, and the pSTLSQ hyperparameters. Finally, we compare pSTLSQ with a grid search over the parameter and with a decoupled method that estimates the parameter before any library term is removed. Both fit the general form of the pSTLSQ iteration introduced in this work but update the parameter at different stages. All three recover the true model at low noise levels, but differ in parameter accuracy, cost, and sensitivity to initialization and the chosen library.

The outline of this work is as follows. Section \ref{sec:preliminaries} introduces the sparse regression formulation of SINDy and reviews the properties of STLSQ. Section \ref{sec:family} introduces and analyzes the proposed pSTLSQ algorithm. Section \ref{sec:numerics} presents the numerical results for fractional ODEs, fractional PDEs, and nonlocal diffusion equations. Section \ref{sec:sensitivity} examines their sensitivity to the hyperparameters, discretization, and noise. Section \ref{sec:comparison} discusses the performance of pSTLSQ against benchmark methods. Section \ref{sec:conclusion} concludes.

\textbf{Notation and conventions.} Throughout this work, vector and matrix norms without a subscript respectively denote Euclidean and spectral norms. For $S \subseteq \{1, \dots, p\}$, we write $\overline S := \{1, \dots, p\} \setminus S$, let $\xi|_S$ denote the subvector with entries indexed by $S$, and let $A_S$ denote the submatrix containing the corresponding columns. Maxima over empty sets of coefficient magnitudes are taken as zero, and minima as $+\infty$. We write $B_\rho(v)$ for the open ball of radius $\rho$ centered at $v$, and $A^\dagger$ for the Moore-Penrose pseudoinverse. Furthermore, $\Cc^r(D)$ with $D$ a closed set means that the function is $\Cc^r$ in some open neighborhood of $D$.

\section{Preliminaries} \label{sec:preliminaries}

We begin by reviewing the SINDy framework and the STLSQ algorithm, which serve as the basis for the pSINDy framework developed in the following section. Throughout this discussion, we use the index conventions in Table \ref{tab:indices}. We also use $k$ to denote the iteration counter of STLSQ -- and for pSTLSQ in Section \ref{sec:family} -- and the meaning will be clear from context.

\begin{table}[H]
    \centering
    \begin{tabular}{ll}
        \toprule
        $i = 1, \dots, p$ & library index \\
        $j = 1, \dots, n_\x$ & initial condition index \\
        $k = 1, \dots, n_t$ & time index in the ODE setting \\
        $k_r = 1, \dots, n_{x^r}$, $r = 0, \dots, d$ & space-time indices in the PDE setting, $x^0 = t$ \\
        $\ell = 1, \dots, n$ & state component index \\
        $q = 1, \dots, N$ & flattened observation index \\
        \bottomrule
    \end{tabular}
    \caption{Index conventions used in the discussion of SINDy and pSINDy.}
    \label{tab:indices}
\end{table}

\subsection{Sparse identification of ordinary and partial differential equations}

\textbf{Scalar ODE discovery.} We begin by considering SINDy for learning a scalar ODE. For the sake of presentation, we first consider learning an autonomous first-order ODE
\begin{align}\label{eq:ODE}
    \dot{u}(t) = f(u(t)), \qquad u(0) = \x,
\end{align}
for a scalar-valued and sufficiently smooth function $u : [0, T] \to \RR, \, T > 0$, where we approximate $f$ as a linear combination of a priori chosen library functions $\t^i : \RR \to \RR$, $i = 1, \dots, p$, i.e.,
\begin{align}\label{eq:f}
    f(u) = \sum_{i = 1}^p\xi^i\t^i(u).
\end{align}
Given data $\{u_{\x_j}(t^k)\}$ on a time grid $t^k, \, k = 1, 2, \dots, n_t$ for possibly different initial conditions $\x_j \in \RR, \, j = 1, 2, \dots, n_\x$, we collect the data in the tensors $\A$ and $\bb$ as
\begin{align}\label{eq:sindy_Ab}
    \A_{i}^{j, k} = \t^i(u_{\x_j}(t^k)), \quad \bb^{j, k} = \Tt_{t^k}(u_{\x_j}),
\end{align}
where $\Tt_{t^k}$ denotes a numerical approximation of the time derivative at time $t^k$. SINDy then balances minimizing the misfit to the data and promoting sparsity by finding library coefficients $\xi$ that solve
\begin{align*}
    \min_{\xi \in \RR^p} \sum_{j = 1}^{n_\x} \sum_{k = 1}^{n_t} |{\textstyle\sum_{i = 1}^p}\A_i^{j, k}\xi^i - \bb^{j, k}|^2 + \lambda^2\|\xi\|_0,
\end{align*}
where $\lambda > 0$ is the thresholding parameter. Finally, ODEs with a different left-hand side than in equation \eqref{eq:ODE} can be learned by modifying $\Tt_{t^k}$ accordingly.

\textbf{Scalar PDE discovery.} The extension to scalar PDEs follows the same principles as in the above ODE case. We illustrate it for the PDE
\begin{align}\label{eq:PDE}
    u_t = \Ff(u), \quad u(0, x) = u_0(x),
\end{align}
for a scalar-valued and sufficiently smooth function $u : [0, T] \times \Omega \rightarrow \RR, \, T > 0$, where $\Omega \subseteq \RR^d$ is the spatial domain. We approximate the partial differential operator $\Ff$ as a linear combination of a priori chosen candidate differential operators,
\begin{align}\label{eq:Ff}
    \Ff(u) = \sum_{i = 1}^{p} \xi^i \Aa^i(u).
\end{align}

Let $x_r^{k_r}$, $k_r = 1, 2, \dots, n_{x^r}$, denote the grid for the coordinate $x^r$, $r = 0, 1, \dots, d$, where $x^0$ denotes the time variable $t$ for ease of notation. Given data $\{u_{\x_j}(x_0^{k_0}, \dots, x_d^{k_d})\}$ for possibly different initial conditions $u_0 = \x_j$, $j = 1, 2, \dots, n_{\x}$ -- with all observations corresponding to a fixed spatial domain and fixed boundary conditions -- we collect the data in the tensors $\A$ and $\bb$ as
\begin{align}\label{eq:sindy_Ab_PDE}
    \begin{aligned}
        \A_i^{j, k_0, \dots, k_d} = \Aa^i(u_{\x_j})(x_0^{k_0}, \dots, x_d^{k_d}), \qquad \bb^{j, k_0, \dots, k_d} = \Tt_{(x_0^{k_0}, \dots, x_d^{k_d})}(u_{\x_j}),
    \end{aligned}
\end{align}
where $\Tt_{(x_0^{k_0}, \dots, x_d^{k_d})}$ denotes a numerical approximation of the time derivative at $(x_0^{k_0}, \dots, x_d^{k_d})$. Derivatives appearing in the candidate operators are likewise approximated numerically from the data. For PDEs with a different left-hand side than in equation \eqref{eq:PDE}, $\Tt_{(x_0^{k_0}, \dots, x_d^{k_d})}$ is modified accordingly. The resulting optimization problem for the library coefficients $\xi$ is
\begin{align*}
    \min_{\xi \in \RR^p} \sum_{j = 1}^{n_\x}\sum_{k_0 = 1}^{n_{x^0}}\cdots\sum_{k_d = 1}^{n_{x^d}} |{\textstyle\sum_{i = 1}^p}\A_i^{j, k_0, \dots, k_d}\xi^i - \bb^{j, k_0, \dots, k_d}|^2 + \lambda^2\|\xi\|_0.
\end{align*}

\textbf{Optimization problem.} The above formulations can be written in matrix-vector form by identifying the observation indices with a single index $q = q(j, k_0, \dots, k_d)$, for example $q(j, k) = (j - 1)n_t + k$ in the ODE case. Defining
\begin{align*}
    A_i^q := \A_i^{j, k_0, \dots, k_d}, \qquad b^q := \bb^{j, k_0, \dots, k_d},
\end{align*}
where $A \in \RR^{N \times p}$, $b \in \RR^N$, and where $N = n_\x n_t$ for ODEs and $N = n_\x n_{x^0}\cdots n_{x^d}$ for PDEs, we then obtain that $\xi$ is a solution to optimization problem \eqref{eq:sindy0}, i.e.,
\begin{align}\label{eq:sindy}
    \min_{\xi \in \RR^p} F(\xi), \qquad F(\xi) := \|A\xi - b\|_2^2 + \lambda^2\|\xi\|_0.
\end{align}
Throughout this work, we assume the following.
\begin{assumption}\label{ass:A}
The matrix $A$ has full column rank, i.e., $\rank(A) = p$.
\end{assumption}
\noindent This is equivalent to linear independence of the sampled library columns, and therefore it depends on both the library and the observations. It requires $N \geq p$ and ensures uniqueness of each least-squares fit. In the numerical experiments throughout this work, the library columns are moreover normalized to have unit Euclidean norm, with the coefficients rescaled accordingly. This allows coefficient magnitudes, and hence the threshold $\lambda$, to be compared on a common scale. We use $A$ and $\xi$ for the scaled quantities and transform the recovered coefficients back when presenting a learned equation.

\textbf{Systems of ODEs and PDEs.} For a system with $n$ state components, the library terms may depend on all components of $u$. In the ODE case, $\t^i : \RR^n \to \RR$ and
\begin{align}\label{eq:sindy_Ab_system}
    \A_i^{j, k} = \t^i(u_{\x_j}(t^k)), \qquad \bb_\ell^{j, k} = \Tt_{t^k}(u_{\x_j}^{\ell}), \qquad \ell = 1, \dots, n.
\end{align}
For PDEs, we similarly use
\begin{align*}
    \A_i^{j, k_0, \dots, k_d} = \Aa^i(u_{\x_j})(x_0^{k_0}, \dots, x_d^{k_d}), \qquad \bb_\ell^{j, k_0, \dots, k_d} = \Tt_{(x_0^{k_0}, \dots, x_d^{k_d})}(u_{\x_j}^{\ell}).
\end{align*}
Flattening the observation indices gives $b_\ell \in \RR^N$. We collect the coefficient vectors in $\Xi = (\xi^1, \dots, \xi^n) \in \RR^{p \times n}$, where $(\xi^\ell)^i$ is the coefficient of library term $i$ in component equation $\ell$, and find $\Xi$ by solving
\begin{align}\label{eq:sindy_system}
    \min_{\Xi \in \RR^{p \times n}} \sum_{\ell = 1}^n\|A\xi^\ell - b_\ell\|_2^2 + \lambda^2\|\Xi\|_0, \qquad \|\Xi\|_0 := \sum_{\ell = 1}^n\|\xi^\ell\|_0.
\end{align}
This problem separates into $n$ scalar regression problems, so STLSQ can be applied to each component independently. Note that this formulation uses a common candidate library for all components.

\subsection{Sequentially thresholded least-squares method}

\textbf{Algorithm.} A local minimizer of problem \eqref{eq:sindy} can be found using the sequentially thresholded least-squares (STLSQ) algorithm, which generates a sequence $\xi_k, \, k \geq 0$ through
\begin{align*}
    \begin{split}
        & \quad \xi_0 \in {\arg \min}_{\xi \in \RR^p} \|A \xi - b\|_2^2 \\
        & \begin{dcases}
            S_k = \{ 1 \leq i \leq p : |\xi_k^i| \geq \lambda\}, \\
            \xi_{k + 1} \in {\arg \min}_{\xi \in \RR^p : \supp(\xi) \subseteq S_k} \|A \xi - b\|_2^2,
        \end{dcases} \quad k \geq 0,
    \end{split}
\end{align*}
where we note that $\xi_k, \, k \geq 0$ is not necessarily unique. However, under Assumption \ref{ass:A}, we have that every $A_{S_k}$ has full column rank and therefore that each optimization problem in the iterative scheme can be solved uniquely and explicitly as
\begin{align*}
    \xi_{k + 1}|_{S_k} = {\arg \min}_{\xi \in \RR^{|S_k|}} \|A_{S_k}\xi - b\|_2^2 = A_{S_k}^\dagger b, \qquad \xi_{k + 1}|_{\overline{S_k}} = 0,
\end{align*}
thus reducing the iterative method to
\begin{align}\label{eq:iteration_scheme_sindy}
    \begin{split}
        & \quad \xi_0 = A^\dagger b \\
        & \begin{dcases}
            S_k = \{ 1 \leq i \leq p : |\xi_k^i| \geq \lambda\}, \\
            \xi_{k + 1}|_{S_k} = A_{S_k}^\dagger b, \qquad \xi_{k + 1}|_{\overline{S_k}} = 0,
        \end{dcases} \quad k \geq 0.
    \end{split}
\end{align}

The supports in iterative scheme \eqref{eq:iteration_scheme_sindy} satisfy $S_{k + 1} \subseteq S_k$. Once $S_k = S_{k + 1}$, the next least-squares problem is unchanged, so $\xi_{k + 2} = \xi_{k + 1}$. Hence the support shrinks strictly until a fixed point is reached, which happens in at most $p$ steps \cite[Theorem 2.1]{zhang2019convergence}. At such a fixed point $\xi_*$ with $S_* := \supp(\xi_*)$, we have
\begin{align*}
    \xi_*|_{S_*} = A_{S_*}^\dagger b, \qquad \xi_*|_{\overline{S_*}} = 0, \qquad \lambda \leq \min_{i \in S_*}|\xi_*^i|.
\end{align*}
Such a point then is a local minimizer of $F$ \cite[Theorem 2.5]{zhang2019convergence}. Under the additional assumption $\|A\|_2 \leq 1$, the objective $F$ is nonincreasing along the iterates, with strict decrease whenever $\xi_{k + 1} \neq \xi_k$ \cite[Theorem 2.4]{zhang2019convergence} (see also Lemma \ref{lem:monotone} and Remark \ref{rem:monotone}). For general $A$, the same result holds when the penalty $\lambda^2\|\xi\|_0$ is replaced by $\lambda^2\|A\|_2^2\|\xi\|_0$.

\textbf{Threshold selection.} Let $S_* \subseteq \{1, \dots, p\}$ denote the true support. The first thresholding step gives $S_0 = S_*$ if and only if
\begin{align}\label{eq:one_step_convergence}
    \max_{i \in \overline{S_*}}|(A^\dagger b)^i| < \lambda \leq \min_{i \in S_*}|(A^\dagger b)^i|.
\end{align}
For noiseless data $b = A\xi_*$, we have $A^\dagger b = \xi_*$, so condition \eqref{eq:one_step_convergence} reduces to $0 < \lambda \leq \min_{i \in S_*}|\xi_*^i|$. Refitting then returns $\xi_1 = \xi_*$, which already satisfies the threshold, giving exact recovery in one iteration as in \cite[Proposition 2.8]{zhang2019convergence}. For noisy or inexact representations, condition \eqref{eq:one_step_convergence} still gives $S_0 = S_*$, but the refit $\xi_1|_{S_*} = A_{S_*}^\dagger b$ may change the coefficient magnitudes. If all refitted coefficients still have magnitude at least $\lambda$, then $S_1 = S_*$ and $\xi_1$ is a fixed point. Otherwise, the next thresholding step removes at least one true term. The true support may also be recovered after several iterations if spurious terms are retained initially and removed later. However, the upper bound in \eqref{eq:one_step_convergence} remains necessary for eventual support recovery since a true term removed at the first thresholding step cannot re-enter. As mentioned earlier, the column normalization puts the library terms on a common scale, allowing their coefficients to be compared using a single threshold $\lambda$. Nevertheless, a threshold satisfying \eqref{eq:one_step_convergence} may not exist if the coefficient magnitudes are not sufficiently separated.

\section{Methods} \label{sec:family}

As explained in Section \ref{sec:intro}, the pSINDy framework aims to discover equations of the form
\begin{align*}
    \Bb_{\a}(u) = \sum_{i = 1}^p \xi^i \Aa^i(u),
\end{align*}
by jointly recovering the parameter $\a \in \RR^m$ and the library coefficients $\xi \in \RR^p$ while promoting sparsity in the latter. We additionally restrict $\a$ to a nonempty, closed, and convex set $\Omega_\a \subseteq \RR^m$ motivated by applications such as discovering fractional differential equations. For a scalar equation, the pSINDy framework then considers the optimization problem
\begin{align}\label{eq:sindy_alpha}
    \min_{\a \in \Omega_\a, \ \xi \in \RR^p} F(\a, \xi), \qquad F(\a, \xi) := \|A\xi - b_\a\|_2^2 + \lambda^2\|\xi\|_0.
\end{align}
Here, $A \in \RR^{N \times p}$ denotes the library matrix constructed in Section \ref{sec:preliminaries}, while $b_\a \in \RR^N$ denotes the discretization of $\Bb_\a$ applied to the given data. We assume throughout that $b_\a$ is continuous on $\Omega_\a$. For a system of $n$ equations, we instead consider
\begin{align}\label{eq:sindy_alpha_system}
    \min_{\a \in \Omega_\a, \;\Xi \in \RR^{p \times n}} F(\a, \Xi), \qquad F(\a, \Xi) := \sum_{\ell = 1}^n\|A\xi^\ell - b_{\a, \ell}\|_2^2 + \lambda^2\|\Xi\|_0,
\end{align}
where $\Xi = (\xi^1, \dots, \xi^n)$ is the coefficient matrix and $b_{\a, \ell}$ the discretization of the parameter-dependent operator in the $\ell$th component equation. The component equations are coupled through shared parameters. If instead each equation depends on its own parameters and these are constrained independently, problem \eqref{eq:sindy_alpha_system} decouples into $n$ problems of the form \eqref{eq:sindy_alpha}.

Since $A$ and $b_\a$ are constructed from the data, we also study how noise affects the local minimizers of $F$. In Appendix \ref{sec:stab}, we show that strict, nondegenerate local minimizers $(\a, \xi)$ with $\a$ in the interior of $\Omega_\a$ are stable under sufficiently small perturbations of $A$ and $b_\a$.

\subsection{Parameterized sequentially thresholded least-squares method}

\textbf{Algorithm.} To solve optimization problem \eqref{eq:sindy_alpha}, we construct a parameterized version of the STLSQ algorithm, which we denote by pSTLSQ. Given an initial parameter $\a_0 \in \Omega_\a$, we initialize the coefficient vector by $\xi_0 \in \arg\min_{\xi \in \RR^p} \|A\xi - b_{\a_0}\|_2^2$ and define the iterative scheme
\begin{align}\label{eq:psindy_naive}
    \begin{split}
        & \quad\a_0 \in \Omega_\a, \quad \xi_0 \in \arg\min_{\xi \in \RR^p} \|A\xi - b_{\a_0}\|_2^2, \\
        & \begin{dcases}
            S_k = \{1 \leq i \leq p : |\xi_k^i| \geq \lambda\}, \\
            (\a_{k + 1}, \xi_{k + 1}) \in \arg\min_{\a \in \Omega_\a, \ \xi \in \RR^p \, : \, \supp(\xi) \subseteq S_k} \|A\xi - b_\a\|_2^2,
        \end{dcases} \qquad k \geq 0,
    \end{split}
\end{align}
in a way similar to the STLSQ method for SINDy. As in the derivation of the STLSQ algorithm in Section \ref{sec:preliminaries}, the update in \eqref{eq:psindy_naive} is equivalent to
\begin{align}\label{eq:psindy_problem_support}
    (\a_{k + 1}, \xi_{k + 1}|_{S_k}) \in \arg\min_{\a \in \Omega_\a, \ \xi \in \RR^{|S_k|}} \|A_{S_k}\xi - b_\a\|_2^2, \qquad \xi_{k + 1}|_{\overline{S_k}} = 0.
\end{align}
By Assumption \ref{ass:A}, $A$ and $A_{S_k}$ have full column rank for every $k \geq 0$. Therefore, the initialization is uniquely given by $\xi_0 = A^{\dagger} b_{\a_0}$ and, for fixed $\a$, the minimizing coefficients in \eqref{eq:psindy_problem_support} are unique and given by $\phi_{S_k}(\a)$, where
\begin{align*}
    \phi_S(\a) := A_S^\dagger b_\a,
\end{align*}
where we write $\phi$ for the full library $\phi_{\{1, \dots, p\}}$ and define $\phi_{\emptyset} := 0$. Substituting the minimizing coefficients $\xi_{k + 1}|_{S_k} = \phi_{S_k}(\a_{k + 1})$ into optimization problem \eqref{eq:psindy_problem_support} reduces the problem to finding $\a_{k + 1}$ through solving
\begin{align}\label{eq:min_G}
    \min_{\a \in \Omega_\a} G_{S_k}(\a), \quad G_{S_k}(\a) := \|P_{S_k} b_\a\|_2^2, \quad P_{S_k} := I - A_{S_k}A_{S_k}^\dagger.
\end{align}
We take $P_\emptyset := I$ and assume that the minimum exists, which holds in particular when $\Omega_\a$ is compact since $b_\a$ is assumed to be continuous. Recalling that $A_{S}A_{S}^\dagger$ denotes the orthogonal projection onto $\text{range}(A_{S}), \, S \subseteq \{1, 2, \dots, p\}$, optimization problem \eqref{eq:min_G} corresponds to finding the parameter $\a_{k + 1} \in \Omega_\a$ that minimizes the squared distance of $b_\a$ from the subspace spanned by the active library columns, i.e., the columns indexed by the active support $S_k$. This elimination of the linear coefficients is the variable projection approach of \cite{golub1973differentiation}. Here $P_{S_k}$ does not depend on $\a$ because the parameter enters only through $b_\a$.

\textbf{Parameter-update schedules.} Solving optimization problem \eqref{eq:min_G} to high accuracy in every iteration can be computationally expensive and is not required for local convergence under suitable assumptions, as we will show in Section \ref{sec:geometry}. As a consequence, we allow each parameter update to consist of a prescribed number of iterations of an iterative method for solving optimization problem \eqref{eq:min_G}. To this end, let $\Ii_S : \Omega_\a \to \Omega_\a$ denote one iteration of such a method applied to $G_S$, and let $\Ii_S^h$, $h \geq 0$, denote $h$ such iterations with $\Ii_S^0$ the identity. Moreover, let $\Ii_S^\infty(\a) \in \arg\min_{\a' \in \Omega_\a}G_S(\a')$ denote an exact solution of optimization problem \eqref{eq:min_G}. Given a schedule $\nu = (\nu_{-1}, \nu_0, \nu_1, \dots)$ with $\nu_k \in \{0, 1, 2, \dots\} \cup \{\infty\}$, we then define the following iterative scheme:
\begin{align}\label{eq:iteration_scheme}
    \begin{split}
        & \quad \a_{-1} \in \Omega_\a, \quad \a_0 = \Ii_{\{1, \dots, p\}}^{\nu_{-1}}(\a_{-1}), \quad \xi_0 = A^{\dagger} b_{\a_0}, \\
        & \begin{dcases}
            S_k = \{ 1 \leq i \leq p : |\xi_k^i| \geq \lambda\}, \\
            \a_{k + 1} = \Ii_{S_k}^{\nu_k}(\a_k), \\
            \xi_{k + 1}|_{S_k} = A_{S_k}^\dagger b_{\a_{k + 1}}, \quad \xi_{k + 1}|_{\overline{S_k}} = 0,
        \end{dcases} \quad k \geq 0.
    \end{split}
\end{align}
We refer to this iterative scheme as the pSTLSQ algorithm, and to $\nu$ as the corresponding schedule. Here, $\nu_{-1}$ is the number of parameter updates on the full library before the first thresholding step, and $\nu_k$, $k \geq 0$, is the number of parameter updates after selecting $S_k$ and before refitting the coefficients. In this work we take $\nu_{-1} = 0$ without loss of generality since any parameter updates before the first thresholding step can be absorbed into the choice of initial parameter $\a_0$. With $\nu_k = 0$, $k \geq 0$, the scheme reduces to STLSQ at the fixed parameter $\a_0$. Taking $\nu_k = \infty, \, k \geq 0$ solves optimization problem \eqref{eq:min_G} exactly in every iteration, and we refer to this case as exact pSTLSQ. In a practical setting, $\nu_k = n_\a, \, k \geq 0$ for some positive integer $n_\a$, and we refer to this case as inexact pSTLSQ.

The following property of pSTLSQ -- in particular, of exact and inexact pSTLSQ -- follows directly from iteration scheme \eqref{eq:iteration_scheme}.

\begin{lemma}\label{lem:nested}
Let $\lambda > 0$. Then the supports $S_k$ of iteration scheme \eqref{eq:iteration_scheme} satisfy $S_{k + 1} \subseteq S_k$ for every $k \geq 0$ and every schedule $\nu$. In particular, the support changes at most $p$ times, i.e., there exists $K \geq 0$ such that $S_k = S_K$ for every $k \geq K$.
\end{lemma}
\begin{proof}
Since $\xi_{k + 1}|_{\overline{S_k}} = 0$ and $\lambda > 0$, we have $S_{k + 1} \subseteq S_k$. Since every strict inclusion removes at least one element from the support and $S_0 \subseteq \{1, 2, \dots, p\}$, we have $S_{k + 1} \neq S_k$ for at most $p$ indices $k$.
\end{proof}

As a consequence, a library term removed in a thresholding step cannot re-enter in a later iteration, so different initializations and schedules can lead to different supports. For exact pSTLSQ, the support stabilizes after at most $p$ iterations provided $G_S$ has a unique minimizer for every support $S \subseteq \{1, 2, \dots, p\}$ as shown below in Proposition \ref{prop:general_iter_scheme_minimizer_global}. For inexact pSTLSQ, the support may remain unchanged for many iterations before further parameter updates cause another term to be removed. Thus, although Lemma \ref{lem:nested} ensures eventual support stabilization, it gives no bound for the iteration at which this occurs.

\textbf{Systems.} For system \eqref{eq:sindy_alpha_system}, pSTLSQ thresholds and refits every component equation separately, and couples the components only through the parameter update. In this setting, $S = (S^1, \dots, S^n)$ denotes a tuple of supports, where $S^\ell \subseteq \{1, 2, \dots, p\}$ is the support of the $\ell$th component equation. For fixed $\a$, the minimizing coefficients are $\xi^\ell|_{S^\ell} = A_{S^\ell}^\dagger b_{\a, \ell}$ and $\xi^\ell|_{\overline{S^\ell}} = 0$, and eliminating them gives the reduced objective
\begin{align}\label{eq:min_G_system}
    G_S(\a) := \sum_{\ell = 1}^n\|P_{S^\ell}b_{\a, \ell}\|_2^2.
\end{align}
Iteration scheme \eqref{eq:iteration_scheme} then applies with the parameter update
\begin{align*}
    \a_{k + 1} = \Ii_{S_k}^{\nu_k}(\a_k), \qquad S_k = (S_k^1, \dots, S_k^n), \qquad S_k^\ell = \{1 \leq i \leq p : |(\xi_k^\ell)^i| \geq \lambda\}.
\end{align*}
The results of this section extend to systems with $S$ interpreted as a tuple, support inclusions interpreted componentwise, and $|S|$ and $p$ replaced by $\sum_{\ell = 1}^n|S^\ell|$ and $np$ in the bounds on support changes and iterations. In the decoupled case, each component can be treated independently.

\subsection{Local minimizers and fixed points}
We first characterize the local minimizers of the pSINDy optimization problem through the reduced objectives $G_S$. We then relate this characterization to the fixed points of pSTLSQ.

\begin{theorem} \label{thm:fixed_point_local_minimizer}
Let $(\a_*, \xi_*) \in \Omega_\a \times \RR^p$ and $S_* := \supp(\xi_*)$. Suppose that Assumption \ref{ass:A} holds, that $\lambda > 0$, and that $b_\a$ is continuous at $\a_*$. Then $(\a_*, \xi_*)$ is a (strict) local minimizer of optimization problem \eqref{eq:sindy_alpha} if and only if $\xi_*|_{S_*} = A_{S_*}^{\dagger}b_{\a_*}$ and $\a_*$ is a (strict) local minimizer of $G_{S_*}$ on $\Omega_\a$.
\end{theorem}

\begin{proof}
We first assume that $\xi_*|_{S_*} = A_{S_*}^{\dagger}b_{\a_*}$ and that $\a_*$ is a local minimizer of $G_{S_*}$ on $\Omega_\a$. Then there exists $\eps_\a > 0$ such that $G_{S_*}(\a) \geq G_{S_*}(\a_*)$ for every $\a \in \Omega_\a \cap B_{\eps_\a}(\a_*)$. We claim that
\begin{align}\label{eq:inclusion_S}
    S_* \subseteq \supp(\xi)
\end{align}
for every $\xi \in B_\eps(\xi_*)$ for $\eps > 0$ sufficiently small. If $S_* = \emptyset$, the claim is trivial. Otherwise, choosing $\eps < \min_{i \in S_*} |\xi_*^i|$ gives that $|\xi^i| \geq |\xi_*^i| - |\xi^i - \xi_*^i| > 0$ for every $i \in S_*$ and every $\xi \in B_\eps(\xi_*)$, which proves inclusion \eqref{eq:inclusion_S}.

Let $(\a, \xi)$ be an arbitrary point satisfying $\a \in \Omega_\a \cap B_{\eps_\a}(\a_*)$ and $\xi \in B_{\eps}(\xi_*)$. By inclusion \eqref{eq:inclusion_S}, we have that either $S_* = \supp(\xi)$ or $S_* \subset \supp(\xi)$. If $\supp(\xi) = S_*$, we have
\begin{align*}
    \|A\xi - b_\a\|_2^2 = \|A_{S_*}\xi|_{S_*} - b_\a\|_2^2 \geq G_{S_*}(\a) \geq G_{S_*}(\a_*),
\end{align*}
where the last inequality follows from the local minimality of $\a_*$. Since $\|\xi\|_0 = |S_*| = \|\xi_*\|_0$ and $\|A\xi_* - b_{\a_*}\|_2^2 = G_{S_*}(\a_*)$, this implies that
\begin{align}\label{eq:local1}
    F(\a, \xi) \geq F(\a_*, \xi_*).
\end{align}
If $S_* \subset \supp(\xi)$, then $\|\xi\|_0 \geq |S_*| + 1$. Since $\|A\xi - b_\a\|_2^2$ is continuous in both $\a$ and $\xi$ at $(\a_*, \xi_*)$, there exists $\kappa > 0$ such that
\begin{align*}
    \left| \|A\xi - b_\a\|_2^2 - \|A\xi_* - b_{\a_*}\|_2^2 \right| < \lambda^2,
\end{align*}
for all $(\a, \xi) \in B_{\kappa}((\a_*, \xi_*)) \cap (\Omega_\a \times \RR^p)$, and thus that
\begin{align}\label{eq:local2}
    F(\a, \xi) > F(\a_*, \xi_*),
\end{align}
for all $(\a, \xi) \in B_{\kappa}((\a_*, \xi_*)) \cap (\Omega_\a \times \RR^p)$ with $S_* \subset \supp(\xi)$. Together, inequalities \eqref{eq:local1} and \eqref{eq:local2} imply that $(\a_*, \xi_*)$ is a local minimizer of $F$ on $\Omega_\a \times \RR^p$. Moreover, if $\a_*$ is a strict local minimizer of $G_{S_*}$, then, after possibly decreasing $\eps_\a$, inequality \eqref{eq:local1} is strict whenever $(\a, \xi) \neq (\a_*, \xi_*)$. Indeed, the second inequality above is then strict for $\a \neq \a_*$, whereas for $\a = \a_*$ and $\xi \neq \xi_*$ the first inequality is strict, since $A_{S_*}$ has full column rank by Assumption \ref{ass:A}, implying that the least-squares solution $\xi_*|_{S_*} = A_{S_*}^\dagger b_{\a_*}$ is unique. Hence, $(\a_*, \xi_*)$ is a strict local minimizer of $F$ on $\Omega_\a \times \RR^p$ in this case.

Conversely, let $(\a_*, \xi_*)$ be a (strict) local minimizer of $F$ on $\Omega_\a \times \RR^p$. For $\xi \in B_\eps(\xi_*)$ with $\xi|_{\overline{S_*}} = 0$, inclusion \eqref{eq:inclusion_S} gives $\supp(\xi) = S_*$, and hence $F(\a_*, \xi) = \|A_{S_*}\xi|_{S_*} - b_{\a_*}\|_2^2 + \lambda^2|S_*|$. Local minimality of $(\a_*, \xi_*)$ therefore implies that $\xi_*|_{S_*}$ is a local minimizer of the function $\z \mapsto \|A_{S_*}\z - b_{\a_*}\|_2^2$. By convexity of the latter function, this local minimizer is also global. Since $A_{S_*}$ has full column rank, this gives $\xi_*|_{S_*} = A_{S_*}^{\dagger}b_{\a_*}$. For $\a \in \Omega_\a$, define $\tilde\xi(\a) \in \RR^p$ by $\tilde\xi(\a)|_{S_*} := \phi_{S_*}(\a)$ and $\tilde\xi(\a)|_{\overline{S_*}} := 0$, so that $\tilde\xi(\a_*) = \xi_*$. Since $b_\a$ is continuous at $\a_*$, we have $\tilde\xi(\a) \to \xi_*$ as $\a \to \a_*$. By inclusion \eqref{eq:inclusion_S} and the definition of $\tilde\xi(\a)$, it follows that $\supp(\tilde\xi(\a)) = S_*$, and hence $F(\a, \tilde\xi(\a)) = G_{S_*}(\a) + \lambda^2|S_*|$, for every $\a \in \Omega_\a$ sufficiently close to $\a_*$. Local minimality of $(\a_*, \xi_*)$ therefore gives $G_{S_*}(\a) \geq G_{S_*}(\a_*)$ for every such $\a$, with strict inequality for $\a \neq \a_*$ in the strict case.
\end{proof}

Letting $(\a_*, \xi_*)$ be a fixed point of pSTLSQ, the coefficient refit thus gives $\xi_*|_{S_*} = A_{S_*}^{\dagger}b_{\a_*}$. By Theorem \ref{thm:fixed_point_local_minimizer}, $(\a_*, \xi_*)$ is therefore a (strict) local minimizer of $F$ if and only if $\a_*$ is a (strict) local minimizer of $G_{S_*}$ on $\Omega_\a$. As a consequence, every fixed point of exact pSTLSQ is a local minimizer of $F$. For inexact pSTLSQ, this conclusion requires further conditions on the parameter update, as discussed in Section \ref{sec:geometry}.

The characterization in Theorem \ref{thm:fixed_point_local_minimizer} also shows that the local minimizers of $F$ are independent of $\lambda > 0$. However, $\lambda$ affects which of these are fixed points of pSTLSQ and which are global minimizers of $F$. Example \ref{ex:global} illustrates this. Moreover, for every support $S$ and every local minimizer $\a$ of $G_S$ on $\Omega_\a$ for which $\phi_S(\a)$ has no zero entries, the point $(\a, \xi)$ with $\xi|_S = \phi_S(\a)$ and $\xi|_{\overline S} = 0$ is a local minimizer of $F$ by Theorem \ref{thm:fixed_point_local_minimizer}.

\subsection{Exact pSTLSQ} \label{sec:analysis}

We first discuss properties of exact pSTLSQ, where we recall that exact pSTLSQ is the pSTLSQ iteration scheme \eqref{eq:iteration_scheme} with $\nu_k = \infty, \, k \geq 0$. In particular, we establish finite termination of exact pSTLSQ under uniqueness of the reduced minimizers. We then illustrate how initialization and nonuniqueness of the reduced minimizers affect the iteration, and why finite termination does not guarantee recovery of a global minimizer of the pSINDy problem.

\textbf{Finite termination under unique minimizers.} The following result shows that exact pSTLSQ reaches a fixed point in at most $p + 1$ iterations provided each reduced optimization problem has a unique minimizer, and that this fixed point is a local minimizer of $F$ by Theorem \ref{thm:fixed_point_local_minimizer} provided $b_\a$ is continuous.

\begin{proposition} \label{prop:general_iter_scheme_minimizer_global}
Assume that $\lambda > 0$ and that Assumption \ref{ass:A} holds. Suppose that, for every support $S \subseteq \{1, 2, \dots, p\}$, the function $G_S$ has a unique minimizer $\a_*(S)$ on $\Omega_\a$. Then exact pSTLSQ reaches a fixed point in at most $p + 1$ steps. If additionally $b_\a \in \Cc(\Omega_\a)$, this fixed point is a local minimizer of optimization problem \eqref{eq:sindy_alpha}.
\end{proposition}

\begin{proof}
We first show that $S_k = S_{k + 1}$ implies $S_k = S_l$ for every $l \geq k$. Indeed, if $S_k = S_{k + 1}$, then
\begin{align*}
    \a_{k + 1} = \arg\min_{\a \in \Omega_\a} G_{S_k}(\a) = \arg\min_{\a \in \Omega_\a} G_{S_{k + 1}}(\a) = \a_{k + 2},
\end{align*}
where the last equality follows from the uniqueness of the minimizer of $G_{S_{k + 1}}$. Hence
\begin{align*}
    \xi_{k + 1}|_{S_k} = A_{S_k}^{\dagger} b_{\a_{k + 1}} = A_{S_{k + 1}}^{\dagger} b_{\a_{k + 2}} = \xi_{k + 2}|_{S_{k + 1}}.
\end{align*}
Since $S_k = S_{k + 1}$, we obtain $\xi_{k + 1} = \xi_{k + 2}$ and therefore $S_{k + 1} = S_{k + 2}$. By induction, we conclude that $S_k = S_l$ for every $l \geq k$. Moreover, by Lemma \ref{lem:nested}, the support changes at most $p$ times. Hence there exists some $l \in \{0, 1, \dots, p\}$ such that $S_l = S_{l + 1}$, and by the earlier induction argument we thus have $S_l = S_{l + q}$ for every $q \geq 0$. Denoting the stabilized support by $S_*$, uniqueness of the minimizer of $G_{S_*}$ implies that $\a_k = \a_*(S_*)$ for every $k \geq l + 1$. Consequently,
\begin{align*}
    \xi_k|_{S_*} = A_{S_*}^{\dagger}b_{\a_*(S_*)},
\end{align*}
for every $k \geq l + 1$, and therefore the iteration has reached a fixed point after at most $p + 1$ steps. If $b_\a \in \Cc(\Omega_\a)$ this fixed point is a local minimizer of $F$ by Theorem \ref{thm:fixed_point_local_minimizer} since $\a_*(S_*)$ minimizes $G_{S_*}$.
\end{proof}

\textbf{Limitations of global recovery.} Proposition \ref{prop:general_iter_scheme_minimizer_global} requires $G_S$ to have a unique minimizer for every support $S \subseteq \{1, 2, \dots, p\}$. Without this assumption, the next iterate may depend on which minimizer is selected. Moreover, even a unique minimizer on the selected support need not give a global minimizer of $F$, since thresholding determines which supports remain available. The following example illustrates how the initialization affects the initial support and how nonuniqueness of a reduced minimizer can lead to a fixed point that is not globally optimal.

\begin{example} \label{ex:global}
Consider recovering $\a \in \Omega_\a = \RR$ and $\xi = (\xi^1, \xi^2)$ in
\begin{align*}
    \Bb_\a(u) = \xi^1x + \xi^2y, \qquad \Bb_\a(u) := \a x + (1 - \a)y, \qquad u(t) = (x(t), y(t)).
\end{align*}
Suppose that the equation is sampled at two points satisfying
\begin{align*}
    u(t_1) = (1, 0), \qquad u(t_2) = (0, 1).
\end{align*}
Then $A = I$ and $b_\a = (\a, 1 - \a)^\top$. Taking $\lambda = \frac{1}{2}$ gives
\begin{align*}
    F(\a, \xi) = (\xi^1 - \a)^2 + (\xi^2 - (1 - \a))^2 + \frac{1}{4}\|\xi\|_0.
\end{align*}
For $\xi \neq 0$, we have $F(\a, \xi) \geq \frac{1}{4}$, whereas $F(\a, 0) = \a^2 + (1 - \a)^2 \geq \frac{1}{2}$. The value $\frac{1}{4}$ is attained precisely at
\begin{align*}
    (\a, \xi) = (1, (1, 0)^\top) \qquad\text{and}\qquad (\a, \xi) = (0, (0, 1)^\top),
\end{align*}
so these are the two global minimizers of $F$. We now apply exact pSTLSQ with $\nu_{-1} = 0$, so that
\begin{align*}
    \xi_0 = (\a_0, 1 - \a_0)^\top.
\end{align*}
Moreover, when the reduced objective $G_{S_k}$ is constant for some $k$, we set $\a_{k + 1} = \a_k$.

If $\a_0 = \frac{3}{4}$, then $\xi_0 = (\frac{3}{4}, \frac{1}{4})^\top$, so the first thresholding step gives $S_0 = \{1\}$. Since
\begin{align*}
    G_{\{1\}}(\a) = (1 - \a)^2,
\end{align*}
the exact parameter update gives $\a_1 = 1$, and the coefficient refit gives $\xi_1 = (1, 0)^\top$. Thus, pSTLSQ reaches one of the two global minimizers after one iteration and remains there. Similarly, if $\a_0 = \frac{1}{4}$, then $S_0 = \{2\}$, and the algorithm reaches the other global minimizer.

If instead $\a_0 = \frac{1}{2}$, then $\xi_0 = (\frac{1}{2}, \frac{1}{2})^\top$. Both coefficients are retained because their magnitudes equal $\lambda$, so $S_0 = \{1, 2\}$, for which
\begin{align*}
    G_{\{1, 2\}}(\a) \equiv 0.
\end{align*}
The exact parameter update can therefore choose any $\a \in \Omega_\a$. With the choice specified above, both the parameter and the coefficient vector $\xi_0$ remain unchanged. Hence, $(\frac{1}{2}, (\frac{1}{2}, \frac{1}{2})^\top)$ is a fixed point. By Theorem \ref{thm:fixed_point_local_minimizer}, this is then a local minimizer, but note that it is not a global minimizer as $F(\a_0, \xi_0) = \frac{1}{2}$ instead of $\frac{1}{4}$. Thus, exact minimization on the retained support does not guarantee a global minimizer of $F$, and the outcome of exact pSTLSQ can depend on the initialization. $\hfill\blacksquare$
\end{example}

\begin{remark}
More generally, for $\lambda > 0$, the global minimum of $F$ is $\min\{\lambda^2, 1/2\}$. Thus, the two one-term solutions are global minimizers when $0 < \lambda < 1/\sqrt{2}$, whereas $(\a, \xi) = (1/2, 0)$ is the unique global minimizer when $\lambda > 1/\sqrt{2}$. At equality, all three are global minimizers. Moreover, the threshold $\lambda$ affects the support selected by pSTLSQ. For example, when $\a_0 = 1/2$, both coefficients are retained if $\lambda \leq 1/2$, whereas both are discarded if $\lambda > 1/2$.
\end{remark}

The example and the remark illustrate that the local minimizer reached by exact pSTLSQ can depend on both the initialization $\a_0$ and the threshold $\lambda$, and does not need to be a global minimizer of $F$. Note that the uniqueness assumption in Proposition \ref{prop:general_iter_scheme_minimizer_global} removes the ambiguity in the exact parameter updates and guarantees finite termination but does not guarantee global optimality.

\begin{remark}\label{prop:general_iter_scheme_minimizer}
Let $(\a_*, \xi_*)$ be a strict local minimizer of $F$, and let $S_* := \supp(\xi_*)$. Suppose that Assumption \ref{ass:A} holds and that
\begin{align}\label{eq:cond_xi}
    \max_{i \notin S_*}|\phi^i(\a_*)| < \lambda < \min_{i \in S_*}|\phi^i(\a_*)|, \qquad \lambda \leq \min_{i \in S_*}|\xi_*^i|.
\end{align}
By continuity and Theorem \ref{thm:fixed_point_local_minimizer}, there exists $\eps > 0$ such that $\a_*$ is the unique minimizer of $G_{S_*}$ on $O_\a := B_\eps(\a_*) \cap \Omega_\a$, and the first threshold selects $S_0 = S_*$ for every $\a_0 \in O_\a$. For $\a_0 \in O_\a$, exact pSTLSQ with parameter minimization restricted to $O_\a$ therefore recovers $(\a_*, \xi_*)$ in one iteration and leaves it fixed afterwards. In other words, local convergence of exact pSTLSQ is guaranteed if one initializes $\a_0$ in $O_\a$ and chooses the minimizer $\a_{k + 1}$ of $G_{S_k}$ within this neighborhood for each iteration $k$. This argument also underlies the local convergence of inexact pSTLSQ, which we will make precise in Theorem \ref{thm:local_inexact}.

The neighborhood in which $S_0 = S_*$ holds can moreover be quantified if $b_\a$ is Lipschitz on $B_\rho(\a_*) \cap \Omega_\a$ with constant $L_b$. Indeed, since $\phi(\a) = A^\dagger b_\a$,
\begin{align*}
    \|\phi(\a_0) - \phi(\a_*)\|_\infty \leq \|\phi(\a_0) - \phi(\a_*)\|_2 \leq L_b\|A^\dagger\|_2\|\a_0 - \a_*\|_2.
\end{align*}
We define
\begin{align*}
    \gamma := \min\left\{ \lambda - \max_{i \notin S_*}|\phi^i(\a_*)|, \min_{i \in S_*}|\phi^i(\a_*)| - \lambda \right\} > 0
\end{align*}
which is the smallest distance between the coefficient magnitudes at $\a_*$ and the threshold $\lambda$ so that perturbations of $\phi(\a_*)$ smaller than $\gamma$ in the infinity norm leave the thresholded support unchanged, and thus $S_0 = S_*$, whenever
\begin{align*}
    \a_0 \in B_\rho(\a_*) \cap \Omega_\a, \qquad L_b\|A^\dagger\|_2\|\a_0 - \a_*\|_2 < \gamma.
\end{align*}

\end{remark}

\subsection{Inexact pSTLSQ}\label{sec:geometry}

We now consider inexact pSTLSQ, where $\nu_k = n_\a$ for $k \geq 0$ with $n_\a$ a positive integer. Each parameter update therefore consists of $n_\a$ steps of an iterative method $\Ii_{S_k}$ for optimization problem \eqref{eq:min_G}. Even if $S_k = S_{k + 1}$, the parameter may continue to change because $\a_{k + 1}$ might not minimize $G_{S_k}$. These parameter changes can moreover cause the refitted coefficients to fall below the threshold, leading to further support reductions in later iterations. The finite-termination argument of Proposition \ref{prop:general_iter_scheme_minimizer_global} therefore no longer applies.

A second difference concerns fixed points. At a fixed point of exact pSTLSQ, $\a_*$ minimizes $G_{S_*}$, so the fixed point is a local minimizer of $F$ by Theorem \ref{thm:fixed_point_local_minimizer}. An inexact parameter update may instead leave $\a_*$ unchanged without minimizing $G_{S_*}$. For example, projected gradient descent leaves $\a_*$ unchanged when $\nabla G_{S_*}(\a_*) = 0$, even if this point is a maximum or saddle point. Thus, a fixed point of inexact pSTLSQ need not be a local minimizer of $F$.

Additionally, if the iterative method contracts near a local minimizer, then the parameter iterates stay near that minimizer when initialized sufficiently close. This allows us to establish local convergence of inexact pSTLSQ without restricting the parameter domain, whereas the argument for exact pSTLSQ in Remark \ref{prop:general_iter_scheme_minimizer} requires minimization within a neighborhood of the local minimizer.

\textbf{Inner iterations.} When proving that $F$ is nonincreasing, we require the inner steps to satisfy
\begin{align}\label{eq:inner_decrease}
    G_S(\Ii_S(\a)) \leq G_S(\a)
\end{align}
for every $\a \in \Omega_\a$ and fixed support $S$. In the numerical examples, we take $n_\a = 1$ and use projected gradient descent for continuously differentiable $b_\a$,
\begin{align}\label{eq:inexact_map}
    \Ii_S(\a) = \pi_{\Omega_\a}\bigl(\a - \eta\nabla G_S(\a)\bigr),
\end{align}
where $\eta > 0$ denotes the learning rate and $\pi_{\Omega_\a}$ is the Euclidean projection onto $\Omega_\a$. Differentiating $G_S(\a) = \|P_Sb_\a\|_2^2$ and using $P_S^\top P_S = P_S$ gives
\begin{align}\label{eq:gradient}
    \nabla G_S(\a) = 2\nabla b_\a^\top P_Sb_\a,
\end{align}
with $\nabla b_\a \in \RR^{N \times m}$ the Jacobian of $b_\a$. If $b_\a \in \Cc^2(\Omega_\a)$ and $\nabla G_S$ is Lipschitz continuous on $\Omega_\a$ with constant $L_S > 0$, then projected gradient descent satisfies condition \eqref{eq:inner_decrease} for $0 < \eta < 2/L_S$ by Lemma \ref{lem:alpha_decreasing_fixedS_GD}(a) in Appendix \ref{ap:b}.

\textbf{Decrease of the objective.} Under the additional condition $\|A\|_2 \leq 1$, the following lemma shows that $F$ is also nonincreasing.

\begin{lemma} \label{lem:monotone}
Let $\lambda > 0$ and $\|A\|_2 \leq 1$, and consider iteration scheme \eqref{eq:iteration_scheme} with a schedule whose parameter updates satisfy $G_{S_k}(\a_{k + 1}) \leq G_{S_k}(\a_k)$ for every $k \geq 0$. Then $F(\a_{k + 1}, \xi_{k + 1}) \leq F(\a_k, \xi_k)$ for every $k \geq 0$.
\end{lemma}

\begin{proof}
Set $S_{-1} := \{1, \dots, p\}$ and $R_k := S_{k - 1} \setminus S_k$. The residual $r_k := A\xi_k - b_{\a_k}$ satisfies $\|r_k\|_2^2 = G_{S_{k - 1}}(\a_k)$ and $A_{S_{k - 1}}^\top r_k = 0$. Let $\tilde\xi$ be obtained from $\xi_k$ by setting the entries indexed by $R_k$ to zero. Since $\supp(\tilde\xi) \subseteq S_k$, $R_k \subseteq S_{k - 1}$, $\|A\|_2 \leq 1$, and $|\xi_k^i| < \lambda$ for every $i \in R_k$, we obtain
\begin{align*}
    G_{S_k}(\a_k) & \leq \|A\tilde\xi - b_{\a_k}\|_2^2 = \|r_k\|_2^2 + \|A_{R_k}\xi_k|_{R_k}\|_2^2 \\
    & \leq G_{S_{k - 1}}(\a_k) + \lambda^2|R_k \cap \supp(\xi_k)|.
\end{align*}
Moreover, $S_k \subseteq \supp(\xi_k) \subseteq S_{k - 1}$ gives $\|\xi_k\|_0 = |S_k| + |R_k \cap \supp(\xi_k)|$. Since $\supp(\xi_{k + 1}) \subseteq S_k$, we therefore obtain
\begin{align*}
    F(\a_{k + 1}, \xi_{k + 1}) \leq G_{S_k}(\a_{k + 1}) + \lambda^2|S_k| \leq G_{S_k}(\a_k) + \lambda^2|S_k| \leq G_{S_{k - 1}}(\a_k) + \lambda^2\|\xi_k\|_0 = F(\a_k, \xi_k).
\end{align*}
\end{proof}

\begin{remark}\label{rem:monotone}
For general $A$, the same proof gives monotonicity with $\lambda^2$ replaced by $\lambda^2\|A\|_2^2$ in the regularization term. This is the relevant case for the numerical experiments. Columns of unit Euclidean norm give $\|A\|_2 \geq 1$, with equality only if the columns are orthonormal.
\end{remark}

The parameter-update condition holds for exact pSTLSQ because each update minimizes $G_{S_k}$. For inexact pSTLSQ, condition \eqref{eq:inner_decrease} ensures this after any number of inner steps. Lemma \ref{lem:monotone} therefore shows that $F$ does not increase in either version of pSTLSQ when $\|A\|_2 \leq 1$.

\textbf{Convergence after support stabilization.} Decrease of the objective does not by itself imply convergence of the iterates. By Lemma \ref{lem:nested}, however, there exists $K \geq 0$ such that $S_k = S_*$ for every $k \geq K$. For $k \geq K$, the method performs $n_\a$ inner steps on the same reduced objective $G_{S_*}$ and refits the coefficients on $S_*$. If $\a_k \to \a_*$ and $b_\a$ is continuous at $\a_*$, then $\xi_k \to \xi_*$, with $\xi_*|_{S_*} = \phi_{S_*}(\a_*)$ and $\xi_*|_{\overline{S_*}} = 0$. Since $|\xi_k^i| \geq \lambda$ for $i \in S_*$ and $k \geq K + 1$, the limit satisfies $|\xi_*^i| \geq \lambda$ for $i \in S_*$, so that $\supp(\xi_*) = S_*$ when $\lambda > 0$. Under Assumption \ref{ass:A}, Theorem \ref{thm:fixed_point_local_minimizer} then shows that $(\a_*, \xi_*)$ is a local minimizer of $F$ if and only if $\a_*$ is a local minimizer of $G_{S_*}$.

If $\Ii_{S_*}^{n_\a}$ is continuous at $\a_*$, taking limits in $\a_{k + 1} = \Ii_{S_*}^{n_\a}(\a_k)$ gives $\a_* = \Ii_{S_*}^{n_\a}(\a_*)$. Together with the coefficient refit, this makes $(\a_*, \xi_*)$ a fixed point of inexact pSTLSQ. This continuity holds for projected gradient descent with $b_\a \in \Cc^1(\Omega_\a)$ because the gradient and projection are continuous. For projected gradient descent under the assumptions of Lemma \ref{lem:alpha_decreasing_fixedS_GD}(a), the fixed-point condition $\a_* = \Ii_{S_*}^{n_\a}(\a_*)$ reduces to $\a_* = \Ii_{S_*}(\a_*)$, since every inner step that changes the parameter strictly decreases $G_{S_*}$. At an interior point, this only gives $\nabla G_{S_*}(\a_*) = 0$ and does not guarantee local minimality of $F$.

\textbf{Local convergence.} Near a local minimizer of $G_{S_*}$, a local contraction of the parameter update, together with the separation condition \eqref{eq:cond_xi}, suffices for convergence.

\begin{theorem} \label{thm:local_inexact}
Let $(\a_*, \xi_*)$ be a local minimizer of optimization problem \eqref{eq:sindy_alpha}, let $S_* := \supp(\xi_*)$, let $\lambda > 0$, assume that $b_\a$ is continuous at $\a_*$, and assume that condition \eqref{eq:cond_xi} holds with $\lambda < \min_{i \in S_*}|\xi_*^i|$. Suppose that Assumption \ref{ass:A} holds and that there exist $\rho_0 > 0$ and $\vartheta \in [0, 1)$ such that
\begin{align}\label{eq:contraction}
    \|\Ii_{S_*}(\a) - \a_*\|_2 \leq \vartheta\|\a - \a_*\|_2
\end{align}
for every $\a \in B_{\rho_0}(\a_*) \cap \Omega_\a$. Then there exists $\rho > 0$ such that, for every integer $n_\a \geq 1$, the inexact pSTLSQ iteration initialized by $\a_0 \in B_\rho(\a_*) \cap \Omega_\a$ and $\xi_0 = A^\dagger b_{\a_0}$ satisfies $S_k = S_*$ for every $k \geq 0$, and $\a_k$ converges linearly to $\a_*$. In particular, $(\a_k, \xi_k) \to (\a_*, \xi_*)$.
\end{theorem}

\begin{proof}
By Theorem \ref{thm:fixed_point_local_minimizer}, $\xi_*|_{S_*} = \phi_{S_*}(\a_*)$. By condition \eqref{eq:cond_xi} with $\lambda < \min_{i \in S_*}|\xi_*^i|$ and continuity of $\phi$ and $\phi_{S_*}$ at $\a_*$, there exists $\rho > 0$ such that
\begin{align}\label{eq:margins_rho}
    \max_{i \in \overline{S_*}}|\phi^i(\a)| < \lambda < \min_{i \in S_*}|\phi^i(\a)|, \qquad \lambda < \min_{i \in S_*}|\phi_{S_*}^i(\a)|,
\end{align}
for every $\a \in B_\rho(\a_*) \cap \Omega_\a$. Decreasing $\rho$ if necessary, we may assume that $\rho \leq \rho_0$, so that \eqref{eq:contraction} holds for every $\a \in B_\rho(\a_*) \cap \Omega_\a$. Taking $\a = \a_*$ shows that $\a_*$ is a fixed point of $\Ii_{S_*}$. In particular, $\Ii_{S_*}$ maps $B_\rho(\a_*) \cap \Omega_\a$ into itself.

We now show by induction that $\a_k \in B_\rho(\a_*) \cap \Omega_\a$ and $S_k = S_*$ for every $k \geq 0$. For $k = 0$, this follows from $\xi_0 = \phi(\a_0)$ and inequalities \eqref{eq:margins_rho}. If it holds for some $k \geq 0$, then $\a_{k + 1} = \Ii_{S_*}^{n_\a}(\a_k) \in B_\rho(\a_*) \cap \Omega_\a$ with $\|\a_{k + 1} - \a_*\|_2 \leq \vartheta^{n_\a}\|\a_k - \a_*\|_2$. Moreover, $\xi_{k + 1}|_{S_*} = \phi_{S_*}(\a_{k + 1})$ and $\xi_{k + 1}|_{\overline{S_*}} = 0$, so that \eqref{eq:margins_rho} gives $S_{k + 1} = S_*$. Hence $\|\a_k - \a_*\|_2 \leq \vartheta^{kn_\a}\|\a_0 - \a_*\|_2$, and $\xi_k \to \xi_*$ follows from continuity of $\phi_{S_*}$ at $\a_*$.
\end{proof}

Theorem \ref{thm:local_inexact} gives local convergence for the inexact iteration without assuming convergence of the parameters. The separation conditions ensure that $S_k = S_*$ throughout the iteration, and \eqref{eq:contraction} then gives linear convergence of the parameter. If $b_\a$ is locally Lipschitz near $\a_*$, then so is $\phi_{S_*}$, and the coefficient errors also decay at a linear rate from $k = 1$ onward. The same proof applies to every schedule with $1 \leq \nu_k < \infty$ for $k \geq 0$, since each parameter update then satisfies $\|\a_{k + 1} - \a_*\|_2 \leq \vartheta^{\nu_k}\|\a_k - \a_*\|_2 \leq \vartheta\|\a_k - \a_*\|_2$. 

\textbf{Contraction for projected gradient descent.} We now examine the contraction condition \eqref{eq:contraction} of Theorem \ref{thm:local_inexact} for projected gradient descent near the local minimizer $(\a_*, \xi_*)$ considered in Theorem \ref{thm:local_inexact}. By Theorem \ref{thm:fixed_point_local_minimizer}, $\a_*$ is a local minimizer of $G_{S_*}$ and hence a fixed point of projected gradient descent. For $b_\a \in \Cc^2(\Omega_\a)$, Lemma \ref{lem:alpha_decreasing_fixedS_GD}(c) then gives that condition \eqref{eq:contraction} holds whenever $\nabla^2G_{S_*}(\a_*) \succ 0$ and $0 < \eta \leq \tilde\eta$.

To examine the positive definiteness of the Hessian, we compute
\begin{align}\label{eq:hessian}
    \nabla^2G_S(\a) = 2(P_S\nabla b_\a)^\top(P_S\nabla b_\a) + 2\sum_{q = 1}^N\nabla^2b_\a^q(P_Sb_\a)^q.
\end{align}
The first term is positive semidefinite, whereas the second term does not need to be so. However, if $b_{\a_*}$ can be represented by the active library columns, i.e.,
\begin{align}\label{eq:exact_representation}
    P_{S_*}b_{\a_*} = 0,
\end{align}
the Hessian at $\a_*$ simplifies to
\begin{align}\label{eq:hess_at_sol}
    \nabla^2 G_{S_*}(\a_*) = 2J_{S_*}^\top J_{S_*},
\end{align}
where $J_{S_*} := P_{S_*}\nabla b_{\a_*}$. In this case, the Hessian is positive definite if and only if $J_{S_*}$ has full column rank. Thus, provided equation \eqref{eq:exact_representation} holds, full column rank of $J_{S_*}$ guarantees the local contraction of projected gradient descent for sufficiently small positive step sizes.

\subsection{Conditions for support recovery}\label{sec:threshold}

We now derive conditions on the threshold $\lambda$ for support recovery by (both exact and inexact) pSTLSQ. To this end, let $\a_* \in \Omega_\a$ be the true parameter value and $S_* \subseteq \{1, \dots, p\}$ the true support.

We first consider the case $\a_0 = \a_*$, for which the first thresholding step coincides with STLSQ at $\a_*$. As in condition \eqref{eq:one_step_convergence}, the initial support $S_0 = \{1 \leq i \leq p : |\phi^i(\a_0)| \geq \lambda\}$ then equals $S_*$ if and only if
\begin{align}\label{eq:cond_lambda}
    \max_{i \in \overline{S_*}}|\phi^i(\a_*)| < \lambda \leq \min_{i \in S_*}|\phi^i(\a_*)|.
\end{align}
If, in addition, $\a_*$ is the unique global minimizer of $G_{S_*}$ on $\Omega_\a$ and $\lambda \leq \min_{i \in S_*}|\phi_{S_*}^i(\a_*)|$, then exact pSTLSQ reaches a fixed point after one iteration.

For $\a_0 \neq \a_*$, a threshold that preserves the true terms at $\a_*$ may still remove one at a transient parameter value $\a_k$. As mentioned earlier, the loss of a true term cannot be reverted by Lemma \ref{lem:nested}, so support recovery requires control of the coefficients along the entire parameter trajectory. Recall that the coefficient iterates satisfy $\xi_k|_{S_{k - 1}} = \phi_{S_{k - 1}}(\a_k), \, k \geq 1$, with $\xi_0 = \phi_{\{1, \dots, p\}}(\a_0)$. Hence, preservation of the true support is equivalent to
\begin{align}\label{eq:lambda_condition_path_upper}
    \lambda \leq \inf_{k \geq 0}\min_{i \in S_*}|\xi_k^i|.
\end{align}
Provided that condition \eqref{eq:lambda_condition_path_upper} holds, all remaining spurious terms are removed at some iteration $K \geq 0$ precisely when
\begin{align}\label{eq:lambda_condition_path_lower}
    \max_{i \in \overline{S_*}}|\xi_K^i| < \lambda.
\end{align}
Thus, eventual recovery of $S_*$ is equivalent to condition \eqref{eq:lambda_condition_path_upper} together with condition \eqref{eq:lambda_condition_path_lower} for some $K \geq 0$. The following proposition gives a stronger sufficient condition that builds upon this, and that additionally assumes $\a_k \to \a_*$ and controls all supports containing $S_*$.

\begin{proposition}\label{prop:lambda_support_recovery}
Let $S_* \subseteq \{1, 2, \dots, p\}$, and consider iteration scheme \eqref{eq:iteration_scheme} with any schedule $\nu$. Suppose that $\a_k \to \a_* \in \Omega_\a$ and that $b_\a$ is continuous at $\a_*$, and let $\mathcal{R}(\a_{0}) \subseteq \Omega_{\a}$ contain all parameter iterates $\{\a_{k}\}_{k \geq 0}$. Define
\begin{align*}
    \Ss_{*} := \left\{ S \subseteq \{1, 2, \dots, p\} : S_{*} \subseteq S \right\}, \qquad \Ss_{*}^{+} := \left\{ S \subseteq \{1, 2, \dots, p\} : S_{*} \subset S \right\},
\end{align*}
and set
\begin{align*}
    L(\a) := \max_{S \in \Ss_{*}^{+}} \max_{i \in S \setminus S_{*}} |\phi_{S}^{i}(\a)|, \qquad U(\a) := \min_{S \in \Ss_{*}} \min_{i \in S_{*}} |\phi_{S}^{i}(\a)|.
\end{align*}
If
\begin{align}\label{eq:lambda_condition_region}
    L(\a_{*}) < \lambda \leq \inf_{\a \in \mathcal{R}(\a_{0})} U(\a),
\end{align}
then $S_{*} \subseteq S_{k}$ for every $k \geq 0$, and there exists $K \geq 0$ such that $S_{k} = S_{*}$ for every $k \geq K$.
\end{proposition}

\begin{proof}
We first show that no true term is removed. Since $\a_{0} \in \mathcal{R}(\a_{0})$ and $\{1, \dots, p\} \in \Ss_{*}$, we have
\begin{align*}
    |\xi_{0}^{i}| = |\phi_{\{1, \dots, p\}}^{i}(\a_{0})| \geq U(\a_{0}) \geq \inf_{\a \in \mathcal{R}(\a_{0})}U(\a) \geq \lambda
\end{align*}
for any $i \in S_{*}$. Hence, $S_{*} \subseteq S_{0}$. Now suppose that $S_{*} \subseteq S_{k - 1}$ for some $k \geq 1$. Then $S_{k - 1} \in \Ss_{*}$ and $\a_{k} \in \mathcal{R}(\a_{0})$, so
\begin{align*}
    |\xi_{k}^{i}| = |\phi_{S_{k - 1}}^{i}(\a_{k})| \geq U(\a_{k}) \geq \inf_{\a \in \mathcal{R}(\a_{0})}U(\a) \geq \lambda
\end{align*}
for any $i \in S_{*}$. Therefore, $S_* \subseteq S_k$, and thus by induction we have that $S_* \subseteq S_k$ for any $k \geq 0$.

It remains to show that all spurious terms are eventually removed. Since $\phi_{S}(\a) = A_{S}^{\dagger}b_{\a}$, continuity of $b_{\a}$ at $\a_{*}$ implies continuity of each $\phi_{S}$, and therefore of $L$, at $\a_{*}$. Hence, $L(\a_{*}) < \lambda$ implies that there exists $\eps > 0$ such that $L(\a) < \lambda$ for every $\a \in B_{\eps}(\a_{*}) \cap \Omega_\a$. Since $\a_{k} \to \a_{*}$, there exists $K \geq 1$ such that $\a_{k} \in B_{\eps}(\a_{*})$ for every $k \geq K$. For every $i \in S_{K - 1} \setminus S_*$, we have
\begin{align*}
    |\xi_K^i| = |\phi_{S_{K - 1}}^i(\a_K)| \leq L(\a_K) < \lambda.
\end{align*}
All remaining spurious terms are therefore removed when $S_K$ is formed. Since no true term is ever removed, nestedness gives $S_k = S_*$ for every $k \geq K$.
\end{proof}

\begin{remark}
For systems, $L$ and $U$ are defined analogously by also taking the maximum and minimum over the component equations and over the tuples of supports containing $S_*$.
\end{remark}

These bounds depend on the target parameter and support, so they describe conditions for recovery rather than provide a practical rule for choosing $\lambda$. The upper bound $\inf_{\a \in \mathcal{R}(\a_{0})}U(\a)$ ensures that the true terms are retained for every parameter value in $\mathcal{R}(\a_{0})$ and every support containing $S_{*}$, and the lower bound $L(\a_{*})$ ensures that, once $\a_{k}$ is sufficiently close to $\a_{*}$, all remaining spurious terms fall below the threshold. Condition \eqref{eq:lambda_condition_region} is sufficient rather than necessary because it controls every support containing $S_{*}$, including supports that may never occur during the pSTLSQ iteration. Without information about the parameter trajectory, one can apply Proposition \ref{prop:lambda_support_recovery} with $\mathcal{R}(\a_0) = \Omega_\a$. However, if it is known a priori that the iterates $\a_k$ remain in a smaller subdomain of $\Omega_\a$, one may instead take $\mathcal{R}(\a_0)$ to be this subdomain as this gives a potentially less restrictive upper bound on $\lambda$ in condition \eqref{eq:lambda_condition_region}.

\section{Numerical examples} \label{sec:numerics}

Having introduced the pSTLSQ algorithm, we now examine its performance on synthetic datasets with noise imposed. In particular, we apply inexact pSTLSQ to fractional ODEs, fractional PDEs, and nonlocal diffusion equations with parameterized kernels. We assess the recovery of the operator parameters and the correct library terms and test the learned models on initial conditions not used for fitting.

\subsection{Experimental setup}\label{sec:data}

\textbf{Data.} For the examples throughout this section, we use synthetic data $u_{\x_j}^{\ell}, \, \ell = 1, 2, \dots, n$, which is obtained by numerically solving the governing equation. Furthermore, we add Gaussian noise, so that the considered data is given by
\begin{align}\label{eq:u_noisy}
    u_{\x_j, \s}^\ell = u_{\x_j}^\ell + \s\varepsilon_{\x_j}^\ell, \qquad \varepsilon_{\x_j}^\ell \sim \mathcal{N}\bigl(0, \hat\s(u_{\x_j}^\ell)^2\bigr),
\end{align}
where $\hat\s( \cdot )$ denotes the sample standard deviation, and where $\s$ measures the noise relative to the variation of each component. To assess the consistency of recovery under noise, we draw $20$ independent realizations of $\varepsilon$ in \eqref{eq:u_noisy} and apply pSTLSQ separately to each noisy dataset, with every fit using all training initial conditions. We reuse these realizations across all nonzero noise levels, changing only the factor $\s$. In each example, we show the realization with the median parameter error (taking the smaller of the two middle values).

We use uniform grids with $400$ time points in $(0, T]$. For the ODEs, we take $T = 20$, $\Delta t = 0.05$, and for the PDEs we take $T = 2$, $\Delta t = 0.005$. Regarding the spatial domain for the PDEs we take $\TT = [0, 2\pi)$ with periodic boundary conditions, and we use a uniform grid with $n_x = 128$ grid points. Further details of the numerical methods are specified later.

\textbf{Error and test data.} We measure the fit and compute the parameter updates using the scaled least-squares error
\begin{align*}
    \overline G(\a, \Xi) := \frac{1}{2Nn}\sum_{\ell = 1}^n \|A\xi^\ell - b_{\a, \ell}\|_2^2,
\end{align*}
with $\Xi = \xi$ for scalar equations. This error averages over observations and components without including the sparsity error. We measure the parameter error by $\|\hat\a - \a_*\|_\infty$, where $\a_*$ is the true parameter and $\hat\a$ is the fitted parameter. When the true support is recovered and $\Xi_* \neq 0$, we measure the coefficient error by the relative error $\|\hat\Xi - \Xi_*\|_F/\|\Xi_*\|_F$, where $\hat{\Xi}$ and $\Xi_*$ denote the fitted and true coefficients. Finally, we test the learned models on an initial condition not used for fitting, with test data at the same noise level as the training data.

\textbf{pSTLSQ implementation.} Throughout this discussion, we take $n_\a = 1$, i.e., each iteration thresholds the coefficients, takes one projected gradient step on the scaled reduced objective $\overline G_{S_k} = G_{S_k}/(2Nn)$, and refits the coefficients on the retained support. The step size $\eta$ therefore corresponds to $\eta/(2Nn)$ in \eqref{eq:inexact_map}. We stop at the first iteration satisfying $\|\a_{k + 1} - \a_k\|_2 < 10^{-3}\|\a_1 - \a_0\|_2$, or after $2000$ iterations. For each example, we specify $\lambda$, $\eta$, and $\a_0$ and keep them fixed across noise levels and realizations. Moreover, we compute $\nabla b_\a$ analytically except for the Caputo derivative, where we use central differences with step $10^{-4}$. If one of the two points $\a \pm 10^{-4}$ lies outside $\Omega_\a$, we replace it by the nearest bound of $\Omega_\a$.

\subsection{Fractional ODEs}\label{sec:ode}

We first consider fractional ODEs involving one or more Caputo derivatives. For $T > 0$ and $f \in \Cc^1((0, T])$ with $f' \in L^1(0, T)$, the Caputo derivative of order $\a \in (0, 1)$ is defined by
\begin{align*}
    D_t^\a f(t) := \frac{1}{\Gamma(1 - \a)}\int_0^t\frac{f'(\tau)}{(t - \tau)^\a}\,\d\tau,
\end{align*}
with $t \in (0, T]$. On a uniform grid $t^k = k\Delta t$, $k = 0, \dots, n_t$, we approximate the derivative by the L1 scheme \cite{langlands2005,lin2007}, i.e.,
\begin{align*}
    (\tilde D_t^\a f)(t^k) := \frac{1}{\Gamma(2 - \a)\Delta t^\a}\sum_{s = 0}^{k - 1}\bigl((s + 1)^{1 - \a} - s^{1 - \a}\bigr)\bigl(f(t^{k - s}) - f(t^{k - s - 1})\bigr), \qquad k = 1, \dots, n_t.
\end{align*}

To generate data, we apply the L1 scheme with time step $2.5 \cdot 10^{-5}$, solve the resulting nonlinear algebraic equations at each time step by Newton's method, and keep every $2000$th value. Note that $b_\a$ depends smoothly on $\a$ in $(0, 1)^m$, as required for the results regarding pSTLSQ derived in Section \ref{sec:family}.

The three examples below involve a single fractional order, an order shared across component equations, and distinct component orders. Throughout these examples, we restrict the parameter to $\Omega_\a = [10^{-6}, 1 - 10^{-6}]^m$.

\textbf{Fractional Riccati equation.} Consider the scalar fractional Riccati equation
\begin{align*}
    D_t^\a u(t) = au(t)^2 + bu(t) + c, \qquad u(0) = \x,
\end{align*}
on $t \in [0, T]$, where $a, b, c \in \RR$. This equation is used as a simple benchmark for pSTLSQ to simultaneously recover the fractional order and the governing equation.

\begin{figure}[!tbp]
    \centering

    \begin{minipage}[c]{0.43\textwidth}
        \centering
        \includegraphics[width=\linewidth]{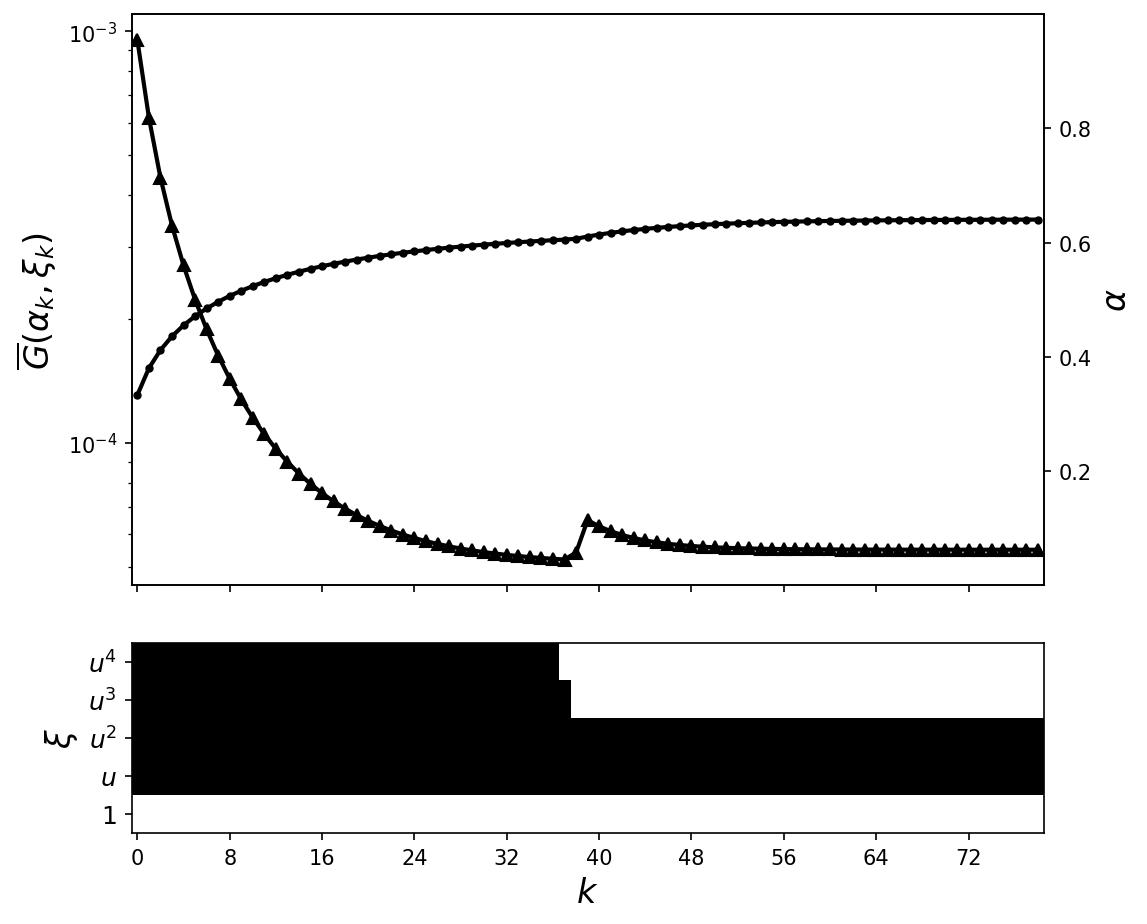}
    \end{minipage}
    \hspace{0.02\textwidth}
    \begin{minipage}[c]{0.35\textwidth}
        \centering
        \includegraphics[width=\linewidth]{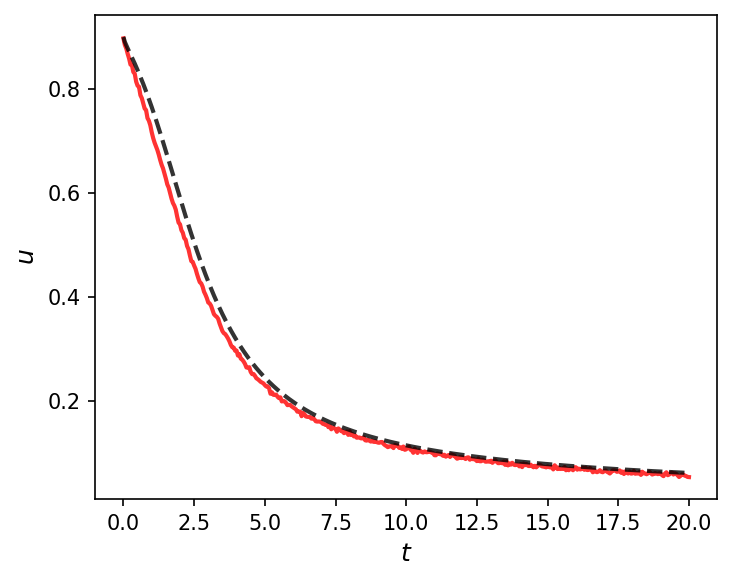}
    \end{minipage}

    \vspace{2mm}

    \begin{minipage}[c]{0.43\textwidth}
        \centering
        \includegraphics[width=\linewidth]{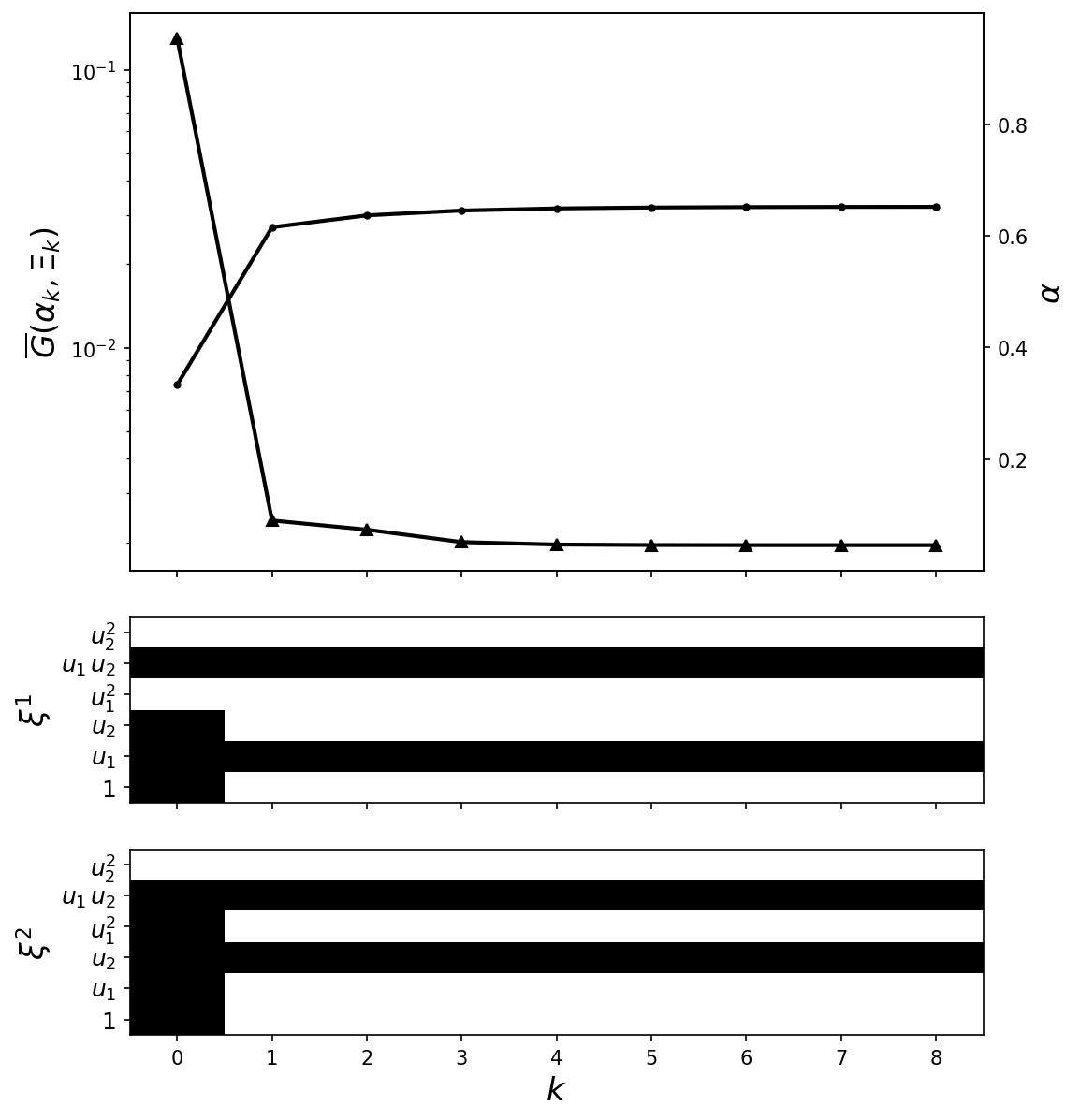}
    \end{minipage}
    \hspace{0.02\textwidth}
    \begin{minipage}[c]{0.35\textwidth}
        \centering
        \includegraphics[width=\linewidth]{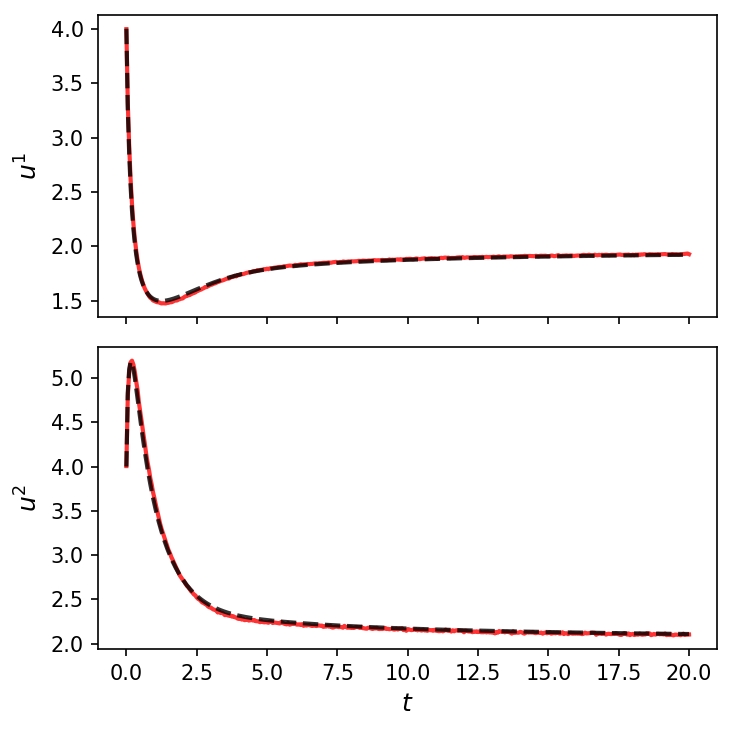}
    \end{minipage}

    \vspace{2mm}

    \begin{minipage}[c]{0.43\textwidth}
        \centering
        \includegraphics[width=\linewidth]{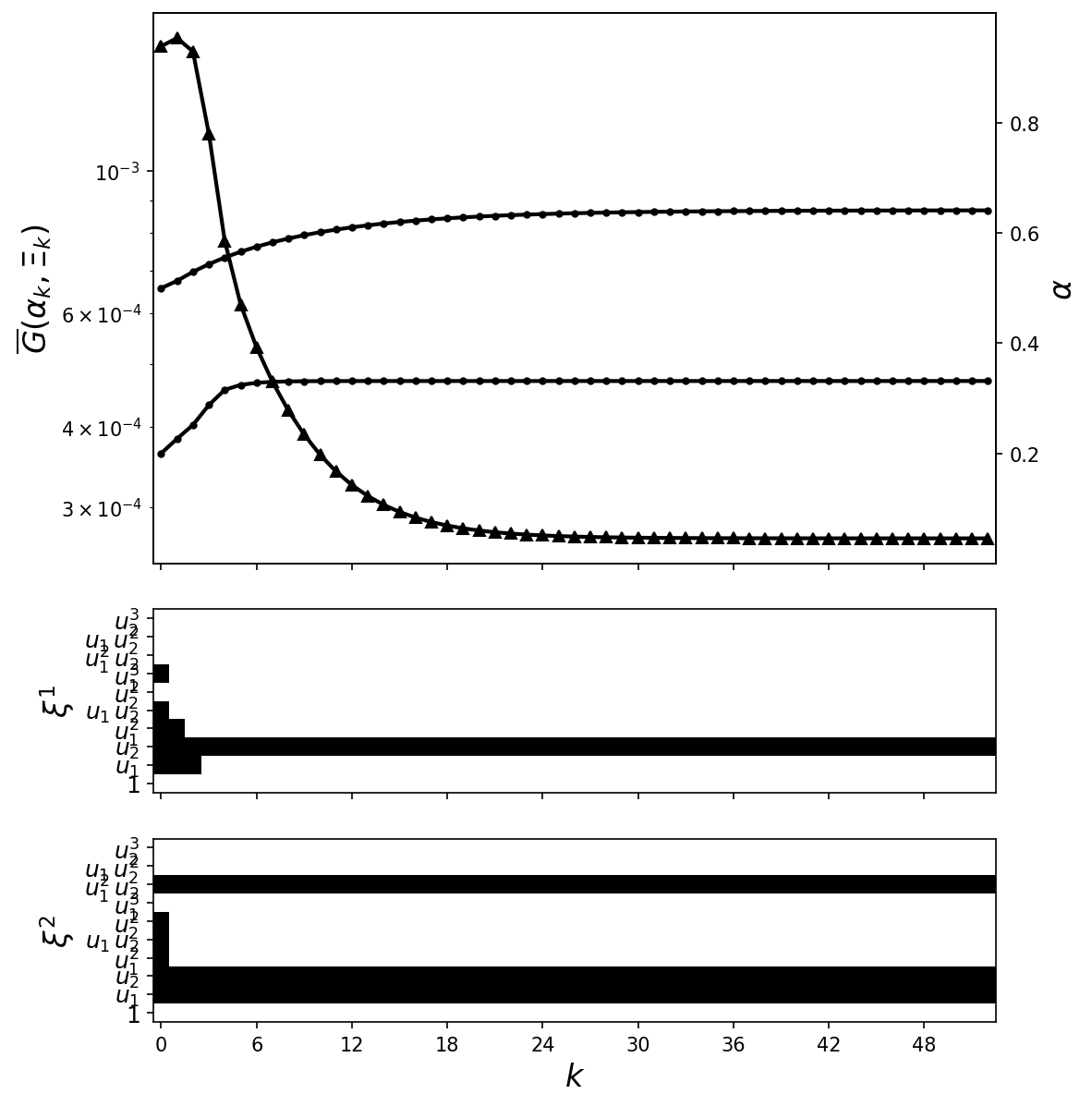}
    \end{minipage}
    \hspace{0.02\textwidth}
    \begin{minipage}[c]{0.35\textwidth}
        \centering
        \includegraphics[width=\linewidth]{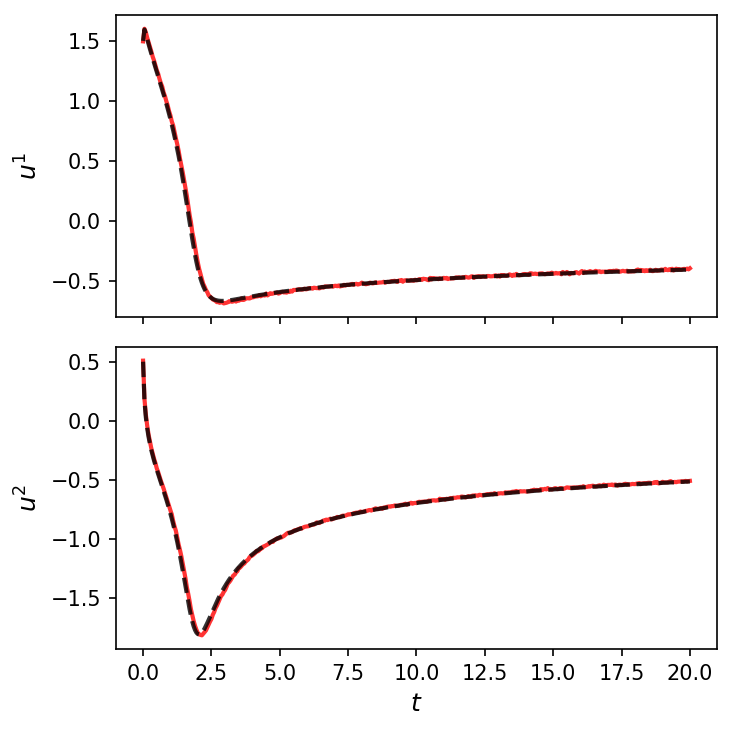}
    \end{minipage}

    \caption{Results of pSTLSQ for the fractional Riccati equation (top row), the fractional Lotka-Volterra system (middle row), and the fractional Van der Pol oscillator (bottom row). For each example, the left panel shows the active support $S_k$, the parameter iterate $\a_k$ (circles), and the scaled least-squares error $\overline{G}(\a_k, \Xi_k)$ (triangles), while the right panel compares the noisy test data (red) with the trajectory of the learned model (black dashed) for the test initial conditions $\x = 0.9$, $\x = (4, 4)$, and $\x = (\frac{3}{2}, \frac{1}{2})$, respectively.}
    \label{fig:fODEs}
\end{figure}

We set $\a_* = \frac{2}{3}$, $a = 1$, $b = -1$, and $c = 0$, and use the training initial conditions $\x = \frac{1}{5}, \frac{2}{5}, \frac{3}{5}, \frac{4}{5}$. We take a library of monomials up to degree $4$, i.e., $\{1, u, \dots, u^4\}$, set $\lambda = 2$ and $\eta = 5$, and initialize pSTLSQ with $\a_0 = \frac{1}{3}$. Under Gaussian noise with $\s = 10^{-2}$, pSTLSQ meets the stopping tolerance after $78$ iterations and recovers
\begin{align*}
    D_t^{0.640}u = -1.084\,u + 1.139\,u^2.
\end{align*}
At this noise level, all $20$ realizations recover the true support, with median parameter error $0.026$. For the test initial condition $\x = \frac{9}{10}$, the noisy test data and the trajectory of the learned model on $[0, T]$ are shown in Figure \ref{fig:fODEs}. In the absence of noise, pSTLSQ recovers the correct library terms, with parameter and coefficient errors up to $\Oo(10^{-2})$. Section \ref{sec:disc_floor} examines the source of these errors.

\textbf{Fractional Lotka-Volterra system.} Consider the fractional Lotka-Volterra system
\begin{align*}
    \begin{dcases}
        D_t^\a x = x(a - by), \\
        D_t^\a y = y(-c + dx),
    \end{dcases} \qquad u(0) = \x,
\end{align*}
on $t \in [0, T]$, where $u(t) = (x(t), y(t))$ and $a, b, c, d \in \RR$. The fractional Lotka-Volterra system arises in biological modeling when population dynamics exhibit memory effects \cite{ahmed2007equilibrium}. Variants with distinct fractional orders for the prey and predator populations have also been proposed \cite{das2011mathematical} and can likewise be identified within the pSINDy framework.

We set $\a_* = \frac{2}{3}$, $a = 1$, $b = \frac{1}{2}$, $c = \frac{3}{2}$, and $d = \frac{3}{4}$, and use the training initial conditions $\x = (\frac{1}{2}, \frac{1}{2}), (\frac{9}{2}, \frac{1}{2}), (\frac{1}{2}, \frac{9}{2}), (\frac{7}{2}, \frac{9}{2})$. We take a library of monomials in $x$ and $y$ of total degree at most $2$, including the constant, set $\lambda = 30$ and $\eta = \frac{1}{4}$, and initialize pSTLSQ with $\a_0 = \frac{1}{3}$. Under Gaussian noise with $\s = 10^{-2}$, pSTLSQ meets the stopping tolerance after $8$ iterations and recovers
\begin{align*}
    D_t^{0.652}\begin{bmatrix}x \\y\end{bmatrix} = \begin{bmatrix}x(1.032 - 0.515\,y) \\y(-1.541 + 0.771\,x)\end{bmatrix}.
\end{align*}
At this noise level, all $20$ realizations recover the true support, with median parameter error $0.014$. For the test initial condition $\x = (4, 4)$, the noisy test data and the trajectory of the learned model on $[0, T]$ are shown in Figure \ref{fig:fODEs}. In the absence of noise, pSTLSQ recovers the correct library terms, with parameter and coefficient errors up to $\Oo(10^{-2})$.

\textbf{Fractional Van der Pol oscillator.} Consider the fractional Van der Pol oscillator
\begin{align*}
    \begin{dcases}
        D_t^{\a^1}x = y, \\
        D_t^{\a^2}y = \mu(1 - x^2)y - x,
    \end{dcases} \qquad u(0) = \x,
\end{align*}
on $t \in [0, T]$, where $u(t) = (x(t), y(t))$ and $\mu > 0$. The Van der Pol oscillator is a classical nonlinear model arising in electrical circuits and self-sustained oscillatory systems \cite{vanderpol1926relaxation,guckenheimer1983nonlinear}. Fractional variants incorporate memory effects, with applications to, e.g., aeroelastic flutter \cite{juarez2023hopf}. Unlike the previous examples, the two component equations have distinct fractional orders, illustrating recovery of vector-valued operator parameters.

We set $\a_* = (\frac{1}{3}, \frac{2}{3})$ and $\mu = 1$, and use the training initial conditions $\x = (\frac{1}{2}, 0), (1, \frac{1}{2}), (2, 0), (0, 1)$. We take a library of monomials in $x$ and $y$ of total degree at most $3$, including the constant, set $\lambda = 5$ and $\eta = 1$, and initialize pSTLSQ with $\a_0 = (\frac{1}{5}, \frac{1}{2})$. Under Gaussian noise with $\s = 10^{-2}$, pSTLSQ meets the stopping tolerance after $52$ iterations and recovers
\begin{align*}
    \begin{dcases}
        D_t^{0.331}x = 1.004\,y, \\
        D_t^{0.641}y = (1.030 - 1.008\,x^2)y - 1.015\,x.
    \end{dcases}
\end{align*}
At this noise level, all $20$ realizations recover the true support, with median parameter error $0.025$. For the test initial condition $\x = (\frac{3}{2}, \frac{1}{2})$, the noisy test data and the trajectory of the learned model on $[0, T]$ are shown in Figure \ref{fig:fODEs}. In the absence of noise, pSTLSQ recovers the correct library terms, with parameter and coefficient errors up to $\Oo(10^{-2})$.

\subsection{Fractional PDEs}\label{sec:pde}

We next consider fractional PDEs involving the fractional Laplacian. For $f \in \Ss(\RR^d)$ and $\a \in (0, 1)$, this operator is defined by
\begin{align}\label{eq:frac_lap}
    (-\Delta)^\a f(x) := \Ff^{-1}\bigl(|\omega|^{2\a}(\Ff f)(\omega)\bigr)(x),
\end{align}
where $\Ss(\RR^d)$ is the Schwartz space and $\Ff$ denotes the Fourier transform.

\begin{remark}\label{rem:fraclap}
For $0 < \a < 1$, the fractional Laplacian also admits the singular-integral representation
\begin{align*}
    (-\Delta)^\a f(x) = C_{d, \a}\,\mathrm{P.V.} \int_{\RR^d} \frac{f(x) - f(y)}{|x - y|^{d + 2\a}}\,\d y,
\end{align*}
where $C_{d, \a} > 0$ is a normalization constant and $\mathrm{P.V.}$ denotes the Cauchy principal value. The Caputo derivative defined in Section \ref{sec:ode} likewise admits a representation using the Laplace transform $\Ll$,
\begin{align*}
    D_t^\a v(t) = \Ll^{-1}\bigl( s^\a(\Ll v)(s) - s^{\a - 1}v(0) \bigr)(t),
\end{align*}
for sufficiently regular functions on $[0, \infty)$ for which the transforms exist. Thus, both operators admit representations through singular kernels and integral transforms. The Caputo derivative depends on the past history, while the fractional Laplacian acts symmetrically in space. See \cite{di2012hitchhikers,kilbas2006theory} for further details. In the next subsection, we build on this kernel formulation and learn PDEs involving operators with different parameterized kernels.
\end{remark}

On the periodic domain $\TT$ specified in Section \ref{sec:data}, the fractional Laplacian is given by
\begin{align*}
    (-\Delta)^\a u(x) = \sum_{k \in \ZZ}|k|^{2\a}\widehat u_k e^{ikx},
\end{align*}
where $\a \in$ $(0, 1]$ and $\widehat u_k$ is the $k$th Fourier coefficient. We discretize this operator by multiplying the discrete Fourier transform of the sampled solution by $|k|^{2\a}$ and applying the inverse transform. For positive orders, this discretization depends smoothly on $\a$, as required for the results regarding pSTLSQ derived in Section \ref{sec:family}.

\begin{remark}
The mode $k = 0$ represents the spatial mean $\widehat u_0 = \overline u$ and has multiplier $|0|^{2\a} = 0$ for every $\a > 0$. For every $k \neq 0$, the multiplier $|k|^{2\a}$ tends to $1$ as $\a \to 0^+$. Consequently,
\begin{align*}
    \lim_{\a \to 0^+}(-\Delta)^\a u = u - \overline u.
\end{align*}
Setting $(-\Delta)^0 = I$ would instead retain the mean and therefore cause $(-\Delta)^\a u$ to be discontinuous at zero whenever $\overline u \neq 0$. Instead of changing the convention $(-\Delta)^0 = I$, we simply omit $\a = 0$ from $\Omega_\a$.
\end{remark}

We consider a fractional heat equation, a multi-order heat equation, and a space-time fractional heat equation. We take $\Omega_\a = [10^{-6}, 1]^m$ in the first two examples and $\Omega_\a = [10^{-6}, 1 - 10^{-6}]^2$ in the third. All three examples use the training initial conditions
\begin{align*}
    \begin{aligned}
        f_1(x) &= \sin(x) + \tfrac{1}{2}\sin(2x), & f_2(x) &= \cos(2x), \\
        f_3(x) &= \sin(3x), & f_4(x) &= \exp\bigl(-4(x - \pi)^2\bigr).
    \end{aligned}
\end{align*}
The test initial condition is $f(x) = \cos(4x) + \cos(6x)$. The equations at hand have analytic solutions in Fourier space, which we evaluate on a uniform spatial grid with $512$ points before restricting the data to the spatial grid specified in Section \ref{sec:data}.

The candidate library consists of all products of degree at most $2$ of $u$, $u_x$, and $u_{xx}$, including the constant, together with $u_t$ in the first two examples. Spatial derivatives are computed by second-order periodic central differences, with $u_{xx}$ obtained by applying the first-derivative approximation twice. We approximate $u_t$ by second-order central differences at interior time points and a one-sided difference at $t = T$ (where we recall that the equation is fitted only at positive times $t > 0$). The Caputo derivative in the third example is approximated by the L1 scheme of Section \ref{sec:ode}.

\textbf{Fractional heat equation.} Consider the space-fractional heat equation
\begin{align}\label{eq:fheat}
    \begin{dcases}
        u_t = -\kappa(-\Delta)^\a u, & \text{in }[0, T] \times \TT, \\
        u(0, \cdot ) = f,
    \end{dcases} \qquad \a \in (0, 1],
\end{align}
where $\kappa > 0$. For $0 < \a < 1$, this equation models anomalous diffusion driven by L\'evy flights, which allow long spatial jumps \cite{applebaum2009levy,bucur2016nonlocal}. At $\a = 1$, it reduces to the classical heat equation. We take $b_\a$ to be the discretization of $(-\Delta)^\a u$, so the true coefficient of the library term $u_t$ is $-1/\kappa$. The solution in Fourier space used to generate the data is
\begin{align}\label{eq:fheat_modes}
    \widehat u_k(t) = \exp\bigl(-\kappa|k|^{2\a}t\bigr)\widehat u_k(0).
\end{align}

\begin{figure}[!tbp]
    \centering

    \begin{minipage}[c]{0.4\textwidth}
        \centering
        \includegraphics[width=\linewidth]{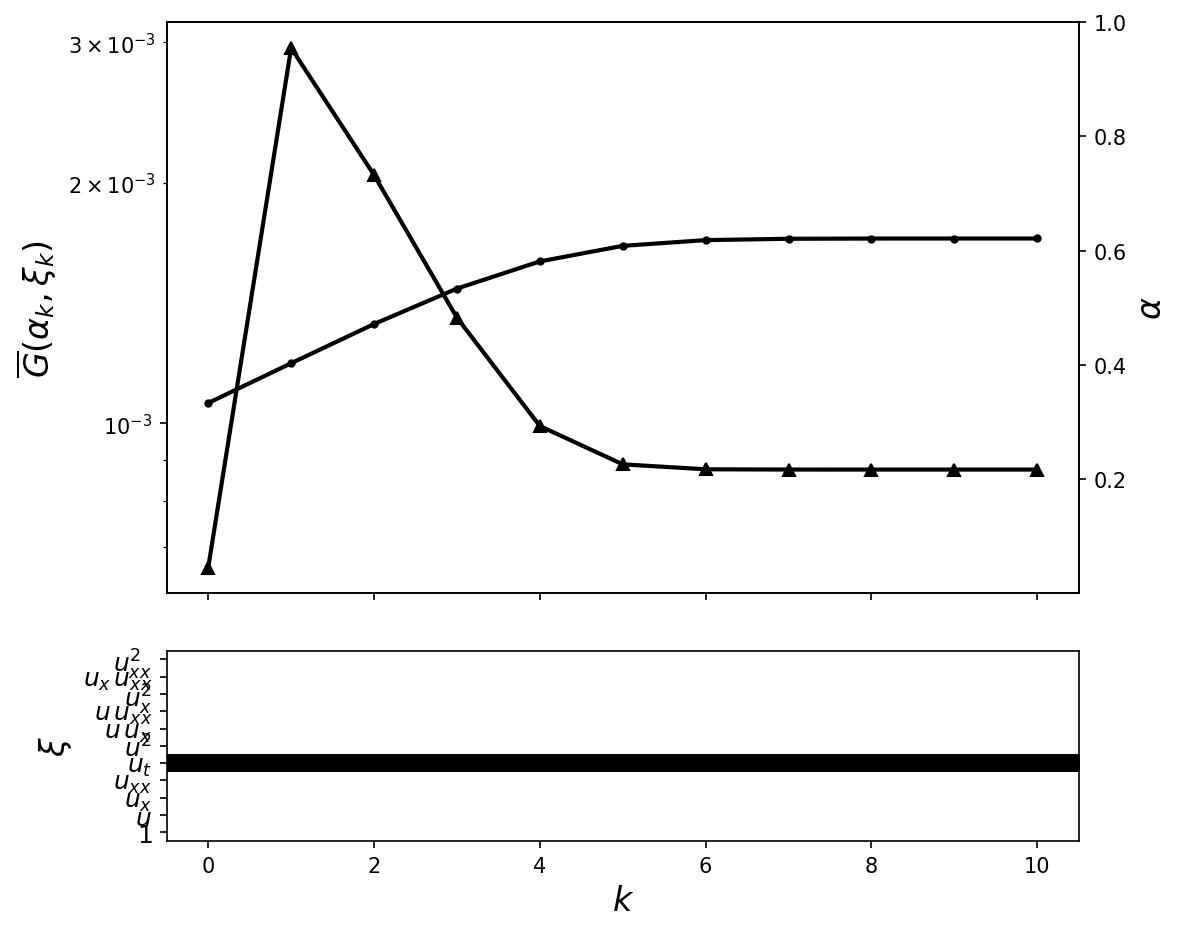}
    \end{minipage}
    \hspace{0.02\textwidth}
    \begin{minipage}[c]{0.55\textwidth}
        \centering
        \includegraphics[width=\linewidth]{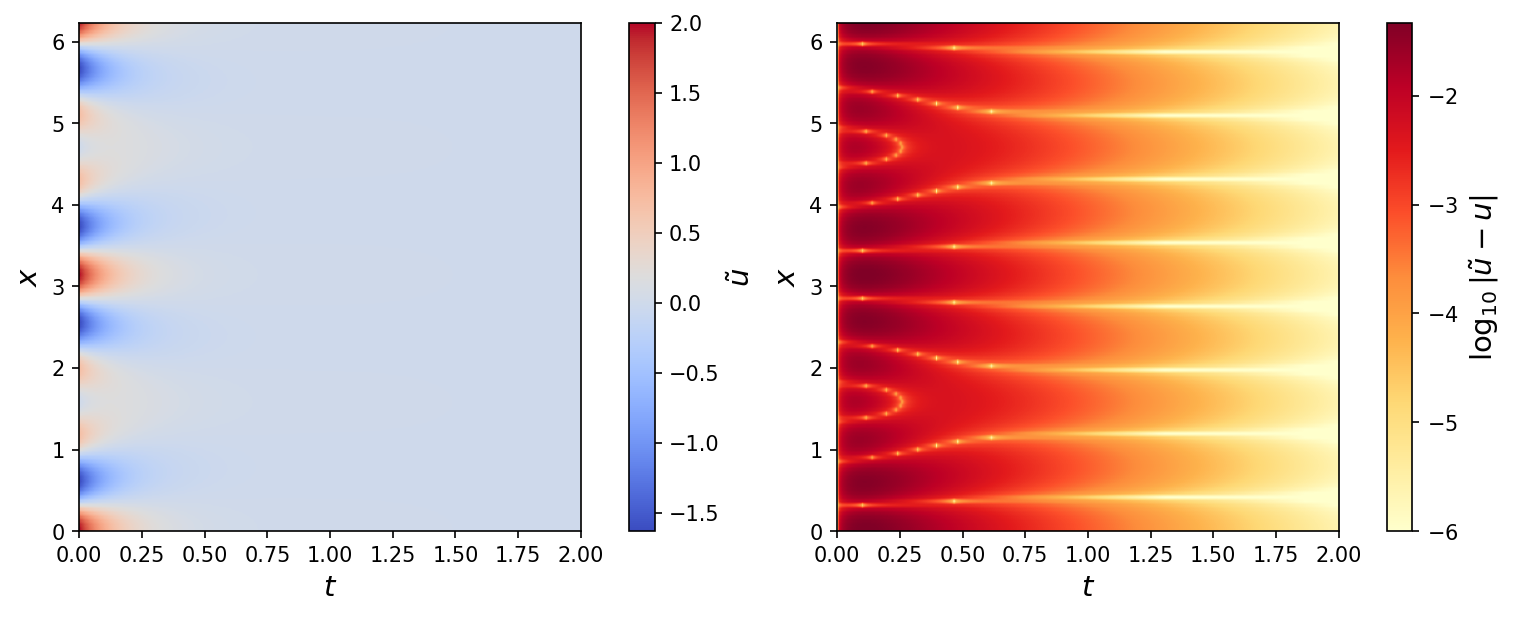}
    \end{minipage}

    \vspace{2mm}

    \begin{minipage}[c]{0.4\textwidth}
        \centering
        \includegraphics[width=\linewidth]{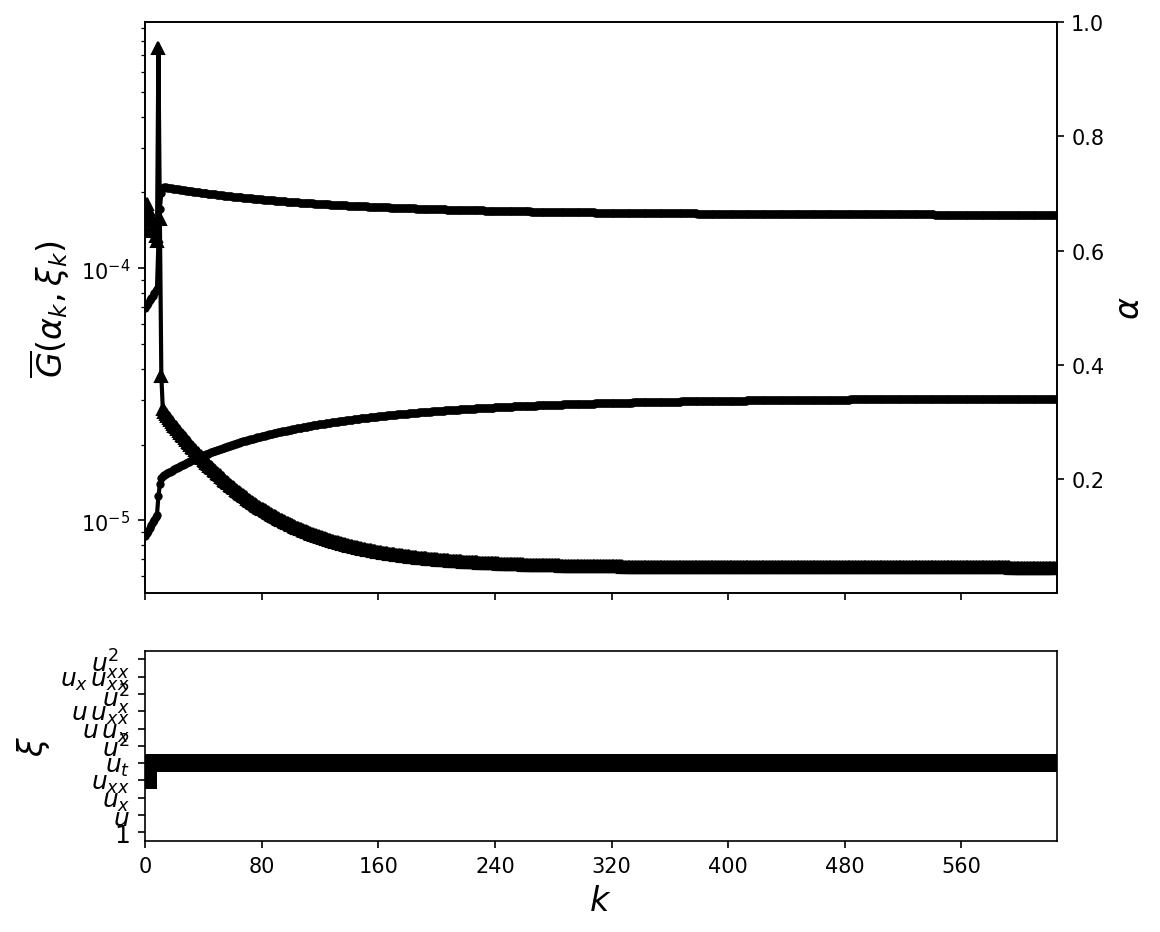}
    \end{minipage}
    \hspace{0.02\textwidth}
    \begin{minipage}[c]{0.55\textwidth}
        \centering
        \includegraphics[width=\linewidth]{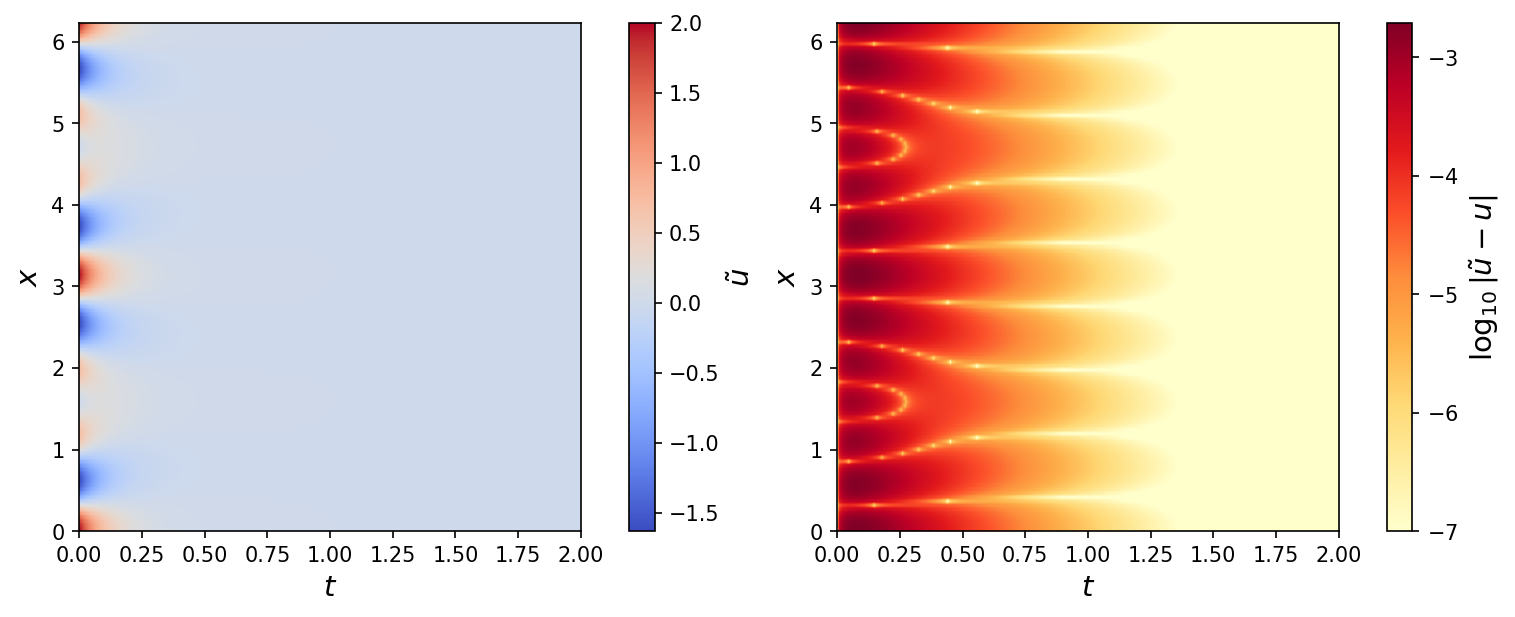}
    \end{minipage}

    \vspace{2mm}

    \begin{minipage}[c]{0.4\textwidth}
        \centering
        \includegraphics[width=\linewidth]{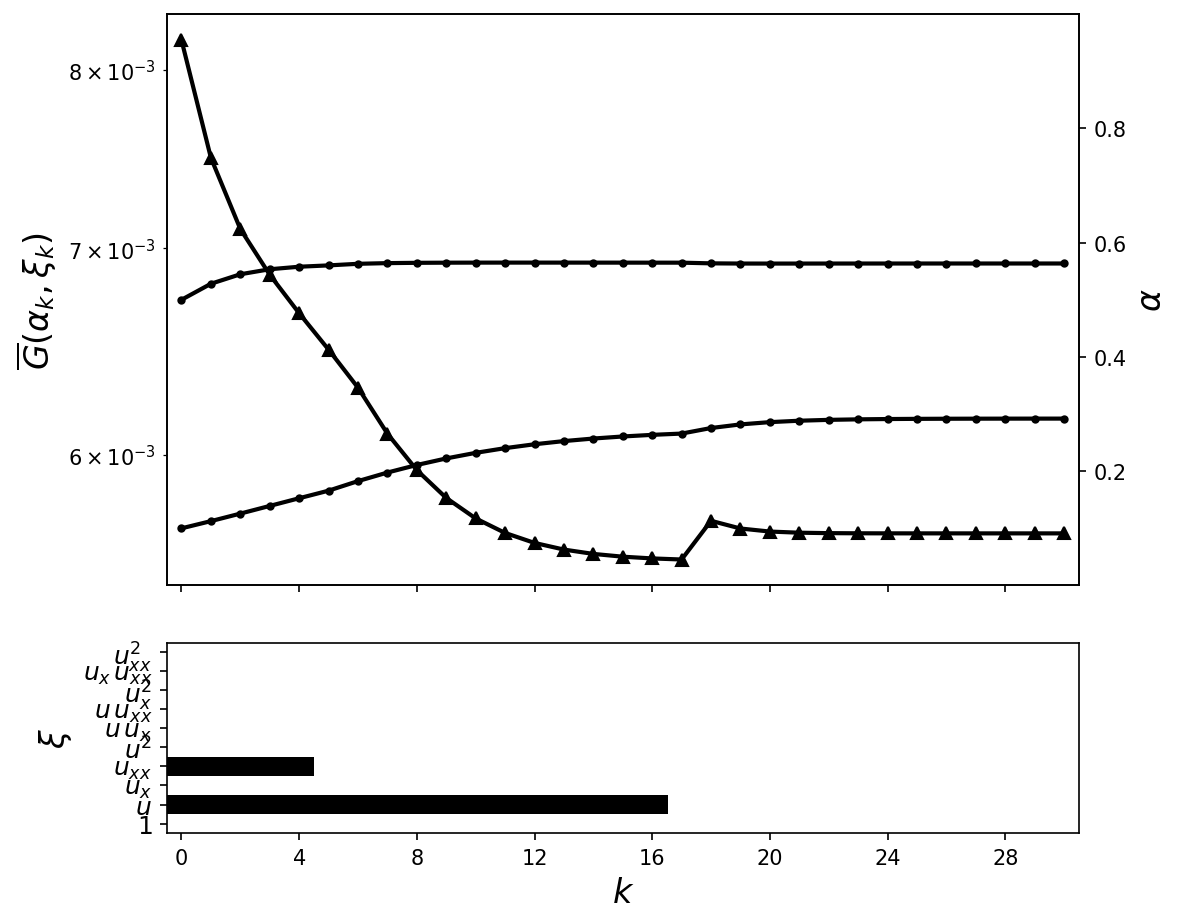}
    \end{minipage}
    \hspace{0.02\textwidth}
    \begin{minipage}[c]{0.55\textwidth}
        \centering
        \includegraphics[width=\linewidth]{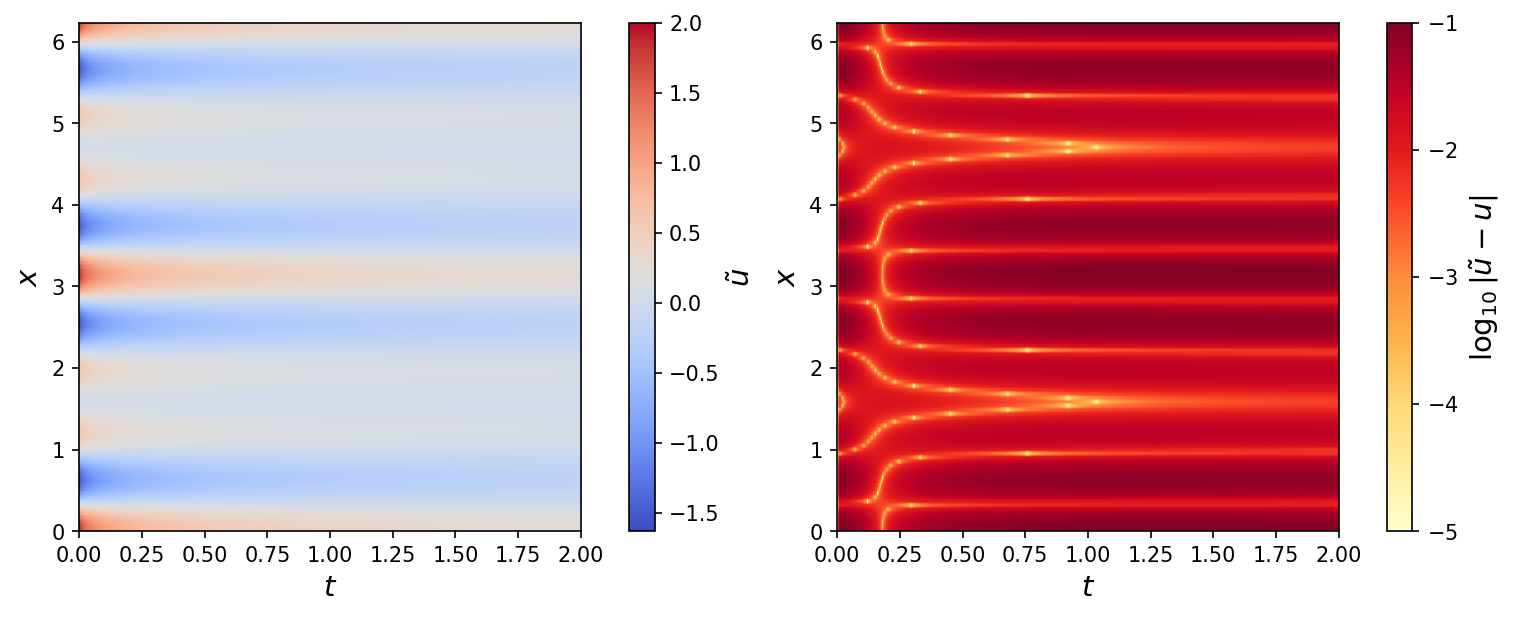}
    \end{minipage}

    \caption{Results of pSTLSQ applied to the fractional heat equation (top row), the multi-order fractional heat equation (middle row), and the space-time fractional heat equation (bottom row). For each example, the left panel shows the active support $S_k$, the parameter iterate $\a_k$ (circles), and the scaled least-squares error $\overline{G}(\a_k, \Xi_k)$ (triangles), while the right panel shows the simulation generated by the learned model and the corresponding absolute error relative to the noiseless test solution, on a logarithmic scale, for the test initial condition $f(x) = \cos(4x) + \cos(6x)$.}
    \label{fig:fPDE}
\end{figure}

We set $\a_* = \frac{2}{3}$ and $\kappa = 1$, take $\lambda = 100$ and $\eta = 5$, and initialize pSTLSQ with $\a_0 = \frac{1}{3}$. Under Gaussian noise with $\s = 10^{-3}$, pSTLSQ meets the stopping tolerance after $10$ iterations and recovers
\begin{align*}
    0.925\,u_t = -(-\Delta)^{0.621}u.
\end{align*}
At this noise level, all $20$ realizations recover the true support, with median parameter error $0.045$. For the test initial condition, the learned solution and its absolute error are shown in Figure \ref{fig:fPDE}. In the absence of noise, pSTLSQ recovers the correct library terms, with parameter and coefficient errors up to $\Oo(10^{-5})$.

\textbf{Multi-order fractional heat equation.} Consider the space-fractional heat equation with two distinct orders
\begin{align*}
    \begin{dcases}
        u_t = -\kappa_1(-\Delta)^{\a^1}u - \kappa_2(-\Delta)^{\a^2}u, & \text{in }[0, T] \times \TT, \\
        u(0, \cdot ) = f,
    \end{dcases}
\end{align*}
where $\kappa_1, \kappa_2 > 0$ and $\a^1 \neq \a^2$. For $0 < \a^1, \a^2 < 1$, this equation describes anomalous diffusion driven by two independent L\'evy motions with different power laws for the jump lengths \cite{chen2010dirichlet}. We prescribe $\kappa_1 = \kappa_2 = 1$ and take $b_\a$ to be the discretization of $[(-\Delta)^{\a^1} + (-\Delta)^{\a^2}]u$, so the true coefficient of the library term $u_t$ is $-1$\footnote{For unequal unknown weights, we can divide the equation by $\kappa_1$ and take $\Bb_\a = (-\Delta)^{\a^1} + \a^3(-\Delta)^{\a^2}$, where $\a^3$ represents the ratio $\kappa_2/\kappa_1$. The true coefficient of $u_t$ is then $-1/\kappa_1$, allowing both weights to be recovered.}. The solution in Fourier space used to generate the data is
\begin{align}\label{eq:multi_heat_modes}
    \widehat u_k(t) = \exp\bigl(-[\kappa_1|k|^{2\a^1} + \kappa_2|k|^{2\a^2}]t\bigr)\widehat u_k(0).
\end{align}

We set $\a_* = (\frac{1}{3}, \frac{2}{3})$, take $\lambda = 100$ and $\eta = 5$, and initialize pSTLSQ with $\a_0 = (\frac{1}{10}, \frac{1}{2})$. Under Gaussian noise with $\s = 10^{-4}$, pSTLSQ meets the stopping tolerance after $625$ iterations and recovers
\begin{align*}
    1.000\,u_t = -(-\Delta)^{0.341}u - (-\Delta)^{0.663}u.
\end{align*}
At this noise level, all $20$ realizations recover the true support, with median parameter error $7.7 \cdot 10^{-3}$. For the test initial condition, the learned solution and its absolute error are shown in Figure \ref{fig:fPDE}. In the absence of noise, pSTLSQ recovers the correct library terms, with parameter and coefficient errors up to $\Oo(10^{-3})$.

\textbf{Space-time fractional heat equation.} Consider the space-time fractional heat equation
\begin{align*}
    \begin{dcases}
        D_t^{\a^2}u = -\kappa(- \Delta)^{\a^1}u, & \text{in }[0, T] \times \TT, \\
        u(0, \cdot ) = f,
    \end{dcases}
\end{align*}
where $\kappa > 0$. This equation combines temporal memory with spatially nonlocal diffusion \cite{metzler2000random}. The spatial Fourier transform reduces it to scalar Caputo equations. To generate the data, we use their solutions \cite{chen2012space}
\begin{align}\label{eq:spacetime_heat_modes}
    \widehat u_k(t) = E_{\a^2}\bigl(-\kappa|k|^{2\a^1}t^{\a^2}\bigr)\widehat u_k(0),
\end{align}
where $E_{\a^2}(z) = \sum_{j = 0}^{\infty}z^j/\Gamma(1 + \a^2 j)$ is the Mittag-Leffler function. Note that for $\a^2 = 1$ we have $E_1(z) = e^z$, recovering the usual exponential evolution for the memoryless case.

We set $\kappa = 1$ and take $b_\a$ to be the discretization of $D_t^{\a^2}u + (-\Delta)^{\a^1}u$, so the true coefficient vector is $\xi_* = 0$. Moreover, we set $\a_* = (\frac{1}{3}, \frac{2}{3})$, take $\lambda = 10$ and $\eta = 1$, and initialize pSTLSQ with $\a_0 = (\frac{1}{10}, \frac{1}{2})$. Under Gaussian noise with $\s = 10^{-2}$, pSTLSQ meets the stopping tolerance after $30$ iterations and recovers
\begin{align*}
    D_t^{0.564}u = -(-\Delta)^{0.292}u.
\end{align*}
At this noise level, all $20$ realizations recover the empty support, with median parameter error $0.10$. For the test initial condition, the learned solution and its absolute error are shown in Figure \ref{fig:fPDE}. In the absence of noise, pSTLSQ also recovers the empty support, with parameter error up to $\Oo(10^{-2})$.

\subsection{Nonlocal diffusion with parameterized kernels}\label{sec:kernel}

We next consider nonlocal diffusion equations of the form
\begin{align*}
    \begin{dcases}
        u_t = \kappa\mathcal{K}_\a u, & \text{in }[0, T] \times \TT, \\
        u(0, \cdot ) = f,
    \end{dcases}
\end{align*}
where $\kappa > 0$ and
\begin{align*}
    (\mathcal{K}_\a u)(x) := \int_\TT K_\a\bigl(d_\TT(x, y)\bigr)\bigl(u(y) - u(x)\bigr)\,\d y.
\end{align*}
Here, $\a$ is the unknown kernel parameter, and $d_\TT(x, y) := \min\{|x - y|, 2\pi - |x - y|\}$ denotes the periodic distance between points $x, y \in [0, 2\pi)$. The kernel weights the differences $u(y) - u(x)$, so larger surrounding values increase $u(x)$ and smaller values decrease it. Operators of this form arise in nonlocal diffusion and peridynamics \cite{du2012analysis,andreu2010nonlocal}.

In the examples below, we consider a regularized power-law kernel and a Gaussian kernel. Moreover, we take $\kappa = 1$, use the training and test initial conditions of Section \ref{sec:pde}, and we generate the data on a uniform spatial grid with $512$ points and restrict them to the spatial grid specified in Section \ref{sec:data}. The candidate library is $\{1, u, u_x, u_{xx}, u_t\}$ for the power-law kernel example, and we omit $u_{xx}$ for the Gaussian kernel for reasons explained below. In the regression, we approximate the integral by the rectangle rule on a uniform grid $\{x^i\}$ with spacing $\Delta x$,
\begin{align*}
    (\mathcal{K}_\a u)(x^i) \approx \Delta x\sum_{j \neq i} K_\a\bigl(d_\TT(x^i, x^j)\bigr)\bigl(u(x^j) - u(x^i)\bigr),
\end{align*}
where the term $j = i$ is zero and therefore omitted. Since the kernel weights depend only on $x^i - x^j$ through the periodic distance, the weighted sum of the $u(x^j)$ values is a convolution, and $(\mathcal{K}_\a u)(x^i)$ can therefore be naturally computed using the discrete Fourier transform. Finally, since all $K_\a$ considered in the examples below are smooth in $\a$, $b_\a$ is also, as required for the results regarding pSTLSQ derived in Section \ref{sec:family}.

\textbf{Power-law kernel.} We take
\begin{align*}
    K_\a(r) = \frac{1}{r^{1 + 2\a} + \delta},
\end{align*}
where $\a \in \Omega_\a = [10^{-3}, 1 - 10^{-3}]$ and where $\delta > 0$ is fixed and keeps the kernel bounded at $r = 0$, while $\a$ controls its decay with distance. We generate the data by evolving each Fourier mode exactly in time using
\begin{align*}
    \widehat u_k(t) = \exp\left(-2\kappa t\int_0^\pi K_\a(r)\bigl(1 - \cos(kr)\bigr)\,\d r\right) \widehat u_k(0),
\end{align*}
for the $512$ modes $k$ of the finer spatial grid. The integral is evaluated by adaptive quadrature to obtain higher accuracy than the rectangle rule used in the regression.

\begin{figure}[!tbp]
    \centering

    \begin{minipage}[c]{0.4\textwidth}
        \centering
        \includegraphics[width=\linewidth]{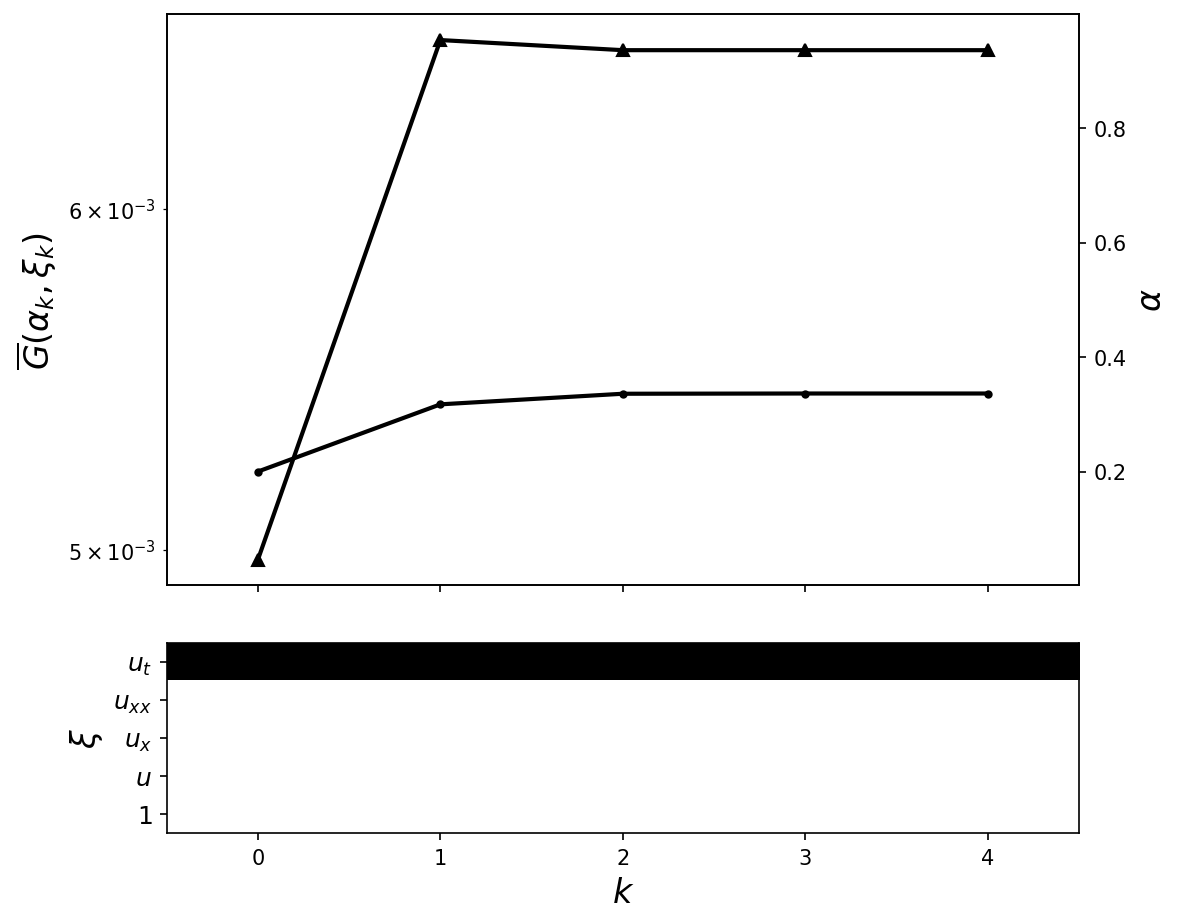}
    \end{minipage}
    \hspace{0.02\textwidth}
    \begin{minipage}[c]{0.55\textwidth}
        \centering
        \includegraphics[width=\linewidth]{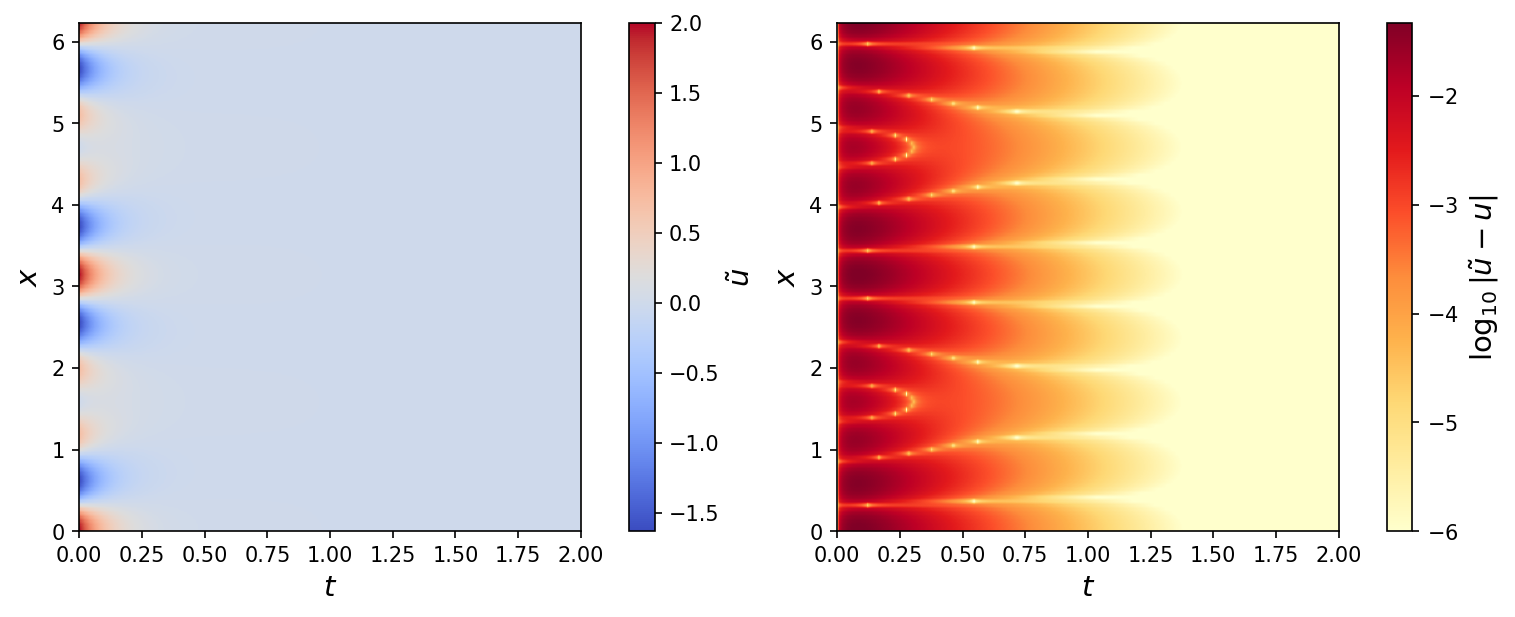}
    \end{minipage}

    \vspace{2mm}

    \begin{minipage}[c]{0.4\textwidth}
        \centering
        \includegraphics[width=\linewidth]{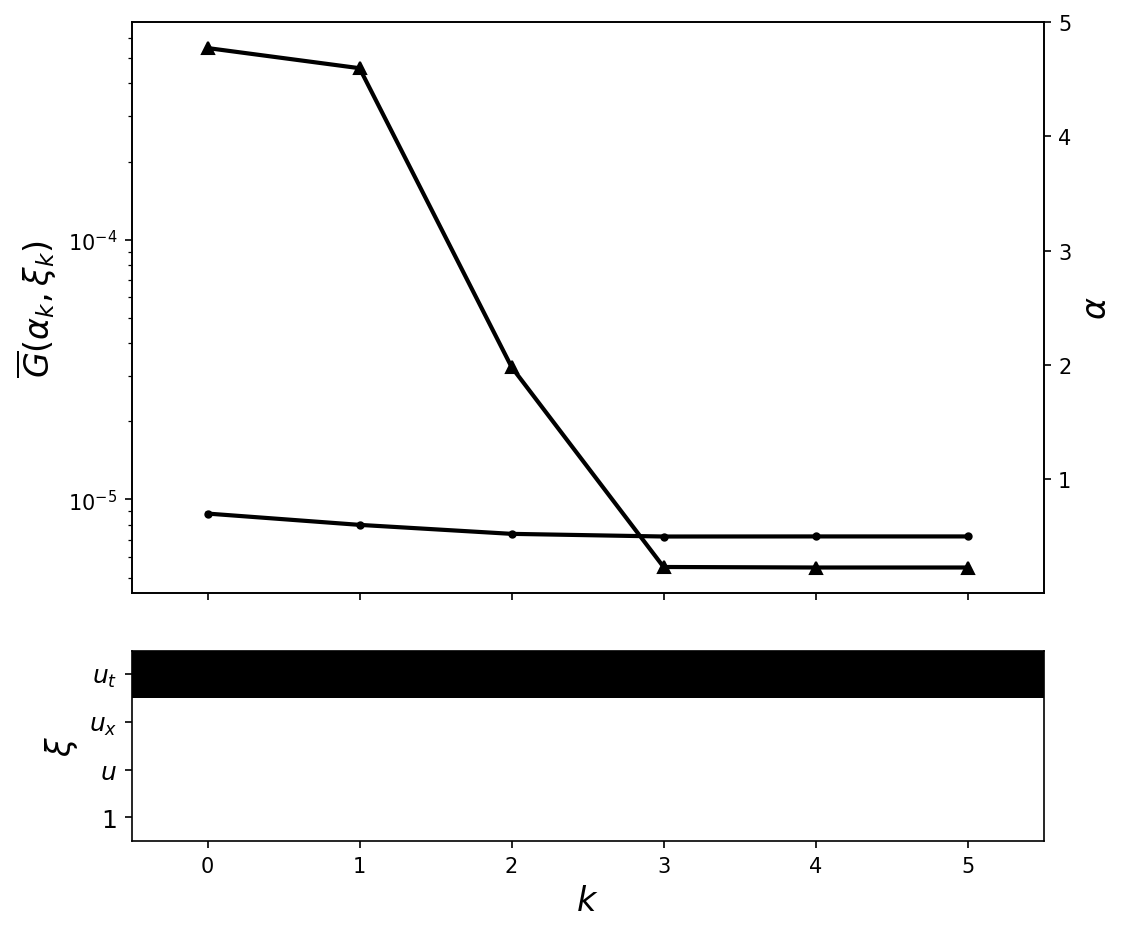}
    \end{minipage}
    \hspace{0.02\textwidth}
    \begin{minipage}[c]{0.55\textwidth}
        \centering
        \includegraphics[width=\linewidth]{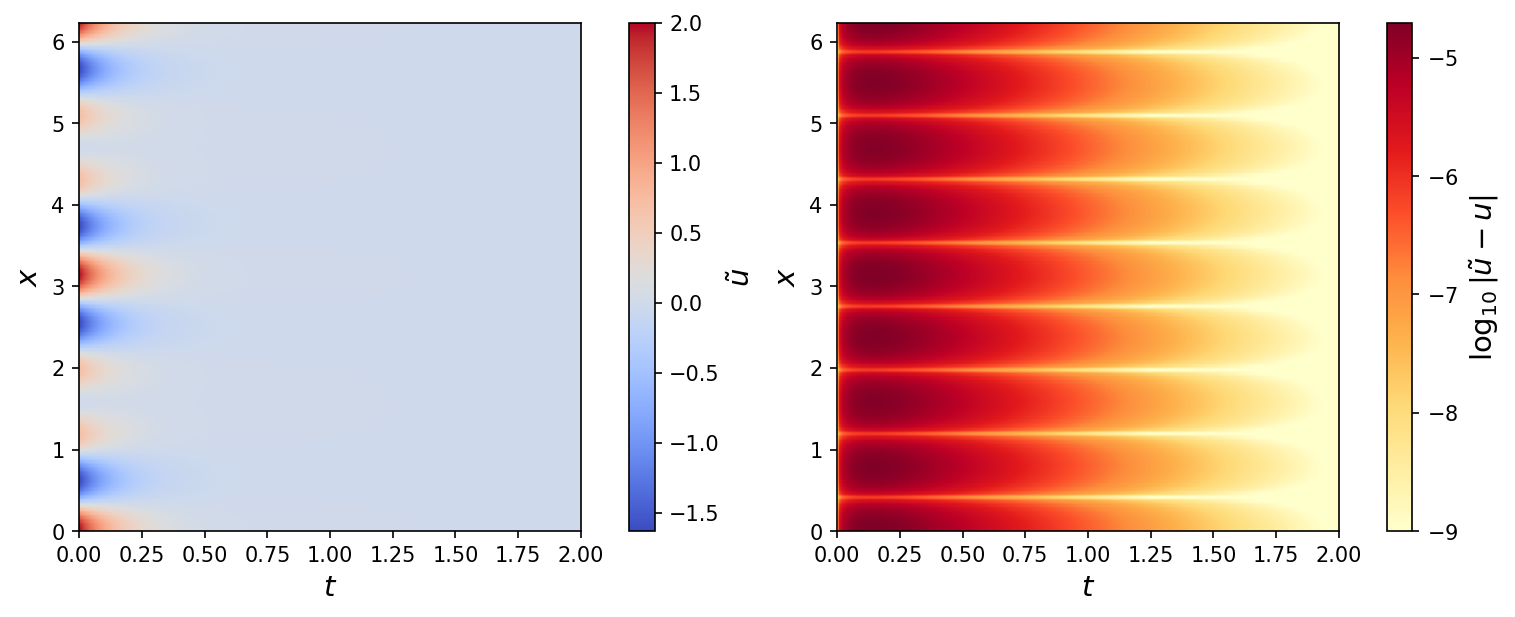}
    \end{minipage}

    \caption{Results of pSTLSQ applied to the nonlocal diffusion equation with a power-law kernel (top row) and a Gaussian kernel (bottom row). For each example, the left panel shows the active support $S_k$, the parameter iterate $\a_k$ (circles), and the scaled least-squares error $\overline{G}(\a_k, \Xi_k)$ (triangles), while the right panel shows the solution generated by the learned model and its absolute error relative to the noiseless test solution, on a logarithmic scale, for the test initial condition $f(x) = \cos(4x) + \cos(6x)$.}
    \label{fig:kernels}
\end{figure}

We set $\a_* = \frac{2}{5}$ and $\delta = 10^{-8}$, take $\lambda = 100$ and $\eta = 5$, and initialize pSTLSQ with $\a_0 = \frac{1}{5}$. Under Gaussian noise with $\s = 5 \cdot 10^{-3}$, pSTLSQ meets the stopping tolerance after $4$ iterations and recovers
\begin{align*}
    0.928\,u_t = \mathcal{K}_{0.336}u.
\end{align*}
At this noise level, all $20$ realizations recover the true support, with median parameter error $0.064$. For the test initial condition, the learned solution and its absolute error are shown in Figure \ref{fig:kernels}. In the absence of noise, pSTLSQ recovers the correct library terms, with parameter and coefficient errors up to $\Oo(10^{-3})$.

\textbf{Gaussian kernel.} We take
\begin{align*}
    K_\a(r) = c_\a\exp\left(-\frac{r^2}{2\a^2}\right),
\end{align*}
where $\a \in \Omega_\a = [10^{-3}, 5]$ controls the width of the kernel. With $c_\a$ fixed, the operator would tend to zero as $\a \to 0^+$, allowing a small least-squares error with $\xi = 0$. To prevent this, we follow the normalization for local diffusion in \cite[Section~5.2]{du2012analysis} and choose
\begin{align*}
    c_\a = \frac{2}{\Delta x\sum_j d_\TT^2(x^j, 0) \exp(-d_\TT^2(x^j, 0)/(2\a^2))}.
\end{align*}
When $\a$ is much smaller than $\Delta x$, this choice of $c_\a$ makes the discrete operator for $(\mathcal{K}_\a u)(x^i)$ approach the usual three-point approximation of $u_{xx}$. Hence, we omit $u_{xx}$ from the library to avoid selecting this term at an incorrect value of $\a$. Section \ref{sec:a0} shows this effect for the fractional heat equation, where $u_{xx}$ is kept in the library. Data for this example are moreover generated by the rectangle rule in space and a fourth-order Runge-Kutta scheme with time step $10^{-3}$ on the finer spatial grid.

We set $\a_* = \frac{1}{2}$, take $\lambda = 80$ and $\eta = 10$, and initialize pSTLSQ with $\a_0 = \frac{7}{10}$. Under Gaussian noise with $\s = 10^{-4}$, pSTLSQ meets the stopping tolerance after $5$ iterations and recovers
\begin{align*}
    1.000\,u_t = \mathcal{K}_{0.500}u.
\end{align*}
At this noise level, all $20$ realizations recover the true support, with median parameter error $4.7 \cdot 10^{-5}$. For the test initial condition, the learned solution and its absolute error are shown in Figure \ref{fig:kernels}. In the absence of noise, pSTLSQ recovers the correct library terms, with parameter and coefficient errors up to $\Oo(10^{-5})$.

\section{Numerical sensitivity analysis}\label{sec:sensitivity}

In this section, we examine how discretization, the stopping criterion, and the choices of threshold $\lambda$, initial parameter $\a_0$, and learning rate $\eta$ affect recovery of the governing equations. We also study the effects of noise and smoothing, the interaction between parameter and support updates, and the influence of symmetric parameterizations. In particular, we focus on recovering the true parameter $\a_*$ and the true support $S_*$, since the coefficients are obtained by least squares once these are fixed. Throughout, $\a_*$ and $S_*$ denote the true parameter and support, respectively, and $\hat\a$ denotes the parameter returned by pSTLSQ. We use different examples from Section \ref{sec:numerics} to illustrate the effects discussed below.

\subsection{Discretization error}\label{sec:disc_floor}

In the absence of noise, pSTLSQ recovers the true support in all examples of Section \ref{sec:numerics}, but small parameter errors remain. This error could be the result of either the error from approximating the data and operators (the discretization error) or the error from stopping the iterations before reaching a minimizer of the reduced objective (the optimization error). To see this, first note that the exact solution satisfies the equation at $\a_*$ with the true coefficients on $S_*$. The continuous least-squares error is therefore zero at $\a_*$, i.e., $\a_*$ is a minimizer of the continuous least-squares problem. However, the numerical data and the discretized operators generally do not satisfy $b_{\a_*} = A\xi_*$ exactly, and as a consequence $\a_*$ might not minimize the discretized objective $\overline G_{S_*}$. If $\a_*$ is a nondegenerate interior minimizer of the pSINDy objective, though, sufficiently small discretization errors give a nearby local minimizer $\a^{\dagger}$ in a way similar to Theorem \ref{thm:noise}. Moreover, the parameter $\hat\a$ returned by pSTLSQ may differ from $\a^{\dagger}$ due to the optimization error. In particular, we can decompose the parameter error as
\begin{align*}
    \hat\a - \a_* = (\a^\dagger - \a_*) + (\hat\a - \a^\dagger),
\end{align*}
corresponding to respectively the discretization and optimization error. Recall that for a small enough discretization error and for $\a_*$ a nondegenerate interior minimizer, $\a^{\dagger}$ is an interior minimizer close to $\a_*$, and linearizing the first-order condition $\nabla\overline G_{S_*}(\a^\dagger) = 0$ around $\a_*$ thus gives
\begin{align*}
    \a^\dagger - \a_* \approx -\bigl[\nabla^2\overline G_{S_*}(\a_*)\bigr]^{-1} \nabla\overline G_{S_*}(\a_*),
\end{align*}
provided the Hessian at $\a_*$ is invertible. Table \ref{tab:disc_floor} compares this linear estimate with the parameter errors returned by pSTLSQ on noiseless data. Except for the multi-order heat equation, the two agree to within $10\%$, suggesting that discretization is the main source of the observed parameter errors in these examples.

\begin{table}[htbp]
    \centering
    \small
    \setlength{\tabcolsep}{4pt}
    \begin{tabular}{@{}lccc@{}}
        \toprule
        & linear estimate of $\a^\dagger - \a_*$ & $\hat\a - \a_*$ & $\operatorname{eig}(\nabla^2\overline G_{S_*}(\a_*))$ \\
        \midrule
        Riccati & $-8.66 \cdot 10^{-3}$ & $-8.88 \cdot 10^{-3}$ & $1.45 \cdot 10^{-2}$ \\
        Lotka-Volterra & $-1.12 \cdot 10^{-2}$ & $-1.13 \cdot 10^{-2}$ & $2.13$ \\
        Van der Pol & $(-1.60, -15.8) \cdot 10^{-3}$ & $(-1.59, -16.0) \cdot 10^{-3}$ & $0.108$, $0.508$ \\
        \midrule
        Fractional heat & $4.32 \cdot 10^{-5}$ & $4.03 \cdot 10^{-5}$ & $0.143$ \\
        Multi-order heat & $(-3.99, 2.98) \cdot 10^{-4}$ & $(-11.7, 6.40) \cdot 10^{-4}$ & $1.55 \cdot 10^{-3}$, $0.147$ \\
        Space-time heat & $(3.33, -10.4) \cdot 10^{-3}$ & $(3.18, -10.7) \cdot 10^{-3}$ & $0.361$, $0.557$ \\
        \midrule
        Power-law kernel & $5.31 \cdot 10^{-3}$ & $5.24 \cdot 10^{-3}$ & $0.154$ \\
        Gaussian kernel & $-8.02 \cdot 10^{-5}$ & $-8.77 \cdot 10^{-5}$ & $0.110$ \\
        \bottomrule
    \end{tabular}
    \caption{Discretization diagnostics for the noiseless fits, where the first column gives the linear estimate of $\a^\dagger - \a_*$, the second the parameter error returned by pSTLSQ, and the last column the eigenvalues of $\nabla^2\overline G_{S_*}(\a_*)$.}
    \label{tab:disc_floor}
\end{table}

The last column of Table \ref{tab:disc_floor} reports the eigenvalues of the Hessian $\nabla^2\overline G_{S_*}(\a_*)$, which describe the local curvature of the reduced objective around $\a_*$. Small eigenvalues indicate directions of weak curvature, in which discretization errors can cause larger parameter displacements and gradient descent may converge slowly. For the multi-order heat equation, the eigenvalues differ by a factor $\Oo(10^2)$, and pSTLSQ stops before reaching the minimizer along the direction of smaller curvature. Hence, the optimization error contributes significantly to the parameter error in this example and helps explain the larger difference between the linear estimate and the parameter error returned by pSTLSQ. With the stricter stopping tolerance $10^{-6}$ instead of $10^{-3}$, pSTLSQ reaches parameter error $3.9 \cdot 10^{-4}$ after $1480$ iterations, in agreement with the linear estimate.

\subsection{Threshold selection}\label{sec:thresh_num}

The threshold $\lambda$ determines which library terms are retained during pSTLSQ. To recover the true support $S_*$, it must be small enough to preserve the true terms while allowing spurious terms to be removed. Proposition \ref{prop:lambda_support_recovery} gives a sufficient condition for this. Suppose that the parameter iterates converge to some $\a_\infty$ (which might not equal the true parameter $\a_*$ due to discretization errors). Then the condition is
\begin{align*}
    L(\a_\infty) < \lambda \leq \inf_{\a \in \mathcal{R}(\a_0)}U(\a),
\end{align*}
where $\mathcal{R}(\a_0)$ contains the entire parameter trajectory. We recall from the discussion in Section \ref{sec:threshold} that the lower bound ensures that spurious terms are eventually removed, while the upper bound prevents true terms from being removed. Without information about the parameter trajectory, one can take $\mathcal{R}(\a_0) = \Omega_\a$, although this may give an unnecessarily restrictive condition.

In our numerical examples, the true parameter $\a_*$ is known. For scalar parameters, we consider the reference interval
\begin{align*}
    I_{\a_*}(\a_0) := [\min\{\a_0, \a_*\}, \max\{\a_0, \a_*\}],
\end{align*}
motivated by parameter iterates moving from $\a_0$ toward $\a_*$ without crossing it. Using $I_{\a_*}(\a_0)$ in place of $\mathcal{R}(\a_0)$, we examine which values of $\lambda$ satisfy
\begin{align}\label{eq:lambda_condition_reference}
    L(\a_*) < \lambda \leq \inf_{\a \in I_{\a_*}(\a_0)}U(\a).
\end{align}
This serves as a numerical diagnostic rather than a direct application of Proposition \ref{prop:lambda_support_recovery}, since convergence to $\a_*$ is not guaranteed and the iterates may leave $I_{\a_*}(\a_0)$. Moreover, an empty interval does not imply that no thresholding parameter $\lambda > 0$ exists that recovers the true support as the bounds also consider supports that may not occur during the iterations.

Computing $L$ and $U$ requires considering every support containing $S_*$. For a simpler diagnostic, we consider only the full library, whose least-squares coefficients are given by $\phi(\a)$. We define
\begin{align*}
    \widetilde U(\a) := \min_{i \in S_*}|\phi^i(\a)|, \qquad \widetilde L(\a) := \max_{i \in \overline{S_*}}|\phi^i(\a)|.
\end{align*}
Since the full library is one of the supports considered in $L$ and $U$, we have $L(\a) \geq \widetilde L(\a)$ and $U(\a) \leq \widetilde U(\a)$. Replacing $L$ and $U$ by these quantities gives a weaker condition on $\lambda$ that can be evaluated using only the full-library coefficients.

Finally, we also examine the coefficients along the actual pSTLSQ trajectory. At iteration $k$, before forming the new support $S_k$, define
\begin{align*}
    U_k := \min_{i \in S_*}|\xi_k^i|, \qquad L_k := \max_{i \in S_{k - 1} \setminus S_*}|\xi_k^i|,
\end{align*}
where $S_{-1} = \{1, \dots, p\}$ and $L_k = 0$ if no spurious terms remain. Since the initial fit uses the full library, we have $U_0 = \widetilde U(\a_0)$ and $L_0 = \widetilde L(\a_0)$. If $\lambda \leq U_k$ at every iteration, no true term is removed. Once $L_k < \lambda$, all remaining spurious terms are removed at that thresholding step. These quantities therefore allow us to examine directly whether the chosen threshold separates true and spurious terms along the computed trajectory.

\begin{figure}[!tbp]
    \centering
    \begin{subfigure}[t]{0.49\textwidth}
        \centering
        \includegraphics[width=\linewidth]{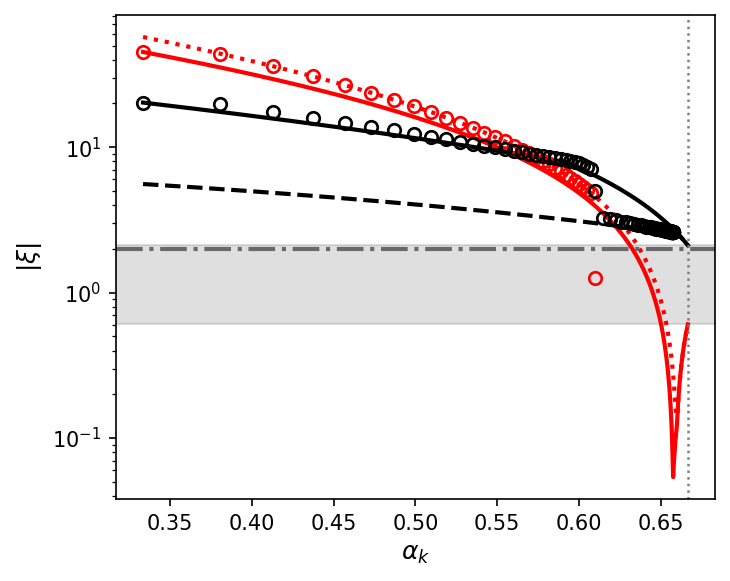}
    \end{subfigure}
    \hfill
    \begin{subfigure}[t]{0.49\textwidth}
        \centering
        \includegraphics[width=\linewidth]{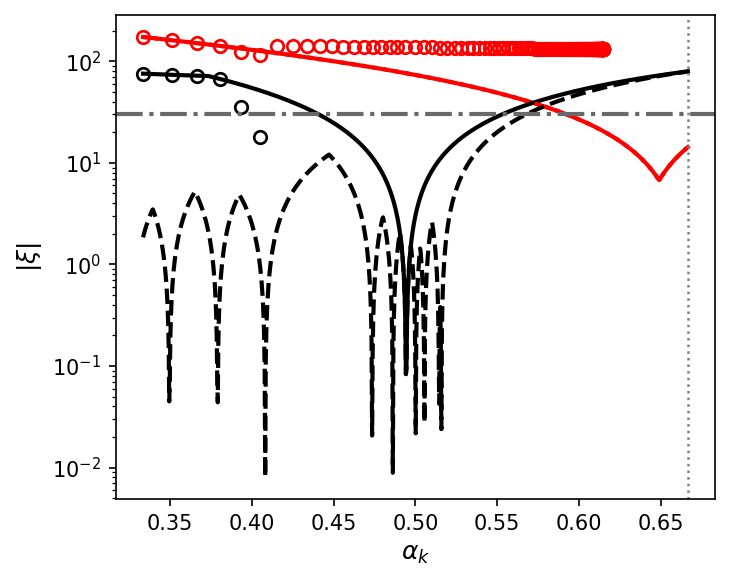}
    \end{subfigure}
    \caption{Threshold diagnostics for the fractional Riccati equation (left) and the fractional Lotka-Volterra system with initial conditions $\x = (1, 1), (2, 2), (1, 2), (2, 1)$ (right). In both examples, $\a_0 = 1/3$, $\a_* = 2/3$, and $\a \in I_{\a_*}(\a_0)$. The black dashed and solid curves show $U$ and $\widetilde U$, respectively, while the red dotted and solid curves show $L$ and $\widetilde L$. The black and red circles show the trajectory-dependent quantities $U_k$ and $L_k$. The horizontal dash-dotted line indicates the selected threshold $\lambda$, and the vertical dotted line indicates $\a_*$. The gray region indicates the interval in condition \eqref{eq:lambda_condition_reference}, bounded below by $L(\a_*)$ and above by the minimum of $U$ over $I_{\a_*}(\a_0)$. For the Lotka-Volterra initial conditions shown, this interval is empty because $U$ vanishes within the reference interval.}
    \label{fig:lambda_condition_Ia}
\end{figure}

For the fractional Riccati equation, we take $\a_0 = 1/3$ and $\a_* = 2/3$. The bounds in condition \eqref{eq:lambda_condition_reference} give the interval $(0.605, 2.11]$, which contains the selected value $\lambda = 2$ as shown in Figure \ref{fig:lambda_condition_Ia}. Along the computed trajectory, the coefficients of the true terms remain above the threshold, while the spurious terms are eventually removed. At $\lambda = 2$, the true model is a local but not the global minimizer of $F$. The model containing only $u$ at $\a \approx 0.80$ has a smaller objective value, and the true model is the global minimizer only for $\lambda < 0.71$ on the noiseless data.

The Lotka-Volterra system illustrates how threshold selection depends on the training data and the candidate library. The right panel of Figure \ref{fig:lambda_condition_Ia} uses the initial conditions $\x = (1, 1), (2, 2), (1, 2), (2, 1)$ instead of those in Section \ref{sec:ode}, where we note that the initial condition $\x = (2, 2)$ is an equilibrium (and thus its constant trajectory provides no information about the fractional order). For these data, $U$ vanishes at several points in $I_{\a_*}(\a_0)$, making the interval in condition \eqref{eq:lambda_condition_reference} empty. With the remaining settings of Section \ref{sec:ode}, pSTLSQ removes the true term $x$ and retains the spurious terms $y$ and $y^2$ in the first equation, while recovering the support of the second equation. The iteration nevertheless meets the stopping criterion at $\hat\a = 0.615$. Using the training initial conditions of Section \ref{sec:ode} instead removes the zeros of $U$ and gives the interval $(6.55, 35.5]$, which contains the selected threshold $\lambda = 30$. For the initial conditions $\x = (1, 1), (2, 2), (1, 2), (2, 1)$, the zeros also disappear when the library is reduced to $\{x, y, xy\}$ for each component equation. Thus, both the training initial conditions and the candidate library affect the separation of true and spurious coefficients.

\subsection{Initial parameter}\label{sec:a0}

\begin{figure}[!tbp]
    \centering
    \begin{subfigure}[t]{0.48\textwidth}
        \centering
        \includegraphics[width=\linewidth]{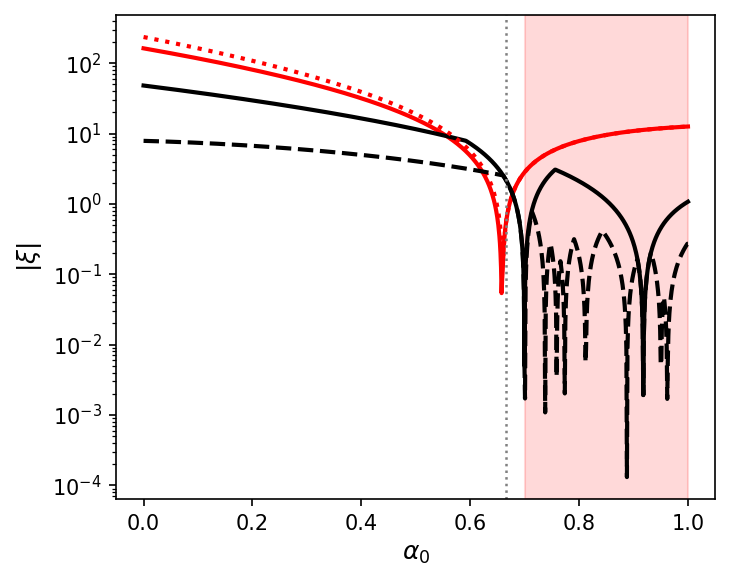}
    \end{subfigure}
    \hfill
    \begin{subfigure}[t]{0.48\textwidth}
        \centering
        \includegraphics[width=\linewidth]{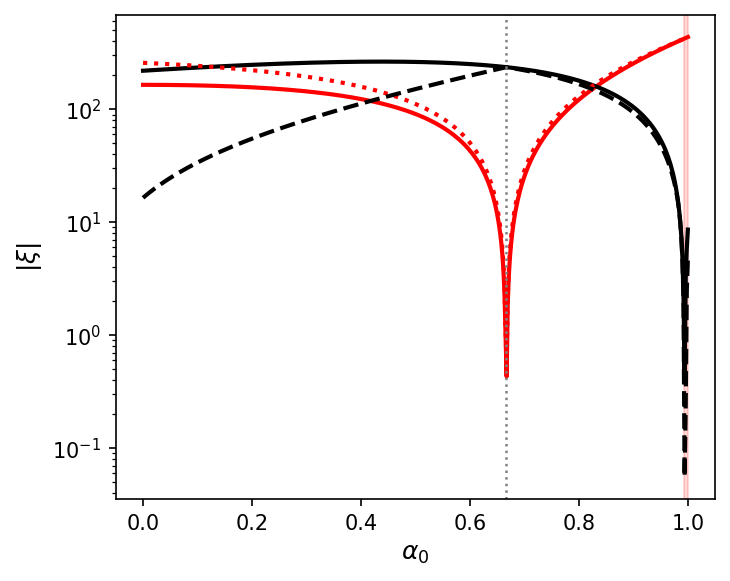}
    \end{subfigure}

    \vspace{0.5em}

    \begin{subfigure}[t]{0.48\textwidth}
        \centering
        \includegraphics[width=\linewidth]{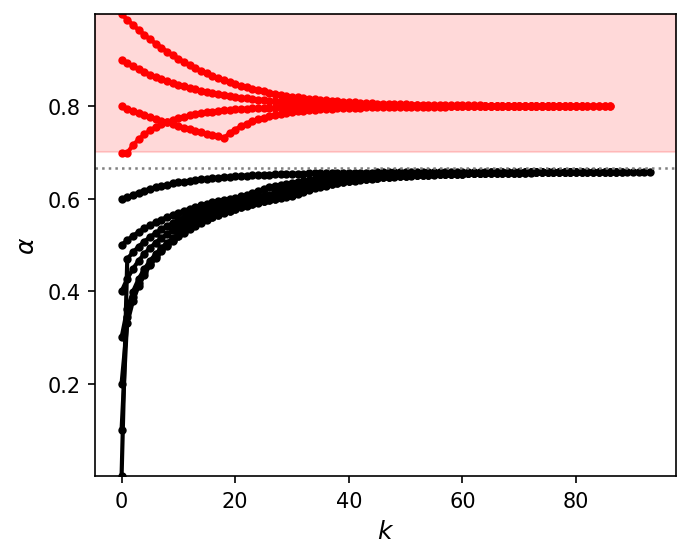}
    \end{subfigure}
    \hfill
    \begin{subfigure}[t]{0.48\textwidth}
        \centering
        \includegraphics[width=\linewidth]{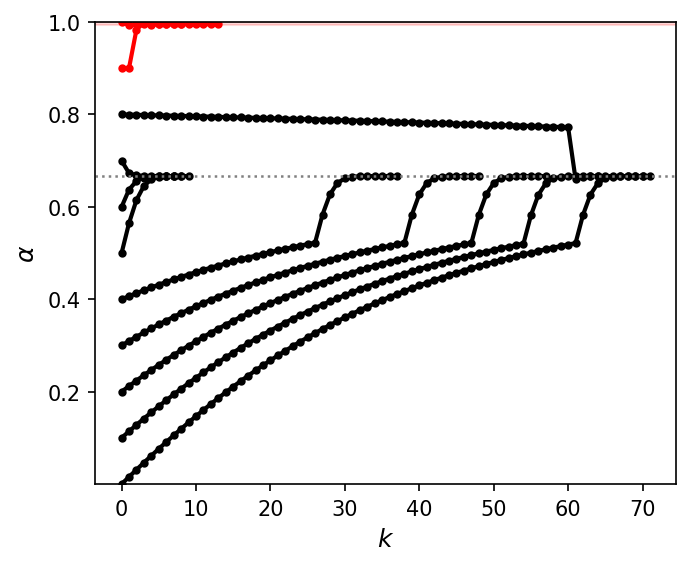}
    \end{subfigure}
    \caption{Dependence of pSTLSQ on the initial parameter for the fractional Riccati equation (left) and fractional heat equation with $u_{xx}$ included in the candidate library (right). In the top panels, the black dashed and solid curves show $U$ and $\widetilde U$, respectively, while the red dotted and solid curves show $L$ and $\widetilde L$. The vertical dotted line indicates $\a_*$, and the shaded region contains initial parameters for which the interval between $\a_0$ and $\a_*$ crosses a zero of $U$. The bottom panels show the corresponding parameter trajectories. Black trajectories recover the true support, whereas red trajectories do not. The horizontal dotted line indicates $\a_*$.}
    \label{fig:alpha0_dependence}
\end{figure}

We next examine how support recovery and parameter convergence depend on the initial parameter $\a_0$, while keeping the other settings fixed. For the fractional Riccati equation, $U$ vanishes at $\a \approx 0.701$ and at further points up to $\a \approx 0.962$. For reference intervals containing these zeros, condition \eqref{eq:lambda_condition_reference} gives no positive threshold $\lambda$. Figure \ref{fig:alpha0_dependence} shows runs from $\a_0 = 10^{-6}, 0.1, 0.2, \dots, 0.9, 1 - 10^{-6}$ with the settings of Section \ref{sec:ode}. All tested initializations with $\a_0 \leq 0.6$ recover the true support, whereas those with $\a_0 \geq 0.7$ terminate near $\hat\a \approx 0.80$ with incorrect supports. At $\a_0 = 0.7$, we have $\phi^i(\a_0) = -0.023$ for the index $i$ corresponding to $u^2$, so the first thresholding step removes this true term for both $\lambda = 2$ and $\lambda = 0.05$. Note that this shortcoming is a result of the dataset at hand, as replacing the training initial conditions by $\x = -0.8, -0.3, 0.5, 0.9$ removes the zeros of $U$, thus allowing a feasible threshold $\lambda$ for any $\a_0 \in \Omega_\a$.

For the fractional heat equation, the candidate library contains $u_{xx}$, and since $(-\Delta)u = -u_{xx}$ for every smooth periodic $u$, the model containing only $u_{xx}$ reproduces $b_\a$ at $\a = 1$ for any data, up to discretization errors. This competing representation makes $u_t$ redundant near $\a = 1$, so that $U$ vanishes there. Starting from $\a_0 = 0.9$ or $\a_0 = 1$, pSTLSQ therefore returns this model at $\hat\a \approx 0.995$, whereas all tested initializations with $\a_0 \leq 0.8$ recover the true support, as shown in Figure \ref{fig:alpha0_dependence}. Omitting $u_{xx}$ removes this competing representation, as in the Gaussian kernel example of Section \ref{sec:kernel}. Without $u_{xx}$, however, the magnitude of the fitted coefficient of $u_t$ at $\a_0 = \frac{1}{3}$ drops to about $90$, below the threshold $\lambda = 100$ chosen for the library with $u_{xx}$, so a smaller threshold or a larger initial parameter is needed in that case, see Section \ref{sec:comparison}.

\begin{figure}[!tbp]
    \centering
    \begin{subfigure}[t]{0.49\textwidth}
        \centering
        \includegraphics[width=\linewidth]{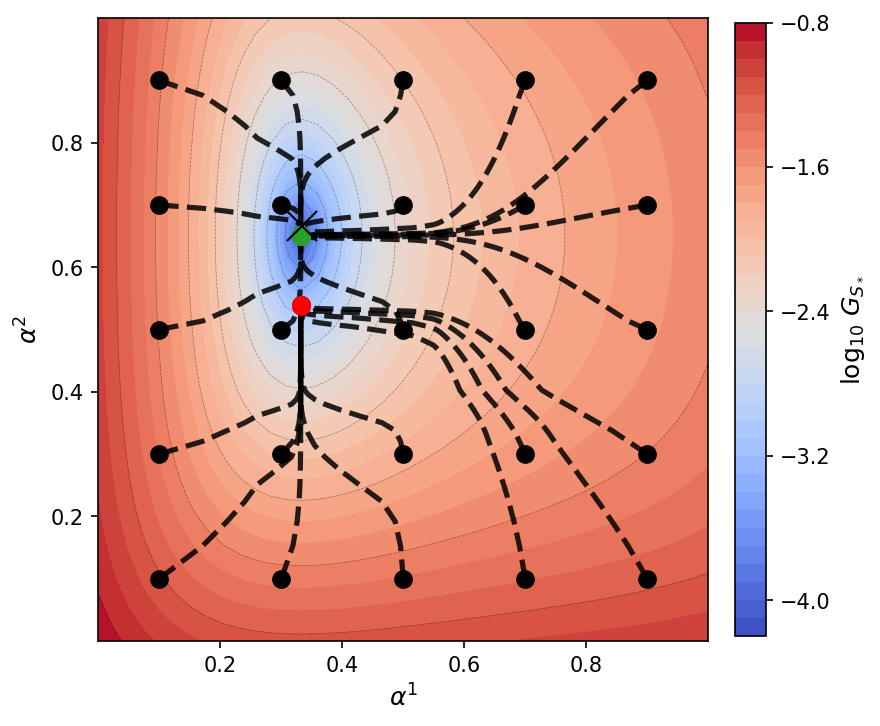}
    \end{subfigure}
    \hfill
    \begin{subfigure}[t]{0.49\textwidth}
        \centering
        \includegraphics[width=\linewidth]{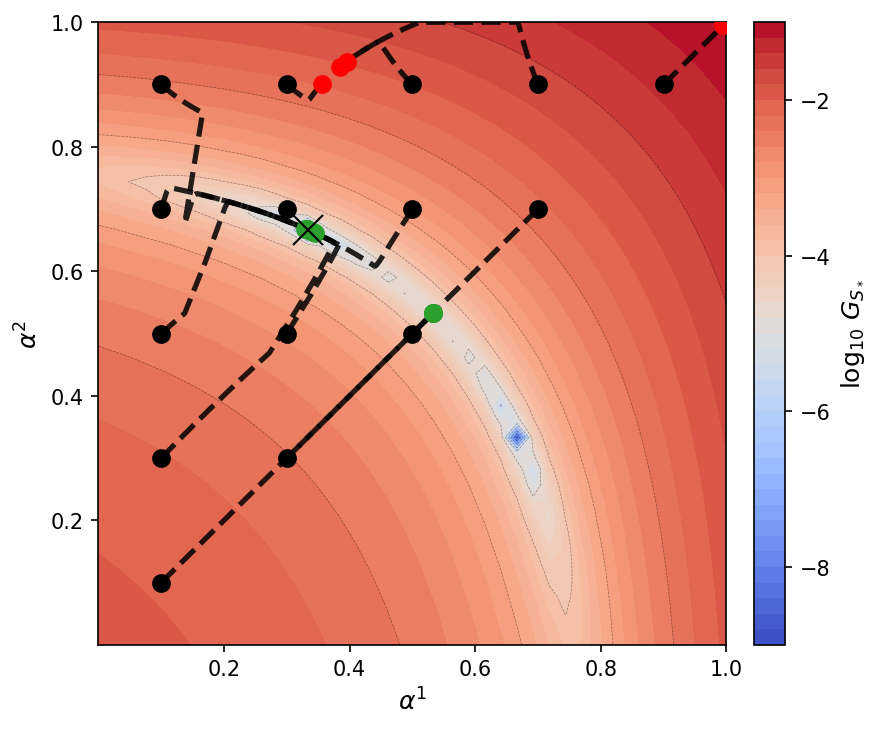}
    \end{subfigure}
    \caption{Parameter trajectories for initial values satisfying $\a_0^1, \a_0^2 \in \{0.1, 0.3, 0.5, 0.7, 0.9\}$ in the fractional Van der Pol example (left) and the multi-order space-fractional heat example (right), superimposed on $\log_{10}\overline G_{S_*}$, computed on the noiseless data. Green endpoints indicate recovery of the true support, red endpoints indicate failure to recover it, and the black cross marks the true parameter. A green endpoint need not coincide with $\a_*$. By symmetry, only trajectories with $\a^2 \geq \a^1$ are shown in the right panel.}
    \label{fig:ak_trajectories_2d}
\end{figure}

Figure \ref{fig:ak_trajectories_2d} shows the parameter trajectories from different initializations for the fractional Van der Pol oscillator and the multi-order heat equation. For the fractional Van der Pol oscillator, the reference intervals $I_{\a_*^2}(\a_0^2)$ for $\a_0^2 \in \{0.1, 0.3\}$ contain zeros of $U$ for the second equation. All these tested initializations return $\hat\a \approx (0.332, 0.539)$, with $x^2y$ removed and the spurious terms $1$ and $x^2$ retained in the second equation. In contrast, all tested initializations with $\a_0^2 \geq 0.5$ recover the true support. For the multi-order heat equation, $u_{xx}$ again provides a competing approximate fit near $\a = (1, 1)$, where $\Bb_\a(u) = -2u_{xx}$ at the continuous level. All tested initializations $\a_0 = (\a_0^1, 0.9)$ with $\a_0^1 \in \{0.3, 0.5, 0.7, 0.9\}$ retain $u_{xx}$ and fail to recover the true support. In addition, equal initial orders cannot separate under projected gradient descent, as discussed in Section \ref{sec:symmetry}.

These results show that initialization affects both support and parameter recovery, and that recovering the true support does not guarantee recovery of the true parameter. It is therefore useful to test several initial parameters and thresholds and compare the learned models through their predictions on independent data.

\subsection{Noise sensitivity and smoothing}\label{sec:failures}

Measurement noise perturbs both the library matrix $A$ and $b_\a$, and numerical differentiation can amplify these perturbations. These perturbations can affect the estimated parameters even when the true support is recovered. At higher noise levels, thresholding can also remove true terms. For the fractional and multi-order heat equations and the two kernel examples, the true support consists only of $u_t$. 
This term is removed when the magnitude of its fitted coefficient falls below $\lambda$, and if the support is then empty the parameter updates act on
\begin{align*}
    G_\emptyset(\a) = \|b_\a\|_2^2,
\end{align*}
whose minimizer might not be close to $\a_*$. In these examples, all tested runs that remove $u_t$ and retain no other terms return a parameter on the boundary of $\Omega_\a$. For the fractional heat equation, this boundary behavior follows from the fact that $G_\emptyset$ is minimized at the smallest allowed order. To see this, note that the operator multiplies each nonzero spatial Fourier mode of the observed data by $|k|^{2\a}$. Parseval's identity therefore expresses $G_\emptyset(\a)$ as a sum of squared Fourier coefficient magnitudes weighted by $|k|^{4\a}$. The Fourier coefficients of the observed data are fixed as $\a$ varies. Since $|k| \geq 1$ for every nonzero integer mode, each weight is nondecreasing in $\a$. Hence, when no library terms remain, the least-squares error is minimized at the smallest allowed order, $\a = 10^{-6}$. Noisy fits can also retain spurious terms, including $u$ for the power-law kernel and higher-degree polynomials for the Riccati equation.

To reduce these effects, we smooth the observations in time before constructing $A$ and $b_\a$. We use a moving average and a Gaussian filter with fixed widths, and a cubic smoothing spline whose smoothing strength is chosen for each time series by generalized cross-validation \cite{wahba1990spline}. For the fractional heat equation, we also apply Gaussian smoothing in both time and space. All methods use the same $20$ noisy datasets at each positive noise level, and the pSTLSQ settings are unchanged.

\begin{figure}[!tbp]
    \centering
    \includegraphics[width=\textwidth]{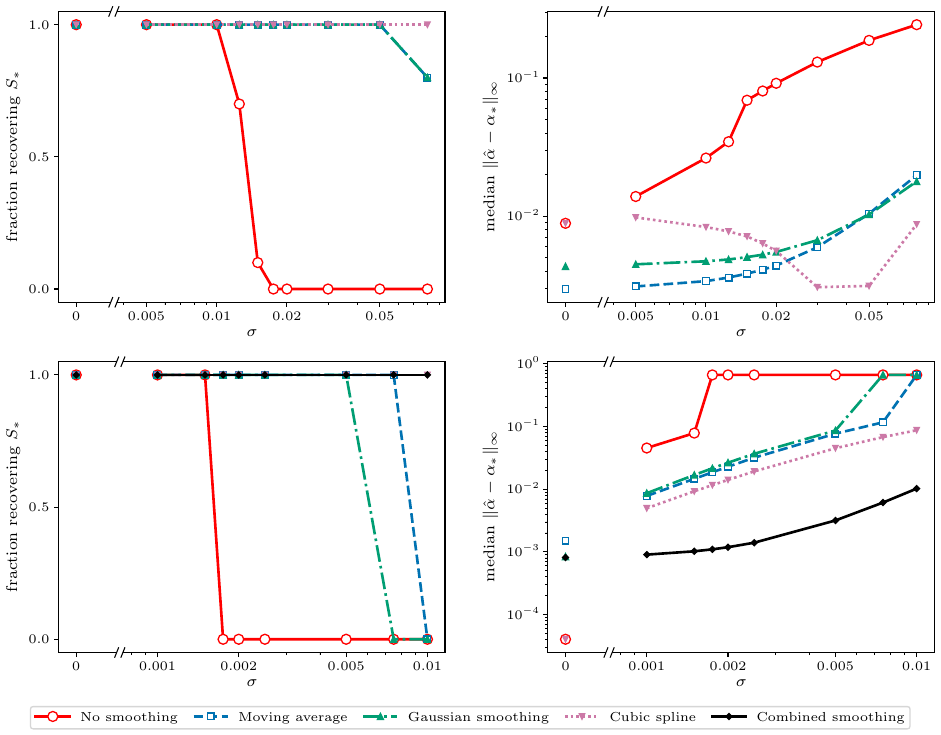}
    \caption{Effect of smoothing on pSTLSQ for the fractional Riccati equation (top row) and the fractional heat equation (bottom row), with $20$ noisy datasets per noise level and a single noiseless dataset at $\s = 0$. The left panels show the fraction of datasets for which pSTLSQ recovers the true support, and the right panels show the median parameter error $\|\hat\a - \a_*\|_\infty$ over all runs. Temporal smoothing uses a moving average over $9$ samples, a Gaussian filter with standard deviation $2$ samples, or a cubic smoothing spline with strength chosen by generalized cross-validation. Combined smoothing applies the Gaussian filter followed by the spatial filter $\widehat u_k \mapsto e^{-(k/16)^2}\widehat u_k$.}
    \label{fig:smoothing}
\end{figure}

Figure \ref{fig:smoothing} shows the resulting improvement in support recovery. We define the collapse level as the smallest tested noise level at which fewer than $10$ of the $20$ runs recover the true support. For the Riccati equation, all three smoothing methods raise this level from $\s = 0.015$ to beyond the largest tested level $\s = 0.08$. At $\s = 0.08$, the moving-average and Gaussian filters each recover the true support in $16$ of $20$ runs, while the spline recovers it in all $20$. For the fractional heat equation, the Gaussian and moving-average filters raise the collapse level from $1.75 \cdot 10^{-3}$ to $7.5 \cdot 10^{-3}$ and $10^{-2}$, respectively. The spline and combined temporal and spatial smoothing recover the true support in all $20$ runs at every tested level up to $\s = 10^{-2}$.

The parameter estimates also improve in the displayed noisy examples. For the Riccati equation at $\s = 0.02$, the median parameter error decreases from $9.2 \cdot 10^{-2}$ without smoothing to between $4.4 \cdot 10^{-3}$ and $5.6 \cdot 10^{-3}$ with temporal smoothing. For the fractional heat equation at $\s = 1.5 \cdot 10^{-3}$, adding spatial smoothing to Gaussian temporal smoothing reduces the median error from $1.7 \cdot 10^{-2}$ to $10^{-3}$. On noiseless data, however, the fixed filters reduce the Riccati parameter error by partly cancelling the discretization error, but increase the error for the fractional heat equation. In these tests, the spline selects interpolation for noiseless data and returns the same estimate as the raw observations. Smoothing therefore improves noise tolerance in these tests, but can also alter parameter estimates when noise is absent.

\subsection{Parameter and support updates}\label{sec:eta}

\begin{figure}[!tbp]
    \centering
    \includegraphics[width=0.8\textwidth]{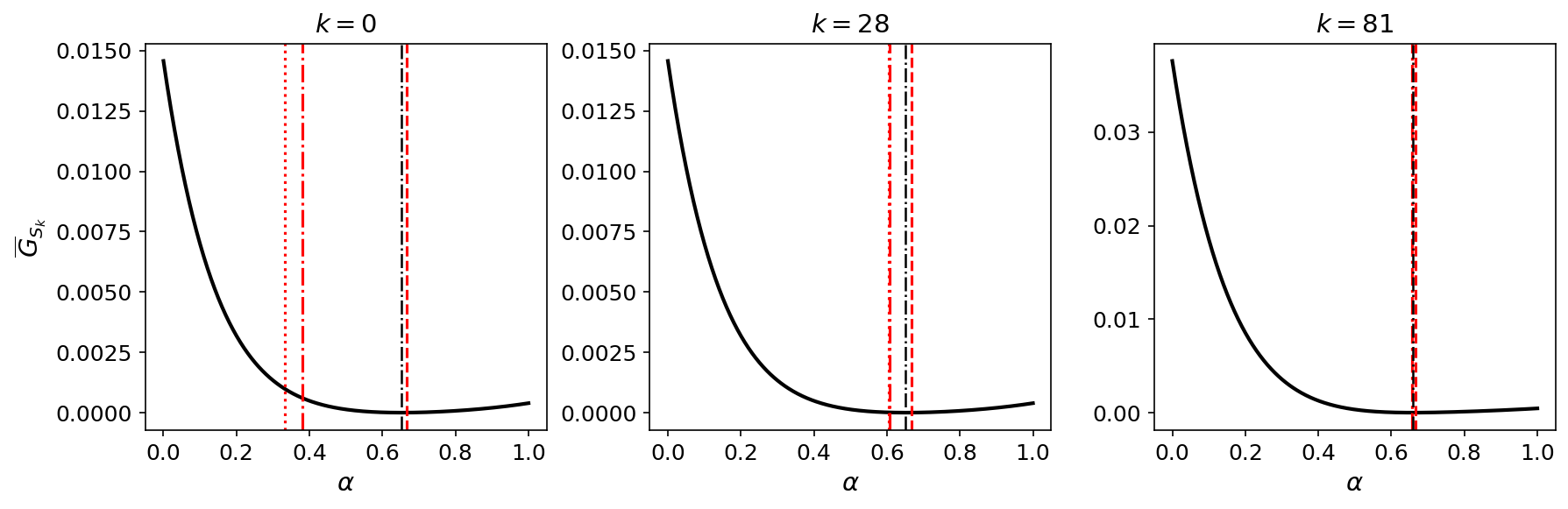}
    \caption{Reduced objectives and parameter updates for the fractional Riccati equation at iterations $k = 0, 28, 81$. The dotted and dash-dotted red lines indicate $\a_k$ and $\a_{k + 1}$, respectively. The dashed red line marks the true parameter $\a_*$, and the dashed black line marks a minimizer of $G_{S_k}$. Here, $n_\a = 1$ and $\eta = 5$.}
    \label{fig:alpha_convergence_1d}
\end{figure}

Recall that the inexact pSTLSQ method used here takes one projected gradient step on $\overline G_{S_k}$ per iteration and refits the coefficients, as illustrated in Figure \ref{fig:alpha_convergence_1d} for the fractional Riccati equation on noiseless data with the settings of Section \ref{sec:ode}. After the initial thresholding, the support changes only when $S_{29}$ and $S_{30}$ are formed, and the remaining iterations then refine the parameter on $S_{30} = S_*$ until the stopping criterion is met. By Lemma \ref{lem:alpha_decreasing_fixedS_GD}(a), a learning rate $0 < \eta < 2/\overline L_S$ ensures that $\overline G_S$ does not increase on a fixed support $S$, where $\overline L_S$ is a Lipschitz constant of $\nabla\overline G_S$. Large learning rates can prevent convergence, as for the Lotka-Volterra system at $\s = 0.05$, where increasing $\eta$ from $\frac{1}{4}$ to $\frac{1}{2}$ causes the iterates to alternate between the bounds of $\Omega_\a$. To test the effect of the number $n_\a$ of gradient steps per iteration, we take $40$ gradient steps in total with $n_\a = 1, 2, 4, 10, 20, 40$. For the fractional Riccati equation on noiseless data, the parameter error decreases from $0.042$ for $n_\a = 20$ and $40$ to $0.027$ for $n_\a = 1$, whereas $n_\a = 20$ gives the smallest error for the multi-order heat equation. In the latter case the weak curvature discussed in Section \ref{sec:disc_floor} slows gradient descent, so that no run comes close to $\a_*$ within $40$ gradient steps.

\subsection{Symmetric parameterizations}\label{sec:symmetry}

The two orders in the multi-order heat equation enter symmetrically, so $G_S(\a^1, \a^2) = G_S(\a^2, \a^1)$. When the orders are equal, the partial derivatives of $G_S$ with respect to the two orders are also equal. A gradient step with the same learning rate for both orders therefore keeps them equal. Since $\Omega_\a$ is invariant under exchanging the orders and its Euclidean projection is unique, projection also keeps the orders equal, and equal initial orders consequently remain equal throughout the iteration. All tested diagonal initializations with $\a_0^1 = \a_0^2 \leq 0.7$ recover the true support but return a parameter near $(0.53, 0.53)$, as shown in Figure \ref{fig:ak_trajectories_2d}. Recovering the two distinct true orders therefore requires initializing them at different values.

\section{Comparison with other methods}\label{sec:comparison}

In this section, we benchmark inexact pSTLSQ, which we simply refer to as pSTLSQ as in Sections \ref{sec:numerics} and \ref{sec:sensitivity}, against two other methods for the pSINDy optimization problem \eqref{eq:psindy0} on the fractional Riccati equation, the fractional Van der Pol oscillator, and the fractional heat equation. All methods use the data, library, discretization of $b_\a$, and threshold $\lambda$ of Section \ref{sec:numerics}, with $20$ datasets tested whenever noise is imposed. We consider both raw data and data smoothed as in Section \ref{sec:failures}, using the smoothing spline for the two ODE examples and the combined temporal and spatial smoothing for the fractional heat equation.

The first benchmark method is a grid search, which applies STLSQ with threshold $\lambda$ at each point of a uniform grid over $\Omega_\a$, with $101$ points for $m = 1$ and $41$ points per parameter for $m = 2$. This method then returns the candidate model with the smallest held-out error, i.e.,
\begin{align*}
    \biggl(\frac{\sum_{\ell = 1}^n\|v_\x^\ell - u_{\x, \s}^\ell\|_2^2}{\sum_{\ell = 1}^n\|u_{\x, \s}^\ell\|_2^2}\biggr)^{1/2},
\end{align*}
where $u_{\x, \s}^\ell$ and $v_\x^\ell$ denote the $\ell$th components of the noisy data \eqref{eq:u_noisy} and of the solution of the candidate model for the test initial condition $\x$ of Section \ref{sec:numerics}, respectively. We select by this error rather than by the pSINDy objective $F$, since the candidate with the smallest $F$ might not give the true support. For example, for the fractional Riccati equation of Section \ref{sec:ode} with $\lambda = 2$, the model containing only $u$ has a smaller objective $F$ than the true model containing $u$ and $u^2$, see Section \ref{sec:thresh_num}, since removing $u^2$ increases the least-squares error by less than the penalty $\lambda^2$ it saves. Selecting by $F$ therefore returns this one-term model instead of the true model on every dataset of this example. However, this selection requires one forward solve per grid point, i.e., $101$ solves for $m = 1$ and $1681$ for $m = 2$, and this number grows exponentially with the number of parameters. This grid search is similar to the approach of \cite{zhang2024sparse}. The second benchmark method is a decoupled method, which separates the parameter estimation from the thresholding by first minimizing the reduced objective over $\a$ with the full library, then applying STLSQ at the resulting parameter, and finally refitting $\a$ on the selected support by minimizing the reduced objective again. Here, the minimizations over $\a$ start from several points of $\Omega_\a$ to reduce the dependence on an initial parameter. We note that both benchmark methods still fit the pSTLSQ iteration scheme \eqref{eq:iteration_scheme}, but with a different schedule than inexact pSTLSQ. In particular, the grid search runs this scheme with $\nu_k = 0$, $k \geq 0$, from each grid point, while the decoupled method uses the schedule $\nu_{-1} = \infty$, followed by $\nu_k = 0$ until the support no longer changes and one update with $\nu_k = \infty$.

Table \ref{tab:comparison} reports the parameter errors at the noise level used for each example in Section \ref{sec:numerics}, where all three methods recover the true support on every dataset. pSTLSQ and the decoupled method have nearly the same errors and underestimate the parameter. Since minimizing the reduced objective $G_{S_*}$ directly on the true support gives the same errors, this bias is not an optimization error as in Section \ref{sec:disc_floor}, but a shift of the minimizer of $G_{S_*}$ due to discretization and noise. A possible explanation for the shift to smaller values is that the noise in $b_\a$ contributes approximately $\|P_{S_*}D_\a\varepsilon\|_2^2$ to $G_{S_*}(\a)$, with $D_\a$ and $\varepsilon$ as in Appendix \ref{sec:stab}, and that the expected value of this term increases with $\a$ since the L1 scheme scales with $\Delta t^{-\a}$ -- where $\Delta t < 1$ -- and the fractional Laplacian multiplies each nonzero Fourier mode by $|k|^{2\a}$. In contrast, the grid search selects the parameter by the held-out error instead of the reduced objective, which gives a smaller error for the fractional Riccati equation and the fractional heat equation, but larger errors for the Van der Pol oscillator. Smoothing the data as in Section \ref{sec:failures} also reduces the bias of pSTLSQ, to $-0.008$ for the fractional Riccati equation and $-0.0009$ for the fractional heat equation.

\begin{table}[!tbp]
    \centering
    \small
    \begin{tabular}{@{}lccc@{}}
    \toprule
    Method & fractional Riccati & fractional Van der Pol & fractional heat \\
    \midrule
    pSTLSQ & $-0.026$ & $(-0.002, -0.025)$ & $-0.045$ \\
    grid search & $-0.007$ & $(-0.008, 0.033)$ & $-0.007$ \\
    decoupled method & $-0.026$ & $(-0.002, -0.025)$ & $-0.045$ \\
    \bottomrule
    \end{tabular}
    \caption{Parameter error $\hat\a - \a_*$ (median) of pSTLSQ, the grid search, and the decoupled method over the $20$ noisy datasets of each example at the noise levels of Section \ref{sec:numerics}. All three methods recover the true support on every dataset.}
    \label{tab:comparison}
\end{table}

Tables \ref{tab:comparison_raw} and \ref{tab:comparison_smooth} compare the support recovery of pSTLSQ and the decoupled method at increasing noise levels, on raw and smoothed data, respectively, and show that neither method is uniformly better. We do not include the grid search here, as these tables test how the order of the parameter updates and the thresholding affects support recovery under noise, whereas the grid search keeps the parameter fixed at each grid point. On raw data, the decoupled method recovers the true support on more datasets for the fractional Riccati equation and the fractional heat equation, and pSTLSQ for the Van der Pol oscillator, while on smoothed data, pSTLSQ recovers it on more datasets for the fractional Riccati equation, and the decoupled method for the Van der Pol oscillator. The largest difference occurs for the smoothed fractional heat data, where pSTLSQ recovers the true support on all datasets at every noise level, whereas the decoupled method returns the model containing only $u_{xx}$ at $\a \approx 0.995$ on every noisy dataset. This is due to the competing representation of Section \ref{sec:a0}, as the decoupled method first fits $\a$ with the full library and thus moves to the parameter value at which $u_{xx}$ reproduces $b_\a$. In contrast, pSTLSQ with $\a_0 = \frac{1}{3}$ removes $u_{xx}$ by thresholding and returns a parameter near $\a_*$. Without $u_{xx}$ in the library, the decoupled method recovers the true support on all smoothed datasets, and so does pSTLSQ with $\lambda = 50$ instead of $\lambda = 100$, or with $\a_0 = \frac{1}{2}$ instead of $\a_0 = \frac{1}{3}$. With the settings of Section \ref{sec:numerics}, pSTLSQ fails in this case because the threshold $\lambda = 100$ was chosen for the library with $u_{xx}$ and exceeds the magnitude of the fitted coefficient of $u_t$ at $\a_0 = \frac{1}{3}$ without $u_{xx}$, as explained in Section \ref{sec:a0}.

\begin{table}[!tbp]
    \centering
    \small
    \begin{tabular}{@{}llccccc@{}}
    \toprule
    & $\s$ & $0.01$ & $0.015$ & $0.0175$ & $0.05$ & $0.08$ \\
    \midrule
    fractional Riccati & pSTLSQ & 20 & 2 & 0 & 0 & 0 \\
     & decoupled method & 20 & 15 & 5 & 0 & 0 \\
    fractional Van der Pol & pSTLSQ & 20 & 20 & 20 & 17 & 6 \\
     & decoupled method & 20 & 20 & 20 & 5 & 0 \\
    \midrule
    & $\s$ & $0.001$ & $0.0015$ & $0.00175$ & $0.005$ & $0.01$ \\
    \midrule
    fractional heat & pSTLSQ & 20 & 20 & 0 & 0 & 0 \\
     & decoupled method & 20 & 20 & 4 & 0 & 0 \\
    \bottomrule
    \end{tabular}
    \caption{Number of datasets (out of $20$) on which pSTLSQ and the decoupled method recover the true support from raw data at noise level $\s$. The fractional heat equation is tested at smaller noise levels than the two ODE examples.}
    \label{tab:comparison_raw}
\end{table}

\begin{table}[!tbp]
    \centering
    \small
    \begin{tabular}{@{}llccccc@{}}
    \toprule
    & $\s$ & $0.01$ & $0.015$ & $0.0175$ & $0.05$ & $0.08$ \\
    \midrule
    fractional Riccati & pSTLSQ & 20 & 20 & 20 & 20 & 20 \\
     & decoupled method & 20 & 20 & 20 & 12 & 2 \\
    fractional Van der Pol & pSTLSQ & 20 & 20 & 20 & 14 & 11 \\
     & decoupled method & 20 & 20 & 20 & 20 & 20 \\
    \midrule
    & $\s$ & $0.001$ & $0.0015$ & $0.00175$ & $0.005$ & $0.01$ \\
    \midrule
    fractional heat & pSTLSQ & 20 & 20 & 20 & 20 & 20 \\
     & decoupled method & 0 & 0 & 0 & 0 & 0 \\
    fractional heat without $u_{xx}$ & pSTLSQ & 0 & 0 & 0 & 0 & 0 \\
     & pSTLSQ, $\lambda = 50$ & 20 & 20 & 20 & 20 & 20 \\
     & decoupled method & 20 & 20 & 20 & 20 & 20 \\
    \bottomrule
    \end{tabular}
    \caption{As Table \ref{tab:comparison_raw}, but for smoothed data. The data of the fractional Riccati equation and the fractional Van der Pol oscillator are smoothed in time by the cubic smoothing spline of Section \ref{sec:failures}, and the data of the fractional heat equation by the combined temporal and spatial smoothing of Section \ref{sec:failures}. All runs use the settings of Section \ref{sec:numerics}, also when $u_{xx}$ is removed from the library, except in the row with $\lambda = 50$.}
    \label{tab:comparison_smooth}
\end{table}

In summary, all three methods recover the true support on the raw data at the noise levels of Section \ref{sec:numerics}, but differ in parameter accuracy, cost, and sensitivity to initialization and the candidate library. The grid search does not depend on an initial parameter and, since it selects by the held-out error instead of the reduced objective, has a smaller error in two of the three examples, but requires one forward solve per grid point. pSTLSQ and the decoupled method both underestimate the parameter under noise, since both minimize the reduced objective, and have similar accuracy when both recover the true support, but fail in different ways. The decoupled method fits the parameter before any term is removed, which makes it sensitive to competing representations as in Section \ref{sec:a0}, whereas pSTLSQ depends on the initial parameter and requires the true coefficients to exceed $\lambda$ at $\a_0$, as discussed in Sections \ref{sec:thresh_num} and \ref{sec:a0}. As both are schedules of iteration scheme \eqref{eq:iteration_scheme}, this shows that the schedule affects which local minimizer is reached.

\section{Discussion}\label{sec:conclusion}

We developed pSTLSQ for the pSINDy optimization problem and tested it on equations involving Caputo derivatives, fractional Laplacians, and parameterized nonlocal kernels. For each support $S$ and parameter $\a$, the library coefficients are obtained by least squares and are unique under the full column rank assumption on $A$. Eliminating these coefficients gives a reduced objective for $\a$, which can be optimized using derivative information from $b_\a$. Although the numerical implementation uses only first-order information, the parameter update can also use higher-order information when available. The iteration is formulated for a general schedule of parameter updates, which includes STLSQ on a parameter grid and a decoupled method as special cases.

Under suitable assumptions, exact pSTLSQ terminates at a local minimizer and inexact pSTLSQ converges locally to such a minimizer. We also derived sufficient conditions for recovery of the true support. These conditions assume parameter convergence and require a threshold that preserves true terms throughout the iteration while eventually removing spurious ones. Whether such a threshold exists depends on the dataset, candidate library, and parameter trajectory. Although the local minimizers of the pSINDy objective are independent of $\lambda>0$, the thresholding parameter affects which are globally optimal, which are fixed points of pSTLSQ, and which supports the iteration can reach. On a fixed grid, we also proved that strict, nondegenerate interior local minimizers retain their support under sufficiently small data perturbations, with their parameters and coefficients varying continuously with the data.

In the numerical experiments, we used inexact pSTLSQ with one projected gradient step per iteration, referred to as pSTLSQ below, which recovered the governing terms and approximated the true parameters at low noise levels for suitable initializations and thresholds. A linear perturbation estimate indicated that, in most examples, the parameter error without added noise was mainly caused by discretization errors. Weak curvature can also slow gradient descent, making the choice of learning rate and stopping tolerance important. Threshold selection depended on the training data, candidate library, and initialization, while competing representations and parameter symmetries also affected recovery. Support changes could occur after many parameter updates, and large learning rates could prevent convergence. Cubic smoothing splines and standard filters improved recovery under noise but could also shift noiseless estimates. In the benchmark of Section \ref{sec:comparison}, pSTLSQ, the grid search, and the decoupled method all recovered the true support on raw data at low noise levels. We note that both benchmark methods are also pSTLSQ with a different schedule. The grid search had smaller parameter errors in two of the three examples, as it selects among its candidate models by a held-out error instead of the reduced objective, but required one forward solve per grid point. pSTLSQ and the decoupled method had similar accuracy when both recovered the true support, and neither was uniformly better at higher noise levels. Their differences show that the schedule determines which local minimizer is reached, with the decoupled method being more sensitive to competing representations in the library and pSTLSQ to the threshold and the initial parameter.

These results suggest several directions for future work. One is to combine pSTLSQ with the fractional weak formulations developed in \cite{thanasutives2026robust}, which avoid the numerical differentiation of noisy data, and to examine to what extent this improves parameter accuracy and support recovery under noise compared to the smoothing of Section \ref{sec:failures}. Another extension is to include unknown parameters $\beta$ in the candidate library \cite{champion2020unified}, i.e., to replace $A$ by a parameterized library matrix $A_\beta$. For fixed $(\a, \beta)$ and support $S$, the least-squares coefficients remain unique if $A_{\beta, S}$ has full column rank, but the parameter update then involves derivatives of both $A_\beta$ and $b_\a$, so that extending the convergence and support recovery analysis requires suitable rank and smoothness assumptions on $A_\beta$. Moreover, the slow convergence under weak curvature observed in Section \ref{sec:disc_floor} motivates comparing projected gradient descent with parameter updates that use curvature information, such as Newton or quasi-Newton methods, or Gauss-Newton methods as commonly used with variable projection \cite{golub1973differentiation}, with respect to both the cost of the updates and their effect on parameter convergence and support selection. An SR3 formulation \cite{zheng2018unified} could also be used for the coefficient estimation, although its relaxed coefficient fit changes the reduced objective and thus requires revisiting the convergence analysis. A further direction is a convergence analysis beyond the local setting of Theorem \ref{thm:local_inexact}.

Finally, the benchmark of Section \ref{sec:comparison} covers only three examples and a few schedules. For a better understanding of pSTLSQ, more schedules need to be tested on more datasets, e.g., schedules between pSTLSQ and the decoupled method with a few parameter updates on the full library before the first thresholding step, together with a comparison with the methods of \cite{zhang2024sparse,fan2023data,xiangnan2025data,thanasutives2026robust} on the same datasets. Moreover, while the grid search is attractive in this work due to its simplicity and comparable performance, its complexity grows exponentially with the number of parameters, which makes it infeasible for high-dimensional parameters, e.g., for a spatially varying fractional order with one parameter per grid point. In contrast, inexact pSTLSQ does not require a parameter grid, and a gradient step requires only the Jacobian $\nabla b_\a \in \RR^{N \times m}$, i.e., one derivative of $b_\a$ per parameter. We therefore propose to test inexact pSTLSQ in this setting, for example by running it from several randomly sampled initial parameters and, as for the grid search, selecting among the resulting models by the held-out error.

\section*{Acknowledgement}

This work was supported in part by the Air Force Office of Scientific Research (FA9550-24-1-0141).

\printbibliography

\appendix

\section{Stability of local minimizers} \label{sec:stab}

In this appendix, we study the stability of local minimizers of the pSINDy objective under perturbations of the data. We prove that strict nondegenerate interior local minimizers are stable under sufficiently small perturbations of the library matrix and $\a$-dependent perturbations of $b_\a$ in the context of measurement noise on a fixed grid. We note that this result concerns the optimization problem itself and does not establish convergence of pSTLSQ to the perturbed minimizer.

\begin{theorem} \label{thm:noise}
Let $\lambda > 0$, let $(\a_*, \xi_*)$ be a strict local minimizer of optimization problem \eqref{eq:sindy_alpha} with $\a_* \in \operatorname{int}(\Omega_\a)$, and let $S_* := \supp(\xi_*)$. Suppose that Assumption \ref{ass:A} holds. Moreover, let $r > 0$ be such that $\overline{B_r(\a_*)} \subseteq \operatorname{int}(\Omega_\a)$ and $b_\a$ is $\Cc^2$ on a neighborhood of $\overline{B_r(\a_*)}$. Given $E \in \RR^{N \times p}$ and $e \in \Cc^2(\overline{B_r(\a_*)};\RR^N)$, define
\begin{align*}
    F(\a, \xi;E, e) & := \|(A + E)\xi - (b_\a + e(\a))\|_2^2 + \lambda^2\|\xi\|_0, \\
    H_{S}(\a, \z;E, e) & := \|(A_{S} + E_{S})\z - (b_\a + e(\a))\|_2^2,
\end{align*}
for $\a \in \overline{B_r(\a_*)}$, $\xi \in \RR^p$, $\z \in \RR^{|S|}$, $S \subseteq \{1, 2, \dots, p\}$, and let
\begin{align*}
    \|e\|_{\Cc^2} := \sup_{\a \in \overline{B_r(\a_*)}}\max\{\|e(\a)\|, \|\nabla e(\a)\|, \|\nabla^2e(\a)\|\}.
\end{align*}
If $\nabla^2_{(\a, \z)}H_{S_*}(\a_*, \xi_*|_{S_*};0, 0) \succ 0$, then there exist constants $\delta, C > 0$ such that the perturbed objective $F( \cdot , \cdot ;E, e)$ has a strict local minimizer $(\tilde\a_*, \tilde\xi_*)$ in a neighborhood of $(\a_*, \xi_*)$ satisfying $\supp(\tilde\xi_*) = S_*$ and
\begin{align}\label{eq:noise_bound}
    \|\tilde\a_* - \a_*\| + \|\tilde\xi_* - \xi_*\| \leq C (\|E\| + \|e\|_{\Cc^2}),
\end{align}
provided $\|E\| + \|e\|_{\Cc^2} < \delta$.
\end{theorem}

\begin{proof}
Throughout this proof, let $\Ee_{S_*} : \RR^{|S_*|} \to \RR^p$ denote the extension operator that pads a vector with zeros outside the support $S_*$, and let $\Cc^2_r := \Cc^2(\overline{B_r(\a_*)};\RR^N)$ with the norm $\| \cdot \|_{\Cc^2}$. We first prove stability of the restricted optimization problem on the fixed support $S_*$, i.e., for the optimization problem with $H_{S_*}( \cdot , \cdot ;E, e)$ as objective. Define
\begin{align*}
    \Gg(\a, \z, E, e) := \nabla_{(\a, \z)}H_{S_*}(\a, \z;E, e)
\end{align*}
for $\a \in B_r(\a_*)$, $\z \in \RR^{|S_*|}$, $E \in \RR^{N \times p}$, and $e \in \Cc^2_r$. The perturbation $e$ enters $\Gg$ only through $e(\a)$ and $\nabla e(\a)$, and the maps $(\a, e) \mapsto e(\a)$ and $(\a, e) \mapsto \nabla e(\a)$ are $\Cc^1$ on $B_r(\a_*) \times \Cc^2_r$. Since $b_\a$ is $\Cc^2$ on a neighborhood of $\overline{B_r(\a_*)}$, it follows that $\Gg$ is $\Cc^1$. Since $(\a_*, \xi_*)$ is a local minimizer of $F( \cdot , \cdot ;0, 0)$ and $S_* = \supp(\xi_*)$, restricting to sufficiently small perturbations with support contained in $S_*$ keeps the sparsity term constant and shows that $(\a_*, \xi_*|_{S_*})$ is a local minimizer of $H_{S_*}( \cdot , \cdot ;0, 0)$. Moreover, since $\a_* \in \operatorname{int}(\Omega_\a)$, the first-order optimality condition for $H_{S_*}( \cdot , \cdot ;0, 0)$ then gives
\begin{align*}
    \Gg(\a_*, \xi_*|_{S_*}, 0, 0) = \nabla_{(\a, \z)}H_{S_*}(\a_*, \xi_*|_{S_*};0, 0) = 0.
\end{align*}
By assumption,
\begin{align*}
    \nabla_{(\a, \z)}^2H_{S_*}(\a_*, \xi_*|_{S_*};0, 0) \succ 0,
\end{align*}
and hence this Hessian is invertible. Applying the implicit function theorem to $\Gg$ then gives neighborhoods $U$ of $(0, 0)$ in $\RR^{N \times p} \times \Cc^2_r$ and $V$ of $(\a_*, \xi_*|_{S_*})$, together with a unique map $\psi \in \Cc^1(U)$ satisfying $\psi(0, 0) = (\a_*, \xi_*|_{S_*})$ and
\begin{align*}
    \Gg(\psi(E, e), E, e) = 0,
\end{align*}
for every $(E, e) \in U$, where $\psi(U) \subseteq V$. Since $\nabla_{(\a, \z)}^2H_{S_*}(\a_*, \xi_*|_{S_*};0, 0) \succ 0$ and $\nabla_{(\a, \z)}^2H_{S_*}$ is continuous in $(\a, \z, E, e)$, we may shrink $U$ and $V$, choosing $V$ convex and contained in $B_r(\a_*) \times \RR^{|S_*|}$, so that
\begin{align*}
    \nabla_{(\a, \z)}^2H_{S_*}(\a, \z;E, e) \succ 0
\end{align*}
for every $(\a, \z) \in V$ and every $(E, e) \in U$. Since $\psi(U) \subseteq V$, the point $(\tilde\a_*, \tilde\z_*) := \psi(E, e)$ is the unique local minimizer of $H_{S_*}( \cdot , \cdot ;E, e)$ in $V$, and this minimum is strict because the Hessian is positive definite. To estimate its distance from $(\a_*, \xi_*|_{S_*})$, we use that $\psi$ is $\Cc^1$, so its derivative $\psi'$ is bounded on a sufficiently small ball around $(0, 0)$. After shrinking $U$ to this ball, we consider the path $t \mapsto \psi(tE, te)$ for $0 \leq t \leq 1$, which joins the unperturbed minimizer $\psi(0, 0)$ to the perturbed minimizer $\psi(E, e)$. The chain rule and the fundamental theorem of calculus then give
\begin{align*}
    \psi(E, e) - \psi(0, 0) = \int_0^1 \psi'(tE, te)[E, e] \,\d t.
\end{align*}
Since $\psi'$ is bounded on $U$, taking norms gives $C > 0$ such that
\begin{align}\label{eq:ift_bound}
    \|\tilde\a_* - \a_*\| + \|\tilde\z_* - \xi_*|_{S_*}\| \leq C(\|E\| + \|e\|_{\Cc^2})
\end{align}
for every $(E, e) \in U$.

We now show that $\tilde\xi_* := \Ee_{S_*}(\tilde\z_*)$ gives the required minimizer of the full objective $F$. Since $\xi_*^i \neq 0$ for every $i \in S_*$, estimate \eqref{eq:ift_bound} allows us to shrink $U$ so that $\tilde\z_*$ has no zero entry, and hence $\supp(\tilde\xi_*) = S_*$. Since $A$ has full column rank, we may also shrink $U$ so that $A + E$ has full column rank. Consequently, for fixed $\a$, $H_{S_*}(\a, \cdot ;E, e)$ is a strictly convex quadratic with least-squares solution $(A_{S_*} + E_{S_*})^\dagger(b_\a + e(\a))$. This solution depends continuously on $\a$ and equals $\tilde\z_*$ at $\tilde\a_*$, so strict local minimality of $(\tilde\a_*, \tilde\z_*)$ for $H_{S_*}( \cdot , \cdot ;E, e)$ implies that $\tilde\a_*$ is a strict local minimizer of $\a \mapsto \min_{\z}H_{S_*}(\a, \z;E, e)$. Theorem \ref{thm:fixed_point_local_minimizer}, applied with $A + E$, $b_\a + e(\a)$, and $\overline{B_r(\a_*)}$ in place of $A$, $b_\a$, and $\Omega_\a$, therefore shows that $(\tilde\a_*, \tilde\xi_*)$ is a strict local minimizer of $F( \cdot , \cdot ;E, e)$ on this restricted parameter domain. Since $\tilde\a_* \in B_r(\a_*) \subseteq \operatorname{int}(\Omega_\a)$, this restriction does not affect local minimality at $\tilde\a_*$. Finally, since both $\xi_*$ and $\tilde\xi_*$ vanish outside $S_*$, we have $\|\tilde\xi_* - \xi_*\| = \|\tilde\z_* - \xi_*|_{S_*}\|$, and estimate \eqref{eq:ift_bound} therefore gives estimate \eqref{eq:noise_bound}. Choosing $\delta > 0$ such that $\|E\| + \|e\|_{\Cc^2} < \delta$ implies $(E, e) \in U$ completes the proof.
\end{proof}

Theorem \ref{thm:noise} can be applied to show that sufficiently small measurement noise preserves a nearby strict local minimizer with the same support on a fixed grid. Both $E$ and $e$ depend on the measurement noise $\varepsilon$ and the chosen discretization, including the mesh sizes, but we suppress these dependencies in the notation and keep the discretization fixed throughout this discussion. Suppose that $u + \varepsilon$ is observed instead of $u$. The library matrix becomes $A + E$, while $b_\a$ becomes $b_\a + e(\a)$, with both perturbations measured relative to the noiseless quantities computed using the same discretization. Since the operators $\Bb_\a$ used in our numerical examples and their discretizations $D_\a$ are linear in $u$, we have
\begin{align*}
    e(\a) = D_\a(u + \varepsilon) - D_\a u = D_\a\varepsilon.
\end{align*}
Local Lipschitz continuity of the discretized library terms gives $\|E\| \leq C_A\|\varepsilon\|$ for sufficiently small $\varepsilon$. If $D_\a$ is $\Cc^2$ in $\a$ on $\overline{B_r(\a_*)}$, then $D_\a$ and its first two parameter derivatives are bounded there. Since $\varepsilon$ is fixed as $\a$ varies, differentiating $e(\a) = D_\a\varepsilon$ only differentiates $D_\a$. Thus, $\|e\|_{\Cc^2} \leq C_D\|\varepsilon\|$, and hence
\begin{align*}
    \|E\| + \|e\|_{\Cc^2} \leq (C_A + C_D)\|\varepsilon\|
\end{align*}
for sufficiently small $\varepsilon$. Estimate \eqref{eq:noise_bound} therefore gives parameter and coefficient changes of order $\|\varepsilon\|$ from the unperturbed minimizer. The constant depends on the discretized operators and the library, and may grow under grid refinement due to noise amplification by numerical differentiation. We note that this does not imply that pSTLSQ recovers $S_*$ from the perturbed data, since at higher noise levels pSTLSQ can lose a true term, as illustrated in Section \ref{sec:failures}, even though the perturbed objective may still have a nearby local minimizer with the correct support.

\section{Properties of projected gradient descent}\label{ap:b}

The following lemma discusses the properties of the projected gradient update \eqref{eq:inexact_map} used in Section \ref{sec:geometry}.

\begin{lemma}\label{lem:alpha_decreasing_fixedS_GD}
Let $\Omega_\a \subseteq \RR^m$ be nonempty, closed, and convex. Fix a support $S \subseteq \{1, \dots, p\}$, assume that $b_\a \in \Cc^2(\Omega_\a)$, and let $\Ii_S$ be given by equation \eqref{eq:inexact_map}.
\begin{itemize}
\item[(a)] If $\nabla G_S$ is Lipschitz continuous on $\Omega_\a$ with constant $L_S > 0$ and $0 < \eta < 2/L_S$, then
\begin{align*}
    G_S(\Ii_S(\a)) \leq G_S(\a) - \left(\frac{1}{\eta} - \frac{L_S}{2}\right) \|\Ii_S(\a) - \a\|_2^2
\end{align*}
for every $\a \in \Omega_\a$. In particular, condition \eqref{eq:inner_decrease} holds, and $G_S(\Ii_S(\a)) = G_S(\a)$ if and only if $\Ii_S(\a) = \a$.

\item[(b)] For every $\eta > 0$, every local minimizer $\a_*$ of $G_S$ on $\Omega_\a$ satisfies $\a_* = \Ii_S(\a_*)$.

\item[(c)] If $\a_*$ is a local minimizer of $G_S$ on $\Omega_\a$ with $\nabla^2G_S(\a_*) \succ 0$, then there exist $\rho > 0$, $\tilde\eta > 0$, and $c > 0$ such that, for every $0 < \eta \leq \tilde\eta$,
\begin{align*}
    \|\Ii_S(\a) - \a_*\|_2 \leq (1 - \eta c)\|\a - \a_*\|_2
\end{align*}
for every $\a \in B_\rho(\a_*) \cap \Omega_\a$. In particular, condition \eqref{eq:contraction} holds for $S_* = S$ with $\rho_0 = \rho$ and $\vartheta = 1 - \eta c \in [0, 1)$.
\end{itemize}
\end{lemma}

\begin{proof}
Parts (a) and (b) are standard properties of projected gradient descent, see \cite[Chapter 10]{beck2017first}. The sufficient decrease estimate in (a) implies its equality statement because $\frac{1}{\eta} - \frac{L_S}{2} > 0$. For (b), every local minimizer on the convex set $\Omega_\a$ satisfies
\begin{align*}
    \nabla G_S(\a_*)^\top(\a - \a_*) \geq 0
\end{align*}
for every $\a \in \Omega_\a$. By the characterization of Euclidean projection, this is equivalent to $\a_* = \Ii_S(\a_*)$ for every $\eta > 0$.

For (c), continuity of the Hessian and its positive definiteness at $\a_*$ give $\rho > 0$ and $0 < c \leq C$ such that
\begin{align*}
    cI \preceq \nabla^2G_S(\a) \preceq CI
\end{align*}
for every $\a \in B_\rho(\a_*) \cap \Omega_\a$. Set $\tilde\eta := 1/C$. For $0 < \eta \leq \tilde\eta$ and $\a \in B_\rho(\a_*) \cap \Omega_\a$, the line segment joining $\a_*$ and $\a$ lies in $B_\rho(\a_*) \cap \Omega_\a$ by convexity. The fundamental theorem of calculus therefore gives
\begin{align*}
    \nabla G_S(\a) - \nabla G_S(\a_*) = M(\a - \a_*), \qquad M := \int_0^1 \nabla^2G_S\bigl(\a_* + t(\a - \a_*)\bigr)\,\d t,
\end{align*}
with $cI \preceq M \preceq CI$. Since Euclidean projection does not increase distances and $\a_* = \Ii_S(\a_*)$ by part (b), we get
\begin{align*}
    \begin{aligned}
        \|\Ii_S(\a) - \a_*\|_2 \leq \|(I - \eta M)(\a - \a_*)\|_2 \leq (1 - \eta c)\|\a - \a_*\|_2.
    \end{aligned}
\end{align*}
Moreover, $0 < \eta c \leq c/C \leq 1$, so $1 - \eta c \in [0, 1)$.
\end{proof}

\end{document}